\documentclass[MNA]{mna} % 

\usepackage[T1]{fontenc}
\usepackage[utf8]{inputenc}

\SHORTTITLE{Propagation of chaos in mean field networks of FHN neurons}

\TITLE{Propagation of chaos in mean field networks of FitzHugh-Nagumo neurons.} %\thanks{}} % \thanks is optional. Insert line breaks with \\

\AUTHORS{%
  Laetitia Colombani\orcidlink{0000-0001-8952-3761} \footnote{Institut Mathématiques de 
  Toulouse - Université de Toulouse 3, 
  France.
  		\BEMAIL{laetitia.colombani.math@gmail.com}}
  \and %% remove this line and below if single author
 Pierre Le Bris\footnote{Laboratoire Jacques-Louis Lions, Sorbonne Université, France.
 	\BEMAIL{pierre.lebris@sorbonne-universite.fr}}}%AUTHORS
\SHORTAUTHOR{L. Colombani P. Le Bris}

\KEYWORDS{FitzHugh-Nagumo; Propagation of Chaos; Mean field; Reflection Coupling} % 
\AMSSUBJ{92C20; 60H10; 60F99} % Edit. Separate items with ;
\SUBMITTED{June 28, 2022} % Edit.
\ACCEPTED{May 3, 2023} % Edit.

\ARXIVID{2206.13291} % Edit.
\HALID{hal-03706842} % Edit.

\VOLUME{3}
\YEAR{2023}
\PAPERNUM{3}
\DOI{10.46298/mna.9748}

\ABSTRACT{In this article, we are interested in the behavior of a fully connected network of 
$N$ neurons, where $N$ tends to infinity. We assume that the neurons follow the stochastic 
FitzHugh-Nagumo model, whose specificity is the non-linearity with a cubic term. We prove a 
result of uniform in time propagation of chaos of this model in a mean-field framework. We 
also exhibit explicit bounds. We use a coupling method initially suggested by  
Eberle~\cite{Eberle_reflection}, and recently extended in~\cite{DEGZ20}, known as the 
\textit{reflection coupling}. We simultaneously construct a solution of the $N$-particle system 
and $N$ independent copies of the non-linear McKean-Vlasov limit in such a way that, 
considering an appropriate semi-metric that takes into account the various possible behaviors 
of the processes, the two solutions tend to get closer together as $N$ increases, uniformly in 
time. The reflection coupling allows us to deal with the non-convexity of the underlying 
potential in the dynamics of the quantities defining our network, and show independence at 
the limit for the system in mean field interaction with sufficiently small Lipschitz continuous 
interactions.}

\DeclareMathOperator{\e}{e}

\begin{document}

%%%%%%%%%%%%%%%%%%%%%%%%%%%%%%%%%%%%%%%%%%%%%%%%%%%%%%%%%%%%%%%%%%%
%%                                                               %%
%% No need for \maketitle.                                       %%
%%                                                               %%
%%%%%%%%%%%%%%%%%%%%%%%%%%%%%%%%%%%%%%%%%%%%%%%%%%%%%%%%%%%%%%%%%%%

%%%%%%%%%%%%%%%%%%%%%%%%%%%%%%%%%%%%%%%%%%%%%%%%%%%%%%%%%%%%%%%%%%%
%%                                                               %%
%% Please replace what follows by the body of your article       %%
%% (up to the bibliography):                                     %%
%%                                                               %%
%%%%%%%%%%%%%%%%%%%%%%%%%%%%%%%%%%%%%%%%%%%%%%%%%%%%%%%%%%%%%%%%%%%

\section{Introduction}

%
%Subsection
%

\subsection{Understanding the model}

Understanding brain activity is both a complex and important challenge in current research. Of 
course, interests are plentiful: characterizing brain functions, unveiling structures and 
links between them, and understanding some phenomena  such as cyclic heartbeat. A way of 
modeling this activity is by considering a very large number of individual neurons with 
interactions. Since the number of neurons in a human brain is around $10^{11}$, and even 
``small'' parts of the brain are thus constituted of a very large number of them, such a 
strategy can be considered coherent. 

The main quantity we study is the membrane potential of the nerve cells: it can ``easily'' be 
observed and its modification characterizes a synapse (an interaction between neurons). 
Neurons regulate their electrical potential. In general, without interaction, the potential evolves 
with time but has quite small changes. Incoming potentials from other neurons are usually 
what make the neuron fire, \textit{i.e.} send action potentials to other neurons. We will here 
focus on a 
homogeneous network of neurons and consider mean-field interactions. This way, each 
neuron will interact with every other one, as it can be the case in small regions of the brain. 
The parameters of the model will be considered the same for each neuron.

A classical model was introduced by Hodgkin and Huxley~\cite{Hodgkin_Huxley_1952} using 
experimental data on the activity of the giant squid axon. It describes the ion exchanges 
$\text{K}^+$, $\text{Na}^+$, and $\text{Cl}^-$ through the membrane and their effects on the 
potential. A simplification of this model is the FitzHugh-Nagumo model, which reduces the 
dimension: from a four-dimensional model (for one neuron) with the Hodgkin-Huxley 
equations, 
we obtain a two-dimensional model, thus yielding a compromise between biological accuracy 
and mathematical simplicity.

The deterministic FitzHugh-Nagumo model for one neuron (or one particle) is given by the following equations
\begin{equation*}
	\left\{
	\begin{array}{ll}
		dX_{t}=(X_{t}-{(X_{t})}^3-C_{t}-\alpha)dt\\
		dC_{t}=(\gamma X_{t}-C_{t}+\beta)dt, 
	\end{array}
	\right.
\end{equation*}
where $X$ is the membrane potential and $C$ is a recovery variable, called the adaptation 
variable. The parameters $\gamma$ and $\beta$ are positive constants that determine the 
duration of excitation and the position of the equilibrium point of this system. Finally, $\alpha 
\in \mathbb{R}$ is the magnitude of a stimulus current (an entrance current in the system). 
Note that the variable $C$ isn't a physical quantity, and is used to allow $X$ to mimic the 
behavior of the potential. This variable $C$ has linear dynamics and provides slower negative 
feedback. 

%Thi13
This deterministic model has been largely studied. In Chapter 7 of~\cite{Thieullen_2013}, 
Thieullen describes the behavior of the solution of one deterministic FitzHugh-Nagumo 
system. She also extends the result in the case of a stochastic FitzHugh-Nagumo system: she 
considers a noise on the dynamics of $X$.

In fact, noise can be introduced in both equations to model different types of randomness: 
when the noise is on the first equation (dynamics of $X$) with a standard deviation 
$\sigma_{X} 
> 0$, it models a noisy presynaptic current. When it is on the second equation (dynamics of 
$C$) with a standard deviation $\sigma_{C} > 0$, it describes a noisy conductance dynamic (a 
noise in the chemical behavior). In general, noise in this model is additive. Various 
mathematical questions can be studied. Some authors choose to focus on the properties of 
the natural macroscopic limit of the model as $N \rightarrow \infty$ when it is clearly defined 
(see system~\eqref{eq:FN_limit}) while others work on properties of the particle system for 
fixed $N$. These models can be quite complicated to study mathematically. The main 
objectives are to characterize the behavior of these models when the number of neurons $N$ 
tends to $+ \infty$ in a mean-field limit, and to prove whether or not there exists an 
equilibrium, a stationary behavior, when $t$ tends to $+ \infty$. The question of the 
synchronization of neurons can also be studied, since it is a phenomenon observed in different 
contexts, such as the generation of respiratory rhythm or complex neurological functionalities. 
It can be characterized as the dissipation of the empirical variance of the system of neurons. 
We refer the reader to~\cite{bossySynchronizationStochasticMean2019} for further discussion 
on the synchronization in neuron models, and especially in the Hodgkin-Huxley model. 

% TRW03
In~\cite{Tuckwell_Rodriguez_Wan_2003}, the authors work on the determination of firing 
times. They consider a stochastic FitzHugh-Nagumo model for one neuron, with Brownian 
noise on $X$, obtain an approximation of firing times, and compare them with numerical 
simulations. 

% TBSST13
In~\cite{TatchimBemmo_SieweSiewe_tchawoua_2013}, Tatchim Bemmo, Siewe Siewe, and 
Tchawoua focus on a quite different stochastic model by considering additive noise $\eta$ on 
the dynamics of $X$, and multiplicative noise $\xi$ on the dynamics of $C$, both defined as 
sinusoidal functions of correlated Brownian motions. They choose to avoid Gaussian noises 
since they are unbounded. They also consider a deterministic and periodic entrance signal in 
the first equation, and observe abrupt transitions of the membrane potential $X$ when the 
intensity of the noise is gradually changed. 

In general, a lot of authors focus on noise on only one variable. 
% LS 18 
In~\cite{Leon_Samson_2018}, León and Samson consider a FitzHugh-Nagumo model with 
noise on $C$ but not on $X$, \textit{i.e.} $\sigma_{X} = 0$, and study the properties of the 
equations 
for one neuron. In particular, they focus on the hypoellipticity of the model, the existence and 
uniqueness of an invariant probability, and a mixing property by establishing a link between 
the model and the class of stochastic damping Hamiltonian systems. They also consider 
neuronal modeling questions and study the generation of spikes according to the parameters 
of the model.
% Uda19,
On the contrary, the article~\cite{Uda_2019} focuses on the stochastic FitzHugh-Nagumo 
model with noise in the dynamics of $X$, and $\sigma_{C}=0$. They study one neuron in a 
periodically forced regime. This study relies on the theory of Markovian Random Dynamical 
Systems. The model is driven by a cosine signal, and Uda studies the spike rate and 
compares it with the probability of a two-point motion of membrane potential. 

However, some do study stochastic models with two noises. 
	Berglund and Landon describe the behavior of the deterministic FitzHugh-Nagumo model 
	for one neuron in~\cite{Berglund_Landon_2012}, and consider the stochastic model, with 
	noise on both equations, to work on the behavior of the interspike interval and the 
	distribution of oscillations of the solution.

As said above, we consider mean-field interactions. These interactions are described by two 
functions $K_{X}$ and $K_{C}$, applied on the difference between two states $((X^{i}_{t}, 
C^{i}_{t}) - (X^{j}_{t}, C^{j}_{t}))$. 
In particular, this type of interaction models electrical synapses.

% BFFT12, MQT 16, BFO19, LP21
In their article~\cite{Baladron_Fasoli_Faugeras_touboul_2012}, Baladron, Fasoli, Faugeras, 
and Touboul study FitzHugh-Nagumo and Hodgkin-Huxley models with mean-field 
interaction, only on $X$. They consider more general interactions, not only applied to the 
difference between two states, modeling chemical synapses and electrical synapses. For the 
FitzHugh-Nagumo model, they consider a noise on $X$ and prove propagation of chaos, 
\textit{i.e.} 
the convergence of the law of $k$ neurons towards the law of $k$ independent solutions of 
the mean-field equations. This article is completed and clarified by the work of Bossy, 
Faugeras, and Talay in~\cite{bossyClarificationComplementMeanField2015}.
Mischler, Quininao, and Touboul consider a FitzHugh-Nagumo model 
in~\cite{Mischler_Quininao_touboul_2016}, with a linear interaction on $X$, and a noise only 
on 
$X$, \textit{i.e.} $\sigma_{C} = 0$ and $K_{X}(z) = \lambda x$. The drift on $X$ is not exactly 
the 
same 
as in the model above but remains similar as it is a cubic function of $X$. They work on the 
properties of a solution of the McKean-Vlasov PDE associated to this model and obtain the 
uniqueness of a global weak solution. Furthermore, they prove that there exists at least one 
stationary solution, and when the interaction is small, the stationary solution is unique and 
exponentially stable. They also exhibit numerical results with open problems, like attractive 
periodic solutions in time. In a similar framework, Luçon and Poquet study the macroscopic 
limit of this mean-field model in~\cite{Lucon_Poquet_2021}, and in particular the periodicity of 
such a system. They analyze the influence of both noise and interaction on the emergence of 
periodic behavior and prove the existence of a periodic solution, exponentially attractive, when 
the parameters satisfy some assumptions and the drift is ``small'' enough with respect to 
interaction and noise. Their approach relies on a slow-fast analysis and Floquet theory.
Results of non-uniform propagation of chaos has also been obtained 
in~\cite{mehriPropagationChaosStochastic2020} by Mehri et al., for stochastic spatially 
structured neuron networks, by applying the Euler 
approximation to the construction of a solution.

% BHV21
This model can be complexified, by considering non-mean-field interaction. In particular, 
Bayrak, Hövel, and Vuksanović work on a stochastic FitzHugh-Nagumo model with a network 
interaction in~\cite{Bayrak_Hovel_Vuksanovic_2021}. Their type of interaction takes into 
account a connectivity coefficient between two neurons and a propagation velocity. 
% Doublon sur LP21

% Spatial model: LW10, Li21, LL19, LX21
Other authors choose to complexify the model by considering stochastic FitzHugh-Nagumo 
with a spatial model. A second spatial derivative of $X$ is added to the dynamics of $X$. 
Various authors study the behavior of such a model and explore the notion of random 
attractors~\cite{Lv_Wang_2010, Li_2021, Li_Li_2019, Li_Xu_2021}.

% Numerical schemes: RS22
Various authors also study numerical schemes for the interacting particle system in the 
stochastic model. In~\cite{Reisinger_Stockinger_2022}, the authors adapt the 
Euler-Maruyama scheme to approximate the solution of the particle system.

%
%Subsection
%

\subsection{Framework and results}
Combining noise and interaction, we work specifically on the following equations, for $1 \leq i \leq N$, where $N$ is the number of neurons
\begin{equation}\label{eq:FN_MF}
	\left\{
	\begin{array}{ll}
		dX^{i,N}_{t}=(X^{i,N}_{t}-{(X^{i,N}_{t})}^3-C^{i,N}_{t}-\alpha)dt + \frac{1}{N} 
		\sum_{j=1}^{N} 
		K_{X}(Z^{i,N}_{t} - Z^{j,N}_{t})dt + \sigma_{X} dB^{i,X}_{t}\\
		dC^{i,N}_{t}=(\gamma X^{i,N}_{t}-C^{i,N}_{t}+\beta)dt + \frac{1}{N} \sum_{j=1}^{N} 
		K_{C}(Z^{i,N}_{t} - Z^{j,N}_{t})dt + \sigma_{C} dB^{i,C}_{t}, \\
	\end{array}
	\right.
\end{equation}
where we denote by $Z^{i,N}_{t}$ the couple $(X^{i,N}_{t}, C^{i,N}_{t})$ to simplify the 
notation. 

We assume ${(B^{i,X}_{t})}_i$ and ${(B^{i,C}_{t})}_i$ to be independent Brownian motions. 
Here, 
we consider two Brownian noises $B^X$ and $B^C$, one on each equation, and thus assume 
that each neuron has its own independent noise and that there is no environmental (or shared) 
noise.

We also assume $K_{X}$ and $K_{C}$ to be Lipschitz continuous and respectively denote 
their 
Lipschitz constants by $L_{X}$ and $L_{C}$.

The goal of this article is to describe the behavior of this network as the number $N$ of 
neurons tends to infinity.

To describe its behavior, we consider the $\mathbb{R}^2$-valued process 
${(\bar{Z}_{t})}_{t\geq 0} ={(\bar{X}_{t},\bar{C}_{t})}_{t\geq0}$ evolving according to the 
following 
non-linear stochastic differential equation of \textit{McKean-Vlasov} type
\begin{equation}\label{eq:FN_limit}
	\left\{
	\begin{array}{ll}
		d\bar{X}_{t}=(\bar{X}_{t}-{(\bar{X}_{t})}^3-\bar{C}_{t}-\alpha)dt+K_{X}\ast\bar{\mu}_{t}{(\bar{Z}_{t})}dt+\sigma_{X}
		 d\bar{B}^{X}_{t}\\
		d\bar{C}_{t}=(\gamma 
		\bar{X}_{t}-\bar{C}_{t}+\beta)dt+K_{C}\ast\bar{\mu}_{t}{(\bar{Z}_{t})}dt+\sigma_{C} 
		d\bar{B}^{C}_{t},
	\end{array}
	\right.
\end{equation}
where $\bar{\mu}_{t}=\text{Law}{(\bar{Z}_{t})}$ is the law at time $t$ of the process 
$(\bar{X}_{t},\bar{C}_{t})$, and $\ast$ denotes the operation of convolution, \textit{i.e.} 
\begin{align*}
	K_{X}\ast\bar{\mu}_{t}(u) = \int K_{X}(u-v) \bar{\mu}_{t}(dv).
\end{align*}
To some extent,~\eqref{eq:FN_MF} can be seen as an approximation of~\eqref{eq:FN_limit} in 
which the operation of convolution is applied to the empirical measure $\mu_{t, 
\text{emp}}=\frac{1}{N}\sum_{i=1}^{N}\delta_{Z^{i,N}_{t}}$, and what we wish to prove is that, 
indeed, the law $\mu^{N}_{t}$ of the network~\eqref{eq:FN_MF} converges in some sense to 
$\bar{\mu}^{\otimes N}_{t}$ (i.e the law of the solution of~\eqref{eq:FN_limit} tensorized $N$ 
times) as $N$ tends to infinity. This phenomenon has been stated under the name 
\textit{propagation of chaos} -an idea motivated by M. Kac~\cite{kac1956}- as it amounts to 
saying that, as the number of particles increases in the system, two particles will become 
``more and more'' independent, their joint law converging towards a tensorized law. The 
notion of 
``propagation'' refers to the fact that proving such convergence at time 0 is sufficient to 
prove it at a later time $t$.

In order to prove the convergence of $\mu^{N}_{t}$ to $\bar{\mu}_{t}^{\otimes N}$, we follow 
the 
coupling method described in recent work by one of the authors in~\cite{GLBM21_Coup}, 
the 
result of which cannot be applied directly here. This method has been put forward by Eberle, 
following earlier works by Lindvall and Rogers~\cite{Lindvall_Rogers_1986}. Before recalling 
the method, let us also mention the recent work~\cite{Schuh22}, which uses a coupling 
approach adapted to a well-chosen distance. 

We consider $r^i_{t}=|\bar{X}^i_{t}-X^{i,N}_{t}|+\delta|\bar{C}^i_{t}-C^{i,N}_{t}|$ with 
$\delta>0$, a constant not yet specified (to prove the first result Theorem~\ref{thm:non_unif}, 
we will consider $\delta = 1$, but we will need a more specific one for 
Theorem~\ref{thm:unif}). 

A natural distance between probability measures is the Wasserstein distance, linked to the theory of optimal transport. For $\mu$ and $\nu$ two probability measures on $\mathbb{R}^d$, we denote
\begin{equation}\label{eq:def_W}
	\mathcal{W}_p(\mu,\nu)=\inf_{X\sim\mu,\ Y\sim\nu}\mathbb{E}{\left(||X-Y||_p^p\right)}^{1/p},
\end{equation}
where $||\cdot||_p$ denotes the usual $L^p$ distance on $\mathbb{R}^d$. It is thus defined 
as the minimum over all possible choices of a pair $(X,Y)$, such that $X$ is distributed 
according to $\mu$ and $Y$ according to $\nu$, of the expectation of the distance between 
$X$ and $Y$. The basic idea behind a coupling method is then that an upper bound on the 
Wasserstein distance between $\mu$ and $\nu$ is given by the construction of any pair of 
random variables distributed according to these probability measures. Thus, instead of 
considering the minimum over all possible couplings, we construct simultaneously two 
solutions of~\eqref{eq:FN_MF} and~\eqref{eq:FN_limit} that will tend to get closer together as 
the number of neurons increases.

Let $\left(\bar{X}^i_{t},\bar{C}^i_{t}\right)$, for $i$ between $1$ and $N$, be $N$ 
independent 
copies of a solution of~\eqref{eq:FN_limit} driven by some independent Brownian motions 
${(\bar{B}^{i,X}_{t})}_{t\geqslant 0}$ and ${(\bar{B}^{i,C}_{t})}_{t\geqslant 0}$. A coupling of 
$\left(\bar{X}^i_{t},\bar{C}^i_{t}\right)$ and $\left(X^{i,N}_{t}, C^{i,N}_{t}\right)$ then follows 
from a 
coupling of the Brownian motions $B$ and $\bar{B}$. 

The first natural choice, popularized by Sznitman~\cite{Saint_Flour_1991}, is the 
\textit{synchronous} coupling and consists in choosing $B=\bar{B}$. By doing so, when 
considering the time evolution of 
${\bar{Z}^i_{t}-Z^{i,N}_{t}=\left(\bar{X}^i_{t}-X^{i,N}_{t},\bar{C}^i_{t}-C^{i,N}_{t}\right)}$, the 
noise 
cancels out. The contraction of a distance between the processes can then only be induced 
by the deterministic drift, as in~\cite{BGM10}, and this usually only holds under rather 
restrictive conditions (in particular the drift should be strongly convex). Nevertheless, in our 
case, the calculation of the evolution of $\bar{X}^i_{t}-X^{i,N}_{t}$ and 
$\bar{C}^i_{t}-C^{i,N}_{t}$ 
(see later) shows that there is still some deterministic contraction when 
$\bar{X}^i_{t}-X^{i,N}_{t}=0$. We can therefore use a synchronous coupling in the vicinity of 
this 
subspace.

Outside of this subspace, we use the noise to get the processes closer together. In the 
direction orthogonal to the contracting space we consider $B=-\bar{B}$, as this maximizes the 
variance of the noise. This yields the \textit{reflection} coupling. Notice however at this stage 
that, because of the symmetry of the noise, there is \textit{a priori} no reason why $r^i_{t}$ 
should 
decrease rather than increase. This invites us to consider $f(r^i_{t})$, with $f$ a concave 
function, so that a random decrease has more effect than a random increase of the same 
value. We will define the function $f$ later.

Finally, we construct a Lyapunov function $H$ to take into account the trend of each process 
to come back to some compact set of $\mathbb R^{2}$. We are then led to the study of a 
suitable distance between the two processes, which will be of the form 
$\rho_{t}:=\frac{1}{N}\sum_{i=1}^{N}f(r^i_{t})(1+\epsilon H(\bar{Z}^i_{t})+\epsilon 
H{(Z^{i,N}_{t})})$, 
where $\epsilon>0$. This quantity controls the usual $L^1$ and $L^2$ distances between the 
two systems and is interesting as, when $r^i_{t}$ is small, $f(r^i_{t})$ tends to decrease 
either 
because of the deterministic drift or the reflection coupling, and when $r^i_{t}$ is big, the 
Lyapunov functions $H$ will tend to decrease. We thus show that $\mathbb{E}\rho_{t}$ 
decays 
exponentially fast. This leads to several constraints on $\delta,\epsilon$ and on the 
parameters involved in the definition of $f$, and we have to prove that it is possible to meet all 
these conditions simultaneously. In reality, the quantity $\rho_{t}$ considered will be a slight 
twist of the one given above (see~\eqref{eq:def_rho}) so as to take into account the 
non-linearity of the process.

As explained, this method requires some noise in the direction orthogonal to the naturally 
contracting subspace. This means, in the description of the method above, that one should 
have $\sigma_{X}>0$ (so that we can use a reflection coupling to bring $\bar{X}^i_{t}$ and 
$X^{i,N}_{t}$ closer together). In the case $\sigma_{X}=0$ and $\sigma_{C}>0$, a 
modification 
of the calculations is necessary. We describe this case and the resulting modifications in the 
computations in Appendix~\ref{app:sigma_{X}_0}.

\begin{assumption}\label{hyp:K}
	$K_{X}$ and $K_{C}$ are Lipschitz continuous, \textit{i.e.}
	\begin{align*}
	\exists L_{X}\geq0,\forall z, z' \in\mathbb{R}^2\quad |K_{X}(z)-K_{X}(z')|\leq L_{X}(\|z-z'\|_1) 
	\\
	\exists L_{C}\geq0,\forall z, z'\in\mathbb{R}^2\quad |K_{C}(z)-K_{C}(z')|\leq L_{C}(\|z-z'\|_1) 
	\\
	K_{X}(0,0)=0\text{ and }K_{C}(0,0)=0.
	\end{align*}
\end{assumption}

Before any result on the propagation of chaos, we prove that both systems~\eqref{eq:FN_MF} 
and~\eqref{eq:FN_limit} have well-defined solutions. 
\begin{proposition}[Existence of solutions]\label{prop:existence_sol}
	Let $K_{X}$ and $K_{C}$ satisfy Assumption~\ref{hyp:K}. We assume the law of 
	$\left((X^{1,N}_0,C^{1,N}_0),\ldots,(X^{N,N}_0,C^{N,N}_0)\right)$ and the law of 
	$(\bar{X}_0, \bar{C}_0)$ have a moment of order 2.
	Then, there exists a unique strong solution for system~\eqref{eq:FN_MF} and a unique 
	strong solution for system~\eqref{eq:FN_limit}. 
\end{proposition}

We denote $\mathcal{W}_{1}$ and $\mathcal{W}_{2}$ the usual $L^1$ and $L^2$ 
Wasserstein distances defined in~\eqref{eq:def_W}.

%Theoreme
\begin{theorem}[Non uniform in time propagation of chaos]\label{thm:non_unif}
	Let $K_{X}$ and $K_{C}$ satisfy Assumption~\ref{hyp:K}. There exist explicit 
	${C_1,C_2>0}$, such that for all probability measures $\mu_0$ on $\mathbb{R}^{2}$ with 
	finite second moment, 
	\[
	\mathcal{W}_{1}\left(\mu^{k,N}_{t},\bar{\mu}_{t}^{\otimes k}\right)\leq C_1 
	e^{C_2t}\frac{k}{\sqrt{N}},
	\]
	for all $k\in\mathbb{N}$, where $\mu^{k,N}_{t}$ is the marginal distribution at time $t$ of 
	the 
	first $k$ neurons $\left((X^{1,N}_{t},C^{1,N}_{t}),\ldots,(X^{k,N}_{t},C^{k,N}_{t})\right)$ of 
	an 
	$N$-particle system~\eqref{eq:FN_MF} with initial distribution $\mu_0^{\otimes N}$, while 
	$\bar \mu_{t}$ is a solution of~\eqref{eq:FN_limit} with initial distribution $\mu_0$.
\end{theorem}
This first theorem is in accordance with the theorem 
from~\cite{Kumar_Neelima_Reisinger_Stockinger_2020} and makes the dependence in $t$ 
explicit. Since its proof is rather quick and provides a good entry point into coupling methods, 
we give it in Subsection~\ref{subsec:NonUniformTimeChaos}.

Our main result consists in removing the time dependency in the previous upper bound. This uniform in time propagation of chaos however requires stronger assumptions on the interaction kernels. 
%Theoreme
\begin{theorem}[Uniform in time propagation of chaos]\label{thm:unif}
	Let $L_{X,\max}$ and $L_{C,\max}$ be two (explicit) universal constants such that 
	$L_{X}\leq 
	L_{X,\max}$ and $L_{C}\leq L_{C,\max}$. Let $\mathcal{C}_{init,exp}>0$ and $\tilde{a}>0$. 
	There is an explicit 
	$c^K>0$ such that, for all $K_{X}$ and $K_{C}$ satisfying Assumptions~\ref{hyp:K} with 
	$L_{X},L_{C}<c^K$, there exist explicit ${B_1,B_2>0}$, such that for all probability 
	measures 
	$\mu_0$ on $\mathbb{R}^{2}$ satisfying 
	${\mathbb{E}_{\mu_0}\left(e^{\tilde{a}\left(|X|+|C|\right)}\right)\leq\mathcal{C}_{init,exp}}$, 
	\[
	\mathcal{W}_{1}\left(\mu^{k,N}_{t},\bar{\mu}_{t}^{\otimes k}\right)\leq B_1\frac{k }{\sqrt{N}} 
	\,,\qquad 
	\mathcal{W}_{2}^2\left(\mu^{k,N}_{t},\bar{\mu}_{t}^{\otimes k}\right)\leq B_2\frac{k 
	}{\sqrt{N}} ,
	\]
	for all $k\in\mathbb{N}$, where $\mu^{k,N}_{t}$ is the marginal distribution at time $t$ of 
	the 
	first $k$ neurons $\left((X^{1,N}_{t},C^{1,N}_{t}),\ldots,(X^{k,N}_{t},C^{k,N}_{t})\right)$ of 
	an 
	$N$-particle system~\eqref{eq:FN_MF} with initial distribution $\mu_0^{\otimes N}$, while 
	$\bar \mu_{t}$ is a solution of~\eqref{eq:FN_limit} with initial distribution $\mu_0$.
\end{theorem}
When we prove uniform in time propagation of chaos, $L_{X,\max}$ and $L_{C,\max}$ are 
\textit{a priori} bounds. Theorem~\ref{thm:unif} above will be true for $L_{X}$ and $L_{C}$ 
sufficiently small: the condition $L_{X}\leq L_{X,\max}$ and $L_{C}\leq L_{C,\max}$ are 
therefore 
not restrictive conditions and are useful in proving some parameters are independent of 
$L_{X}$ and $L_{C}$. Lemma~\ref{lem:Lya_limit} below shows that one can for instance 
consider 
$L_{X,\max}=4$ and $L_{C,\max}=\frac{1}{5}$. 
Furthermore $c^K$, that controls both interactions $K_{X}$ and $K_{C}$, is explained in 
Subsection~\ref{subsec:hyp_{C}onstants}. 

The main interest of obtaining uniform in time estimates is that it allows the study and 
comparison of the long-time behavior of the particle system and its nonlinear limit. As 
previously mentioned, this work follows the method described in~\cite{GLBM21_Coup}. 
Beyond the result of uniform in time propagation of chaos for the FitzHugh-Nagumo model, 
which is in itself an interesting result, the present work is also a testimony to the robustness of 
the coupling method. 

The reader will find \textbf{an index containing the notation}, constants, and parameters for reference at the end of the document.

\subsection{Existence of solutions}\label{subsec:Existence}
First of all, we prove Proposition~\ref{prop:existence_sol}, i.e existence of strong solutions of 
systems~\eqref{eq:FN_MF} and~\eqref{eq:FN_limit}, under Assumption~\ref{hyp:K}. Let's 
denote, for $\kappa \in \mathbb{R}^+$, 
\begin{align*}
	g_{\kappa}(x) = \left\{ 
	\begin{array}{ll}
		-\kappa^3 & \text{ if } x < -\kappa\\
		x^3 & \text{ if } x \in [-\kappa,\kappa]\\
		\kappa^3 & \text{ if } x > \kappa.
	\end{array}
	\right.
\end{align*}
$g_{\kappa}$ is Lipschitz and is bounded. 

Thus, it's well known (see Chapter 4~\cite{ikeda_89}) that the following system (under 
Assumption~\ref{hyp:K} and the assumption that the initial condition has a moment of order 2)
\begin{equation}
	\left\{
	\begin{array}{ll}
		dX^{i,N,\kappa}_{t}=&(X^{i,N,\kappa}_{t}-g_\kappa(X^{i,N,\kappa}_{t})-C^{i,N,\kappa}_{t}-\alpha)dt
		 + \frac{1}{N} \sum_{j=1}^{N} K_{X}(Z^{i,N,\kappa}_{t} - Z^{j,N,\kappa}_{t}) \\
		&+ \sigma_{X} dB^{i,X}_{t}\\
		dC^{i,N,\kappa}_{t}=&(\gamma X^{i,N,\kappa}_{t}-C^{i,N,\kappa}_{t}+\beta)dt + 
		\frac{1}{N} 
		\sum_{j=1}^{N} K_{C}(Z^{i,N,\kappa}_{t} - Z^{j,N,\kappa}_{t}) + \sigma_{C} dB^{i,C}_{t}, \\
	\end{array}
	\right.
\end{equation}
for $1 \leq i \leq N$, has a strong and unique solution that we denote 
${(X^{i,N,\kappa}_{t},C^{i,N,\kappa}_{t})}_{1 \leq i \leq N}$.

In consequence, for a fixed $\kappa \in \mathbb{R}^+$, there exists a strong solution of 
system~\eqref{eq:FN_MF} until time 
\[T_\kappa = \sup\left\{ t, \forall i, \forall s\leq t, 
X^{i,N,\kappa}_{s} \leq \kappa \text{ and } C^{i,N,\kappa}_{s} \leq \kappa \right\},
\]
and the 
solution coincides with the solution of the system with $g_\kappa$. 

We have the following Lemma
\begin{lemma}\label{lem:carreintegrable_{X}CK}
	If, for each $i \leq N$, $\mathbb{E}(|X^{i,N,\kappa}_0|^2) < +\infty$ and 
	$\mathbb{E}(|C^{i,N,\kappa}_0|^2) < +\infty$, then for all $t\geq0$ there exists 
	$\mathcal{C}_{t} < \infty$ such that, for each $i \leq N$
	\begin{equation}\label{eq:carreintegrable_{X}CK}
		\mathbb{E} \left(|X^{i,N,\kappa}_{t}|^2+|C^{i,N,\kappa}_{t}|^2\right) \leq \mathcal{C}_{t}. 
	\end{equation}
\end{lemma}
The proof relies on the Lyapunov function defined in the next Section and is given in 
Appendix~\ref{subsec:preuve_lya}.

Then, by denoting $T_{\infty}$ the explosion time of a solution of system~\eqref{eq:FN_MF}
\begin{align*}
	T_{\infty} = \inf \left\{ t, \exists i, \forall A>0, X^{i,N,\kappa}_{t}>A \text{ or } 
	C^{i,N,\kappa}_{t} > A \right\}
\end{align*}
we obtain that $\forall t \in \mathbb{R}^+, \mathbb{P}(T_{\infty}\leq t) = 0$ and 
$\mathbb{P}(\bar{T}_{\infty}\leq t) = 0$. Eventually, there exists a unique and strong solution 
of system~\eqref{eq:FN_MF}.

The existence and uniqueness of a solution of~\eqref{eq:FN_limit} is known from the Theorem 
3.3 from~\cite{Reis_Salkeld_tugaut_2019}, under the assumption that the law of the initial 
point $(\bar{X}_0, \bar{C}_0)$ has a moment of order 2. We only have to prove that 
Assumptions 3.2~\cite{Reis_Salkeld_tugaut_2019} are satisfied. We define, for all $t \in 
\mathbb{R}^+$, $z =(x,c) \in \mathbb{R}^2$ and for all probability distribution $\nu$ with a 
finite variance
\begin{align*}
	b(t, z, \nu) = \begin{pmatrix} 
		x - x^3 -c - \alpha + K_{X} \ast \nu(z)\\ 
		\gamma x - c + \beta + K_{C} \ast \nu (z)
	\end{pmatrix}\quad \text{ and }\quad 
	\sigma(t, z, \nu)= \begin{pmatrix} 
		\sigma_{X}\\ 
		\sigma_{C}
	\end{pmatrix}.
\end{align*}
$\sigma$ is a constant function, so it clearly satisfies the various conditions. 

For $t \in \mathbb{R}^+$, $z, z'$ in $\mathbb{R}^2$, and $\nu$ a probability measure
\begin{align*}
	\langle z-z', &~ b(t, z, \nu) - b(t, z', \nu) \rangle \\
	= & (x-x') \left[ (x -x') - (x^3-x'^3) -(c-c') + K_{X} \ast \nu(z) - K_{X} \ast \nu(z') \right]\\
	&+ (c-c') \left(\gamma (x-x') - (c-c') + K_{C} \ast \nu(z) - K_{C} \ast \nu(z')\right)\\
	=& {(x -x')}^2 - {(x -x')}^2 (x^2 + xx'+x'^2) + (\gamma -1)(c-c')(x-x') -{(c-c')}^2 \\
	&+ (x-x')(K_{X} \ast \nu(z) - K_{X} \ast \nu(z') ) + (c-c')(K_{C} \ast \nu(z) - K_{C} \ast 
	\nu(z')).
\end{align*}
Since $x^2 + xx'+x'^2 \geq 0$, the second term is non-positive. $K_{X}$ and $K_{C}$ are 
Lipschitz continuous  functions, so the last line is clearly bounded by $\|z-z'\|_2^2$ up to a 
multiplicative constant. Then, there exists a constant $L$ such that 
\[\langle z-z', b(t, z, \nu) - b(t, z', \nu) \rangle \leq L \|z-z'\|_2^2.\]
Since $K_{X}$ and $K_{C}$ are Lipschitz continuous functions, we also have, for all 
probability 
distributions $\nu$ and $\nu'$ with a finite variance, 
\[\|b(t, z, \nu) - b(t, z, \nu') \|_2 \leq L \mathcal{W}_2(\nu, \nu').\]
Eventually, since $b$ is locally Lipschitz continuous with polynomial growth, each Assumption 
is satisfied and Theorem 3.3~\cite{Reis_Salkeld_tugaut_2019} can be applied. Note that we 
could also apply Proposition 2.19 from~\cite{luconMeanFieldLimit2014}: assumptions are the 
same, and it gives a result for interaction depending on a spatial position. 

To complete the Lemma~\ref{lem:carreintegrable_{X}CK}, we also give the following
\begin{proposition}\label{prop:carreintegrable_{X}C}
	If, for each $i \leq N$, $\mathbb{E}(|X^{i,N}_0|^2) < +\infty$ and 
	$\mathbb{E}(|C^{i,N}_0|^2) < +\infty$, then for all $t\geq0$ there exists $\mathcal{C}_{t} < 
	\infty$ such that, for each $i \leq N$
	\begin{equation}\label{eq:carreintegrable_{X}C}
		\mathbb{E} \left(|X^{i,N}_{t}|^2+|C^{i,N}_{t}|^2\right) \leq \mathcal{C}_{t}. 
	\end{equation}
\end{proposition}
and
\begin{proposition}\label{prop:carreintegrable_barXC}
	If $\mathbb{E}(|\bar{X}_0|^2) < +\infty$ and $\mathbb{E}(|\bar{C}_0|^2) < +\infty$, then 
	there exist $C_{0,1}$ and $C_{0,2}$ such that
	\begin{equation}\label{eq:carreintegrable_barXC}
		\mathbb{E} \left(|\bar{X}_{t}|^2+|\bar{C}_{t}|^2\right) \leq C_{0,1}e^{C_{0,2}t}.
	\end{equation}
\end{proposition}
The proof is very similar to the proof of Lemma~\ref{lem:carreintegrable_{X}CK} and can be 
found in Appendix~\ref{subsec:preuve_lya}.

%
%Subsection
%

\subsection{Quick result: non uniform in time propagation of 
chaos}\label{subsec:NonUniformTimeChaos}
We start by proving Theorem~\ref{thm:non_unif}, a non uniform in time propagation of chaos, 
as it highlights the basic strategy behind a coupling argument. Some of the following 
expressions will be used in the proof of Theorem~\ref{thm:unif}, in 
Section~\ref{sec:proof_{t}hm_unif}.

We consider a synchronous coupling between ${(Z^{i,N}_{t})}_i$ and ${(\bar{Z}^i_{t})}_i$, 
\textit{i.e.} for each $1 \leq i \leq N$,  we choose $\tilde{B}^{i,X}_{t} = B^{i,X}_{t}$ and 
$\tilde{B}^{i,C}_{t} = B^{i,C}_{t}$. We have
\begin{equation*}
	\left\{
	\begin{array}{ll}
		dX^{i,N}_{t}=(X^{i,N}_{t}-{(X^{i,N}_{t})}^3-C^{i,N}_{t}-\alpha)dt+\frac{1}{N}\sum_{j=1}^{N}K_{X}(Z^{i,N}_{t}-Z^{j,N}_{t})dt+\sigma_{X}dB^{i,X}_{t}\\
		dC^{i,N}_{t}=(\gamma 
		X^{i,N}_{t}-C^{i,N}_{t}+\beta)dt+\frac{1}{N}\sum_{j=1}^{N}K_{C}(Z^{i,N}_{t}-Z^{j,N}_{t})dt+\sigma_{C}dB^{i,C}_{t}
	\end{array}
	\right.
\end{equation*}
and
\begin{equation*}
	\left\{
	\begin{array}{ll}
		d\bar{X}^i_{t}=(\bar{X}^i_{t}-{(\bar{X}^i_{t})}^3-\bar{C}^i_{t}-\alpha)dt+K_{X}\ast\bar{\mu}_{t}(\bar{Z}^i_{t})dt+\sigma_{X}
		 dB^{i,X}_{t}\\
		d\bar{C}^i_{t}=(\gamma 
		\bar{X}^i_{t}-\bar{C}^i_{t}+\beta)dt+K_{C}\ast\bar{\mu}_{t}(\bar{Z}^i_{t})dt+\sigma_{C} 
		dB^{i,C}_{t},
	\end{array}
	\right.
\end{equation*}
with $\bar{\mu}_{t}$ the law of $\bar{Z}^1_{t}$.
The method is the following
\begin{itemize}
	\item we compute the time evolution of 
	$\mathbb{E}\left(|X^{i,N}_{t}-\bar{X}^i_{t}|+|C^{i,N}_{t}-\bar{C}^i_{t}|\right)$ using Ito's 
	formula,
	\item we control the difference between the drifts $\frac{1}{N}\sum_{j\neq 
	i}K(\bar{Z}^i_{t}-\bar{Z}^j_{t})$ and $K\ast\bar{\mu}_{t}(\bar{Z}^i_{t})$ using some form of 
	the law of large numbers. This is where the convergence rate $\sqrt{N}$ appears,
	\item and we conclude using Gronwall's lemma.
\end{itemize}

\begin{description}
	\item[Time evolution:]
	We have,
	\begin{align*}
		d(X^{i,N}_{t}-\bar{X}^i_{t})= 
		&\Bigg((X^{i,N}_{t}-\bar{X}^i_{t})-\left({{(X^{i,N}_{t})}}^3-{(\bar{X}^i_{t})}^3\right)-(C^{i,N}_{t}-\bar{C}^i_{t})\\
		&+\frac{1}{N}\sum_{j=1}^{N}K_{X}(Z^{i,N}_{t}-Z^{j,N}_{t})-K_{X}\ast\bar{\mu}_{t}(\bar{Z}^i_{t})\Bigg)dt.
	\end{align*}
	We denote 
\[	
		\text{sign}(x)=
		\begin{cases}
			\dfrac{x}{|x|}\text{ if }x\neq0,\\ 0\text{ otherwise,}
		\end{cases}
\]
	and obtain, using Ito's formula for a twice continuously differentiable approximation of the absolute value and usual convergence lemmas (see Lemma~\ref{lem:ito_L1} below), 
	\begin{align}
		d&|X^{i,N}_{t}-\bar{X}^i_{t}|\nonumber\\
		=& \left(\text{sign}(X^{i,N}_{t}-\bar{X}^i_{t}) (X^{i,N}_{t}-\bar{X}^i_{t})- 
		\text{sign}(X^{i,N}_{t}-\bar{X}^i_{t}) 
		\left({{(X^{i,N}_{t})}}^3-{(\bar{X}^i_{t})}^3\right)\right.\nonumber \\
		&-\text{sign}(X^{i,N}_{t}-\bar{X}^i_{t})(C^{i,N}_{t}-\bar{C}^i_{t}) \nonumber\\
		&+\left. 
		\text{sign}(X^{i,N}_{t}-\bar{X}^i_{t})\frac{1}{N}\sum_{j=1}^{N}K_{X}(Z^{i,N}_{t}-Z^{j,N}_{t})-K_{X}\ast\bar{\mu}_{t}(\bar{Z}^i_{t})\right)dt\nonumber
		 \\
		\leq&\left(|X^{i,N}_{t}-\bar{X}^i_{t}|-\left|{{(X^{i,N}_{t})}}^3-{(\bar{X}^i_{t})}^3\right|+\left|C^{i,N}_{t}-\bar{C}^i_{t}\right|
		 \right. \nonumber \\ 
		& \left. 
		+\left|\frac{1}{N}\sum_{j=1}^{N}K_{X}(Z^{i,N}_{t}-Z^{j,N}_{t})-K_{X}\ast\bar{\mu}_{t}(\bar{Z}^i_{t})\right|\right)dt.
		 \label{eq:dX}
	\end{align}
	Similarly, 
	%\begin{align*}
	%d(C^{i,N}_{t}-\bar{C}^i_{t})=&\left(\gamma(X^{i,N}_{t}-\bar{X}^i_{t})-(C^{i,N}_{t}-\bar{C}^i_{t})+\frac{1}{N}\sum_{j=1}^{N}K_{C}(Z^{i,N}_{t}-Z^{j,N}_{t})-K_{C}\ast\bar{\mu}_{t}(\bar{Z}^i_{t})\right)dt,
	%\end{align*}
	%and we obtain
	\begin{align}
		d|&C^{i,N}_{t}-\bar{C}^i_{t}| \nonumber \\ 
		& 
		\leq\left(\gamma\left|X^{i,N}_{t}-\bar{X}^i_{t}\right|-|C^{i,N}_{t}-\bar{C}^i_{t}|+\left|\frac{1}{N}\sum_{j=1}^{N}K_{C}(Z^{i,N}_{t}-Z^{j,N}_{t})-K_{C}\ast\bar{\mu}_{t}(\bar{Z}^i_{t})\right|\right)dt.
		 \label{eq:dC}
	\end{align}
	Thus, denoting $r^i_{t}=|X^{i,N}_{t}-\bar{X}^i_{t}|+|C^{i,N}_{t}-\bar{C}^i_{t}|$ (i.e 
	considering $\delta = 1$) we obtain,
	\begin{align*}
		dr^i_{t}\leq&\left((1+\gamma)|X^{i,N}_{t}-\bar{X}^i_{t}|-\left|{{(X^{i,N}_{t})}}^3-{(\bar{X}^i_{t})}^3\right|\right.\\
		&+\left|\frac{1}{N}\sum_{j=1}^{N}K_{X}(Z^{i,N}_{t}-Z^{j,N}_{t})-K_{X}\ast\bar{\mu}_{t}(\bar{Z}^i_{t})\right|\\
		& \left. 
		+\left|\frac{1}{N}\sum_{j=1}^{N}K_{C}(Z^{i,N}_{t}-Z^{j,N}_{t})-K_{C}\ast\bar{\mu}_{t}(\bar{Z}^i_{t})\right|\right)dt.
	\end{align*}
	\item[Difference of the drifts:] Let us now consider these last two terms
	\begin{align*}
		&\left|\frac{1}{N}\sum_{j=1}^{N}K_{X}(Z^{i,N}_{t}-Z^{j,N}_{t})-K_{X}\ast\bar{\mu}_{t}(\bar{Z}^i_{t})\right|
		 \leq 
		\left|\frac{1}{N}\sum_{j=1}^{N}K_{X}(\bar{Z}^i_{t}-\bar{Z}^j_{t})-K_{X}\ast\bar{\mu}_{t}(\bar{Z}^i_{t})\right|
		 \\
		&\hspace{5cm}+ 
		\left|\frac{1}{N}\sum_{j=1}^{N}K_{X}(Z^{i,N}_{t}-Z^{j,N}_{t})-\frac{1}{N}\sum_{j=1}^{N}K_{X}(\bar{Z}^i_{t}-\bar{Z}^j_{t})\right|
		 . \\
	\end{align*}
	The first sum can be decomposed, using Assumption~\ref{hyp:K}, into
	\begin{align*}
		\frac{1}{N}\left|\sum_{j=1}^{N} 
		K_{X}(Z^{i,N}_{t}-Z^{j,N}_{t})-K_{X}(\bar{Z}^i_{t}-\bar{Z}^j_{t})\right| \leq& 
		\frac{L_{X}}{N}\sum_{j=1}^{N} \|Z^{i,N}_{t}-Z^{j,N}_{t}- (\bar{Z}^i_{t}-\bar{Z}^j_{t})\|_1\\
		%\leq & \frac{L_{X}}{N}\sum_{j=1}^{N} \left(\|Z^{i,N}_{t}-\bar{Z}^i_{t}\|_1 
		%+\|Z^{j,N}_{t}-\bar{Z}^j_{t}\|_1 \right)\\
		\leq & L_{X} r^i_{t} + \frac{L_{X}}{N}\sum_{j=1}^{N} r^j_{t}.
	\end{align*}
	Similarly, we obtain 
	\begin{align*}
		\left|\frac{1}{N}\sum_{j=1}^{N} K_{C}(Z^{i,N}_{t}-Z^{j,N}_{t})-K_{C} 
		\ast\bar{\mu}_{t}(\bar{Z}^i_{t})\right| \leq & L_{C} r^i_{t} + \frac{L_{C}}{N}\sum_{j=1}^{N} 
		r^j_{t} 
		\\
		+ & 
		\left|\frac{1}{N}\sum_{j=1}^{N}K_{C}(\bar{Z}^i_{t}-\bar{Z}^j_{t})-K_{C}\ast\bar{\mu}_{t}(\bar{Z}^i_{t})\right|.
	\end{align*}
	Hence, we get
	\begin{multline*}
		dr^i_{t}%\leq&\left((1+\gamma)|X^{i,N}_{t}-\bar{X}^i_{t}| 
		%-\left|{{(X^{i,N}_{t})}}^3-{(\bar{X}^i_{t})}^3\right| +(L_{X}+L_{C})\left(r^i_{t} + 
		%\frac{1}{N}\sum_{j=1}^{N} r^j_{t}\right)\right.\\
		%&\left.+\left|\frac{1}{N}\sum_{j=1}^{N}K_{X}(\bar{Z}^i_{t}-\bar{Z}^j_{t})-K_{X}\ast\bar{\mu}_{t}(\bar{Z}^i_{t})\right|+\left|\frac{1}{N}\sum_{j=1}^{N}K_{C}(\bar{Z}^i_{t}-\bar{Z}^j_{t})-K_{C}\ast\bar{\mu}_{t}(\bar{Z}^i_{t})\right|\right)dt\\
		\leq \left( (1+\gamma) r^i_{t} +(L_{X}+L_{C})\left(r^i_{t} + \frac{1}{N}\sum_{j=1}^{N} 
		r^j_{t}\right)\right.\\
	+\left|\frac{1}{N}\sum_{j=1}^{N}K_{X}(\bar{Z}^i_{t}-\bar{Z}^j_{t})-K_{X}\ast\bar{\mu}_{t}(\bar{Z}^i_{t})\right|\\
		+\left.\left|\frac{1}{N}\sum_{j=1}^{N}K_{C}(\bar{Z}^i_{t}-\bar{Z}^j_{t})-K_{C}\ast\bar{\mu}_{t}(\bar{Z}^i_{t})\right|\right)dt.
	\end{multline*} 
	By considering the expectation, since $\mathbb{E}(r^j_{t}) = \mathbb{E}(r^i_{t})$ for each 
	$j$, by exchangeability of the particles, we have
	\begin{align*}
		d \mathbb{E}(r^i_{t})\leq &\left( (1+\gamma + 2 L_{X}+ 2L_{C}) \mathbb{E}(r^i_{t}) + 
		\mathbb{E}\left[\left|\frac{1}{N}\sum_{j=1}^{N}K_{X}(\bar{Z}^i_{t}-\bar{Z}^j_{t})-K_{X}\ast\bar{\mu}_{t}(\bar{Z}^i_{t})\right|\right]\right.\\
		&\left.+ 
		\mathbb{E}\left[\left|\frac{1}{N}\sum_{j=1}^{N}K_{C}(\bar{Z}^i_{t}-\bar{Z}^j_{t})-K_{C}\ast\bar{\mu}_{t}(\bar{Z}^i_{t})\right|\right]\right)dt.
	\end{align*} 
	Now, we bound the interaction part. We begin with $K_{X}$. By Cauchy-Schwarz inequality, 
	we can write
	\begin{multline*}
		\mathbb{E} 
		\left[\left|\frac{1}{N}\sum_{j=1}^{N}K_{X}(\bar{Z}^i_{t}-\bar{Z}^j_{t})-K_{X}\ast\bar{\mu}_{t}(\bar{Z}^i_{t})\right|
		\right]\\
		\leq\mathbb{E}{\left(\left|\frac{1}{N}\sum_{j=1}^{N}K_{X}(\bar{Z}^i_{t}-\bar{Z}^j_{t})-K_{X}\ast\bar{\mu}_{t}(\bar{Z}^i_{t})\right|^2\right)}^{1/2}
	\end{multline*}
	We notice that ${(\bar{Z}^j_{t})}_j$ are i.i.d with law $\bar{\mu}_{t}$. Let's denote 
	$\bar{Z}_{t}$ a 
	generic random variable of law $\bar{\mu}_{t}$ independent of $\bar{Z}^i_{t}$. What is 
	more, 
	$K_{X}\ast\bar{\mu}_{t}(\bar{Z}^i_{t})= \int K_{X}(\bar{Z}^i_{t} - z) \bar{\mu}_{t}(dz) = 
	\mathbb{E}[K_{X}(\bar{Z}^i_{t} - \bar{Z}_{t})| \bar{Z}^i_{t}]$. Hence
	\begin{align*}
		\mathbb{E}&\left(\mathbb{E}\left(\left|\frac{1}{N-1}\sum_{j \neq i} 
		K_{X}(\bar{Z}^i_{t}-\bar{Z}^j_{t})-K_{X}\ast\bar{\mu}_{t}(\bar{Z}^i_{t})\right|^2\Big|\bar{Z}^i_{t}\right)\right)\\
		&\hspace{2cm}= \mathbb{E}\left( \text{Var} \left(\frac{1}{N-1}\sum_{j \neq i} 
		K_{X}(\bar{Z}^i_{t}-\bar{Z}^j_{t}) \Big|\bar{Z}^i_{t}\right) \right)\\
		%&\hspace{2cm}= \frac{1}{N-1} \mathbb{E}\left( \text{Var} \left( 
		%K_{X}(\bar{Z}^i_{t}-\bar{Z}_{t}) \Big|\bar{Z}^i_{t}\right) \right)\\
		&\hspace{2cm}\leq \frac{L_{X}^2}{N-1} \mathbb{E}\left( \text{Var} \left( \| 
		\bar{Z}^i_{t}-\bar{Z}_{t} \|_1 \Big|\bar{Z}^i_{t}\right) \right). 
	\end{align*}
	Since 
	\begin{align*}
		\mathbb{E} \left[ \text{Var} \left( \| \bar{Z}^i_{t}-\bar{Z}_{t} \|_1 \Big|\bar{Z}^i_{t}\right) 
		\right] \leq& \mathbb{E} \left[ \mathbb{E} \left( \| \bar{Z}^i_{t}-\bar{Z}_{t} \|_1^2 
		\Big|\bar{Z}^i_{t}\right) \right]\\
		\leq & \mathbb{E} \left[ \mathbb{E} \left( 2 \| \bar{Z}^i_{t}\|_1^2 + 2\|\bar{Z}_{t} \|_1^2 
		\Big|\bar{Z}^i_{t}\right)\right] 
		\leq 4 \mathbb{E}(\| \bar{Z}_{t}\|_1^2) , 
	\end{align*}
	we obtain 
	\begin{align*}
		\mathbb{E}\left(\mathbb{E}\left(\left|\frac{1}{N-1}\sum_{j \neq i} 
		K_{X}(\bar{Z}^i_{t}-\bar{Z}^j_{t})-K_{X}\ast\bar{\mu}_{t}(\bar{Z}^i_{t})\right|^2\Big|\bar{Z}^i_{t}\right)\right)
		&\leq \frac{4 L_{X}^2}{N-1} \mathbb{E}(\| \bar{Z}_{t}\|_1^2). 
	\end{align*}
	We now want to control 
	$\mathbb{E}\left(\left|\frac{1}{N}\sum_{j=1}^{N}K_{X}(\bar{Z}^i_{t}-\bar{Z}^j_{t})-K_{X}\ast\bar{\mu}_{t}(\bar{Z}^i_{t})\right|^2\right)$.
	 We decompose it into 
	\begin{align*}
		\mathbb{E}&\left(\left|\frac{1}{N}\sum_{j=1}^{N}K_{X}(\bar{Z}^i_{t}-\bar{Z}^j_{t})-K_{X}\ast\bar{\mu}_{t}(\bar{Z}^i_{t})\right|^2\right)\\
		=& 
		\mathbb{E}\left(\left|\frac{N-1}{N}\frac{1}{N-1}\sum_{j=1}^{N}K_{X}(\bar{Z}^i_{t}-\bar{Z}^j_{t})-\left(\frac{N-1}{N}
		 + \frac{1}{N}\right)K_{X}\ast\bar{\mu}_{t}(\bar{Z}^i_{t})\right|^2\right)\\
		\leq& 2 \frac{{(N-1)}^2}{N^2} 
		\mathbb{E}\left(\left|\frac{1}{N-1}\sum_{j=1}^{N}K_{X}(\bar{Z}^i_{t}-\bar{Z}^j_{t})- 
		K_{X}\ast\bar{\mu}_{t}(\bar{Z}^i_{t})\right|^2\right)\\
		&+ \frac{2}{N^2} 
		\mathbb{E}\left(|K_{X}\ast\bar{\mu}_{t}(\bar{Z}^i_{t})|^2\right).
	\end{align*}
	Since 
	\begin{align*}
		\mathbb{E}\left(|K_{X}\ast\bar{\mu}_{t}(\bar{Z}^i_{t})|^2\right) = & \mathbb{E}\left(\left| 
		\mathbb{E}\left( K_{X}(\bar{Z}^i_{t} - \bar{Z}_{t} )| \bar{Z^{i}_{t}} \right) \right|^2\right)\\
		\leq & L_{X}^2 \mathbb{E}\left(\mathbb{E}\left( \|\bar{Z}^i_{t} - \bar{Z}_{t} \|_1^2 | 
		\bar{Z^{i}_{t}} \right) \right)
		\leq 4 L_{X}^2 \mathbb{E}(\| \bar{Z}_{t}\|_1^2),
	\end{align*}
	we obtain
	\begin{align}
		\mathbb{E}&\left(\left|\frac{1}{N}\sum_{j=1}^{N}K_{X}(\bar{Z}^i_{t}-\bar{Z}^j_{t})-K_{X}\ast\bar{\mu}_{t}(\bar{Z}^i_{t})\right|^2\right)
		 \nonumber\\
		&\hspace{2cm}\leq {\left(\frac{N-1}{N}\right)}^2 \frac{4 L_{X}^2}{N-1} 
		\mathbb{E}(\|\bar{Z}_{t}\|_1^2) + \frac{4 L_{X}^2}{N^2} \mathbb{E}(\| \bar{Z}_{t}\|_1^2)
		\leq \frac{8 L_{X}^2}{N} \mathbb{E}(\| \bar{Z}_{t}\|_1^2), 
		\label{eq:controle_interactionmoyenne}
	\end{align}
	and finally 
	\begin{align*}
		\mathbb{E} 
		\left[\left|\frac{1}{N}\sum_{j=1}^{N}K_{X}(\bar{Z}^i_{t}-\bar{Z}^j_{t})-K_{X}\ast\bar{\mu}_{t}(\bar{Z}^i_{t})\right|
		 \right]\leq& {\left(\frac{8 L_{X}^2}{N} \mathbb{E}(\| \bar{Z}_{t}\|_1^2)\right)}^{1/2}.
	\end{align*}
	Similarly, we have 
	\begin{align*}
		\mathbb{E} 
		\left[\left|\frac{1}{N}\sum_{j=1}^{N}K_{C}(\bar{Z}^i_{t}-\bar{Z}^j_{t})-K_{C}\ast\bar{\mu}_{t}(\bar{Z}^i_{t})\right|
		 \right]\leq& {\left(\frac{8 L_{C}^2}{N} \mathbb{E}(\| \bar{Z}_{t}\|_1^2)\right)}^{1/2},
	\end{align*}
	which yields
	\begin{align*}
		d \mathbb{E}(r^i_{t})\leq &\left( (1+\gamma + 2 L_{X}+ 2L_{C}) \mathbb{E}(r^i_{t}) + 
		\sqrt{8 
		L_{X}^2 + 8 L_{C}^2} {\left(\frac{1}{N} \mathbb{E}(\| 
		\bar{Z}_{t}\|_1^2)\right)}^{1/2}\right)dt.
	\end{align*} 
	
	Then using Proposition~\ref{prop:carreintegrable_barXC}, we obtain 
	\begin{align*}
		d \mathbb{E}(r^i_{t})\leq &\left( (1+\gamma + 2 L_{X}+ 2L_{C}) \mathbb{E}(r^i_{t}) + 
		\frac{\sqrt{8 L_{X}^2 + 8 L_{C}^2} \sqrt{2C_{0,1}} 
		}{\sqrt{N}}e^{\frac{1}{2}C_{0,2}t}\right)dt
	\end{align*}
	\item[Conclusion:]
	We have thus obtained 
	\begin{align*}
		d&\left(\mathbb{E}(r^i_{t}) + \sqrt{\frac{16( L_{X}^2 + L_{C}^2) 
		C_{0,1}}{N}}\frac{1}{1+\gamma + 2 L_{X}+ 
		2L_{C}-\frac{C_{0,2}}{2}}e^{\frac{1}{2}C_{0,2}t} 
		\right)\\
		&\hspace{1cm}\leq (1+\gamma + 2 L_{X}+ 2L_{C}) \\
		&\hspace{1.5cm} \times  \left(\mathbb{E}(r^i_{t}) + \sqrt{\frac{16( L_{X}^2 + L_{C}^2) 
		C_{0,1}}{N}}\frac{1}{1+\gamma + 2 L_{X}+ 
		2L_{C}-\frac{C_{0,2}}{2}}e^{\frac{1}{2}C_{0,2}t} 
		\right)dt ,
	\end{align*}
	and Gronwall's lemma yields
	\begin{align*}
		\mathbb{E}(r^i_{t}) +& \sqrt{\frac{16( L_{X}^2 + L_{C}^2) C_{0,1}}{N}}\frac{1}{1+\gamma 
		+ 
		2 
		L_{X}+ 2L_{C}-\frac{C_{0,2}}{2}}e^{\frac{1}{2}C_{0,2}t}\\
		\leq & e^{(1+\gamma + 2L_{X}+2L_{C})t} \\
		&\times\left[\mathbb{E}(r^i_0)+\sqrt{\frac{16( L_{X}^2 + 
		L_{C}^2) C_{0,1}}{N}}\frac{1}{1+\gamma + 2 L_{X}+ 2L_{C}-\frac{C_{0,2}}{2}} \right],
	\end{align*}
	thus
	\begin{align*}
		\mathbb{E}(r^i_{t}) \leq C_1 e^{C_2 t} \frac{1}{\sqrt{N}}.
	\end{align*}
	Let $\mu_0$ a measure on $\mathbb{R}^2$, $\mu^{k,N}_{t}$ the marginal distribution at 
	time $t$ of the first $k$ neurons $\left(Z^{1,N}_{t}, \dots , Z^{k,N}_{t}\right)$ of an 
	$N$-particle system~\eqref{eq:FN_MF} with initial distribution $\mu_0^{\otimes N}$, and 
	$\bar \mu_{t}$ is a solution of~\eqref{eq:FN_limit} with initial distribution $\mu_0$.
	We obtain for the $L^1$ Wasserstein distance
	\begin{align*}
		\mathcal{W}_1(\mu^{k,N}_{t}, \bar{\mu}^{\otimes k}_{t}) = & \inf \left\lbrace \mathbb{E}[ \| 
		Z^{(k)}- \bar{Z}^{(k)} \|_1], \mathbb{P}_{Z^{(k)}} = \mu^{k,N}_{t}, 
		\mathbb{P}_{\bar{Z}^{(k)}} = \bar{\mu}^{\otimes k}_{t} \right\rbrace \\
		\leq& \inf \left\lbrace \mathbb{E}\left[ \sum_{i=1}^k r^i_{t}\right] , 
		\mathbb{P}_{{(Z^{i,N}_{t})}_i} = \mu^{k,N}_{t} , \mathbb{P}_{{(\bar{Z}^i_{t})}_i} = 
		\bar{\mu}^{\otimes k}_{t} \right\rbrace\\
		\leq & k \mathbb{E}(r^1_{t}) \\
		\leq & C_1 e^{C_2 t} \frac{k}{\sqrt{N}}.
	\end{align*}
	We hence obtain Theorem~\ref{thm:non_unif}.
\end{description}

%
%
%Section
%
%

\section{Preliminaries}

In this section, before tackling the proof by the coupling method of the uniform in time 
propagation of chaos, we gather the various technical lemmas and construct the necessary 
objects.

%
%
%Subsection
%
%

\subsection{Notation}

To construct the Lyapunov functions (which allow us to bound the moments of the processes 
and show that they tend to come back to some compact set), we begin by introducing the 
generators of the processes.

For $h: \mathbb{R}^{2N} \rightarrow \mathbb{R}$, for all ${(z_i)}_{1\leq i \leq 
N}={(x_i,c_i)}_{1\leq 
i \leq N} \in \mathbb{R}^{2N}$, the generator of~\eqref{eq:FN_MF} is
\begin{align*}
	\mathcal{L}^{N} h(z_1,\ldots,z_N) =&\sum_{i=1}^{N}\mathcal{L}^{i,N}h,
\end{align*}
where 
\begin{align*}
	\mathcal{L}^{i,N}h(z_1,\ldots,z_N) = &\left(x_i-x_i^3-c_i-\alpha+\frac{1}{N}\sum_{j=1}^{N} 
	K_{X}(z_i -z_j)\right) \partial_{x_i}h \\
	&+ \left(\gamma x_i-c_i+\beta+\frac{1}{N}\sum_{j=1}^{N} K_{C}(z_i -z_j)\right) 
	\partial_{c_i}h 
	\\
	&+\frac{\sigma_{X} ^2}{2} \partial^2_{x_i,x_i}h+\frac{\sigma_{C} ^2}{2} \partial^2_{c_i,c_i}h 
	.
\end{align*}
For $h: \mathbb{R}^{2} \rightarrow \mathbb{R}$, for all $z=(x,y) \in \mathbb{R}^2$, the 
generator of~\eqref{eq:FN_limit} for a given distribution $\mu$ is
\begin{align*}
	\mathcal{L}_{\mu}h(x,c)=&\left(x-x^3-c-\alpha+K_{X}\ast\mu(z)\right) 
	\partial_{x}h+\left(\gamma x-c+\beta+K_{C}\ast\mu(z)\right) \partial_{c}h\\
	&+\frac{\sigma_{X} ^2}{2}\partial^2_{xx}h+\frac{\sigma_{C} ^2}{2}\partial^2_{cc}h.
\end{align*}
In particular, we notice that for fixed ${(z_i)}_{1\leq i \leq N} \in {(\mathbb{R}^{2})}^{N}$, if we 
consider the empirical measure $\{\mu_{\text{emp}} = \frac{1}{N} \sum_j 
\delta_{z_j}\}$, we have for all $h: \mathbb{R}^{2} \rightarrow \mathbb{R}$ and $\bar{z} \in 
\mathbb{R}^2$,
\begin{align*}
	\mathcal{L}_{\mu_{\text{emp}}}h(\bar{z})=&\left(\bar{x}-\bar{x}^3-\bar{c}-\alpha+K_{X}\ast\mu_{\text{emp}}(\bar{z})\right)
	 \partial_{x} h +\left(\gamma \bar{x}- \bar{c}+\beta+K_{C}\ast\mu_{\text{emp}}(\bar{z})\right) 
	\partial_{c} h \\
	&+\frac{\sigma_{X} ^2}{2}\partial^2_{xx}h+\frac{\sigma_{C} ^2}{2}\partial^2_{cc}h\\
	=& \left(\bar{x}-\bar{x}^3-\bar{c}-\alpha+ \frac{1}{N} \sum_{j=1}^{N} K_{X}(\bar{z} 
	-z_j)\right) 
	\partial_{x} h\\
	& +\left(\gamma \bar{x}- \bar{c}+\beta+ \frac{1}{N} \sum_{j=1}^{N} K_{C}(\bar{z}- z_j) \right) 
	\partial_{c} h +\frac{\sigma_{X} ^2}{2}\partial^2_{xx}h+\frac{\sigma_{C} 
	^2}{2}\partial^2_{cc}h.
\end{align*}
In this case, if we consider $\bar{z}=z_i$ for a specific $i$ and we denote $\bar{h}^i: 
(z_1,\ldots, z_N) \rightarrow h(z_i)$, then
\begin{align*}
	\mathcal{L}_{\mu_{\text{emp}}} h(z_i) =& \left(x_i-x_i^3-c_i-\alpha+ \frac{1}{N} 
	\sum_{j=1}^{N} K_{X}(z_i-z_j)\right) \partial_{x} h \\
	&+\left(\gamma x_i- c_i+\beta+\frac{1}{N} \sum_{j=1}^{N} K_{C}(z_i -z_j) \right) \partial_{c} 
	h 
	+\frac{\sigma_{X} ^2}{2}\partial^2_{xx}h+\frac{\sigma_{C} ^2}{2}\partial^2_{cc}h\\
	&= \mathcal{L}^{i,N} \bar{h}^i (z_1, \dots, z_N).
\end{align*}

%
%
%Subsection
%
%

\subsection{First Lyapunov function}

Let $H: \mathbb{R}^2 \rightarrow \mathbb{R}$ be defined by
\begin{equation}\label{eq:def_H}
	H(z) = H(x,c)=\frac{1}{2}\gamma x^2+\beta x+\frac{1}{2}c^2+\alpha c+H_0,
\end{equation}
with
\begin{align*}
	H_0=\frac{\beta^2}{\gamma}+\alpha^2,
\end{align*}
where $\gamma$, $\beta$ and $\alpha$ are the parameters of the system~\eqref{eq:FN_MF}.

%Lemma

\begin{lemma}\label{prop:H}
	\begin{description}
		\item[(i)] For all $x,c\in\mathbb{R}$, we have $H(x,c)\geq\frac{\gamma}{4}x^2+\frac{c^2}{4}\geq0$,
		\item[(ii)] For all $x,c\in\mathbb{R}$, we have 
		$H(x,c)\geq\frac{1}{2\max(\gamma,1)}\left({(\gamma x+\beta)}^2+{(c+\alpha)}^2\right)$,
		\item[(iii)] For all $\delta>0$ there is $\mathcal{C}_{r,H}>0$ such that for all 
		$x,x',c,c'\in\mathbb{R}$, we have 
		\[{\left(|x-x'|+\delta|c-c'|\right)}^2\leq\mathcal{C}_{r,H}{(H(x,c)+H(x',c'))},\]
		\item[(iv)] A direct consequence of the previous point is that for all $B\in\mathbb{R}$, $\lambda>0$ and $\delta>0$, there is $R\geq0$ such that, for $x,x',c,c'\in\mathbb{R}$ satisfying $|x-x'|+\delta|c-c'|\geq R$, we have $H(x,c)+H(x',c')\geq\frac{80B}{\lambda}$. An explicit value of $R$ is given by $R=\sqrt{\frac{1280(1+\delta^2)B}{\lambda\min(\gamma,1)}}$.
	\end{description}
\end{lemma}
The first two points are consequences of direct calculations. The last two points are proved in 
Appendix~\ref{subsec:preuve_lemme_prop_H}. The constant $\mathcal{C}_{r,H}$  has been 
thus named because it ensures the control of the modified Euclidean distance $r$, precisely 
defined in~\eqref{eq:def_r}, by the function $H$.

\begin{lemma}[Lyapunov's property of $H$]\label{lem:Lya_limit}
	Let $\lambda\in\mathbb{R}$ such that
	\begin{equation}\label{eq:cond_lambda_lya}
		\frac{L_{X}}{8}+L_{C}\left(2+\frac{1}{8}\right)<1-\frac{\lambda}{2},
	\end{equation}
	then, for $H$ defined in~\eqref{eq:def_H}, there exists $B>0$ such that for all $(\bar{x}, 
	\bar{c}) \in \mathbb{R}^2$, for all probability distribution $\mu$ on $\mathbb{R}^2$, 
	\begin{align}\nonumber
		\mathcal{L}_{\mu}H(\bar{z}) \leq B&+\left(\alpha_{X} L_{X}+\beta_{X} L_{C}\right) 
		\left(\mathbb{E}_{\mu}{(|X|)}^2-\bar{x}^2\right)\\
		\label{eq:dyn_H}&+\left(\alpha_{C} L_{X}+ \beta_{C} 
		L_{C}\right)\left(\mathbb{E}_{\mu}{(|C|)}^2-\bar{c}^2\right)-\lambda H(\bar{z}).
	\end{align}
	Moreover, for all ${(z_i)}_{1\leq i \leq N} \in \mathbb{R}^{2N}$, by denoting $H:(z_1, \dots, 
	z_N) \mapsto H(z_i)$,
	\begin{align}
		\mathcal{L}^{i,N} H(z_1, \dots, z_N) \leq& B+\left(\alpha_{X} L_{X}+\beta_{X} L_{C}\right) 
		\left({\left(\frac{1}{N} \sum_{j=1}^{N} |x_j| \right)}^2 - x_i^2 \right) \nonumber \\
		&+ \left(\alpha_{C} L_{X}+ \beta_{C} L_{C}\right) 
		\left({\left(\frac{1}{N}\sum_{j=1}^{N}|c_j|\right)}^2 -c_i^2 \right) -\lambda 
		H(z_i)\label{eq:dyn_H_part},
	\end{align}
	with
	\begin{equation*}
		\alpha_{X}=\frac{\gamma}{2}+\frac{1}{2},\quad\beta_{X}=\frac{17}{2},\quad 
		\alpha_{C}=\frac{1}{16},\quad\beta_{C}=\frac{1}{2}+\frac{1}{32} . 
	\end{equation*}
	We refer to $H$ as a Lyapunov function, as it ensures that the processes tend to come back to a compact set.
\end{lemma}

We refer to Appendix~\ref{subsec:preuve_lya} for the proof of this lemma and of the following Proposition.

\begin{proposition}\label{prop:Lya}
	We have
	\begin{equation}\label{eq:gronwall_particles}
		\mathcal{L}^{N}\left(\frac{1}{N}\sum_{i=1}^{N}H\left(Z^{i,N}_{t}\right)\right)\leq 
		B-\lambda 
		\left(\frac{1}{N}\sum_{i=1}^{N}H\left(Z^{i,N}_{t}\right)\right),
	\end{equation}
\end{proposition}
A direct consequence of~\eqref{eq:dyn_H} is 
\begin{equation}\label{eq:gronwall_H}
	\mathbb{E}H\left(\bar{Z}^i_{t}\right)\leq \mathbb{E}H\left(\bar{Z}^i_0\right)+\int_0^{t}\left( 
	B-\lambda \mathbb{E}H\left(\bar{Z}^i_s\right)\right)ds,
\end{equation}
and a consequence of~\eqref{eq:gronwall_particles} is
\begin{equation}\label{eq:gronwall_H_MF}
	\left(\frac{1}{N}\sum_{i=1}^{N}\mathbb{E}H\left(Z^{i,N}_{t}\right)\right)\leq 
	\left(\frac{1}{N}\sum_{i=1}^{N}\mathbb{E}H\left(Z^{i,N}_0\right)\right)+\int_0^{t}\left(B-\lambda
	 \frac{1}{N}\sum_{i=1}^{N}\mathbb{E}H\left(Z^{i,N}_s\right)\right)ds.
\end{equation}
From~\eqref{eq:gronwall_H_MF} we obtain bounds on the moments of 
$\left|X^{i,N}_{t}\right|^2$ and $\left|C^{i,N}_{t}\right|^2$, and from~\eqref{eq:gronwall_H} 
Proposition~\ref{prop:carreintegrable_barXC} on the second moments of $\bar{X}^i_{t}$ and 
$\bar{C}^i_{t}$. The proof is given in Appendix~\ref{subsec:preuve_lya}. It also yields the 
following result
\begin{lemma}\label{lem:borne_unif_moment_2}
	Provided the interaction kernels satisfy~\eqref{eq:cond_lambda_lya},
	and that $\mathbb{E}(|\bar{X}_0|^2) < +\infty$ and $\mathbb{E}(|\bar{C}_0|^2) < +\infty$, then there exists $\mathcal{C}_{init,2}$ such that for all $t\geq0$
	\begin{equation*}
		\mathbb{E} \left(|\bar{X}_{t}|^2+|\bar{C}_{t}|^2\right) \leq \mathcal{C}_{init,2}.
	\end{equation*}
\end{lemma}

From now on, we consider $\lambda>0$ satisfying~\eqref{eq:cond_lambda_lya} (and use the 
\textit{a priori} bounds $L_{X,\max}$ and $L_{C,\max}$ to ensure the existence of such a 
$\lambda$).

%
%Subsection
%

\subsection{Modification of the function}\label{subsec:tildeH}
Let $\mathcal{C}_{init,exp}>0$, $\tilde{a}>0$ and consider an initial measure $\mu_0$ on $\mathbb{R}^{2}$ which satisfies ${\mathbb{E}_{\mu_0}\left(e^{\tilde{a}\left(|X|+|C|\right)}\right)\leq\mathcal{C}_{init,exp}}$, where $(X,C)$ is distributed according to $\mu_0$.

For technical reasons, we need a greater restoring force by the Lyapunov function than the 
one given in Lemma~\ref{lem:Lya_limit}. We thus modify it to obtain estimates such 
as~\eqref{eq:dyn_{t}ilde_H_non_lin} and~\eqref{eq:gen_part_sum} below.

Let $a>0$, such that $a \leq \tilde{a}/ \left(4\sqrt{2}\max{\left(\sqrt{\gamma}, 1\right)}\right)$. 
This choice of $a$ is only necessary for further Propositions and Lemmas, in 
Section~\ref{sec:proof_{t}hm_unif}. 

Let us consider for all $z \in \mathbb{R}^2$, 
\begin{equation} \label{eq:deftildeH}
	\tilde{H}(z)=\int_{0}^{H(z)}\exp\left(a\sqrt{u}\right)du=\frac{2}{a^2}\exp\left(a\sqrt{H(z)}\right)\left(a\sqrt{H(z)}-1\right)+\frac{2}{a^2}.
\end{equation}
Direct calculations yield the following technical lemma.
\begin{lemma}\label{lem:control_{t}ilde_H}
	We have, for all $z\in\mathbb{R}^2$
	\begin{align}
		H(z)\exp\left(a\sqrt{H(z)}\right)\geq \tilde{H}(z)\geq& \exp\left(a\sqrt{H(z)}\right) 
		-\frac{2}{a^2}\left(\exp\left(\frac{a^2}{2}\right)-1\right)\label{eq:control_{t}ilde_H},\\
		\frac{2}{a}\sqrt{H(z)}\exp\left(a\sqrt{H(z)}\right)\geq 
		\tilde{H}(z)\geq&\frac{1}{a}\sqrt{H(z)}\exp\left(a\sqrt{H(z)}\right)-\frac{1}{a^2}\left(e-2\right)\label{eq:control_{t}ilde_H_2},\\
		\tilde{H}(z)\geq& H(z) . \label{eq:control_{t}ilde_H_3}
	\end{align}
\end{lemma}
We may calculate, using Lemma~\ref{prop:H} and Equation~\eqref{eq:dyn_H}
\begin{align}
	\mathcal{L}_{\mu}\left(\tilde{H}\right)=&\exp\left(a\sqrt{H}\right)\mathcal{L}_{\mu} 
	H+\frac{1}{2}\frac{a}{2\sqrt{H}}\exp\left(a\sqrt{H}\right)\left(|\sigma_{X}\partial_{X}H|^2+|\sigma_{C}\partial_{C}H|^2\right)\nonumber\\
	%=&\exp\left(a\sqrt{H}\right)\mathcal{L}_{\mu} 
	%H+\frac{a}{4\sqrt{H}}\exp\left(a\sqrt{H}\right)\left(\sigma_{X}^2\left(\gamma 
	%x+\beta\right)^2+\sigma_{C}^2\left(c+\alpha\right)^2\right)\nonumber\\
	\leq& 
	\exp\left(a\sqrt{H}\right)\left(B+\left(\alpha_{X}L_{X}+\beta_{X}L_{C}\right)\mathbb{E}_{\mu}{(|X|)}^2\right.\\
	&\left.+\left(\alpha_{C}L_{X}+
	 \beta_{C}L_{C}\right)\mathbb{E}_{\mu}{(|C|)}^2-\lambda H\right)\nonumber\\
	&+\frac{1}{2}\max\left(\sigma_{X}^2,\sigma_{C}^2\right)\max\left(\gamma,1\right)a\sqrt{H}\exp\left(a\sqrt{H}\right)\nonumber\\
	\leq& 
	\exp\left(a\sqrt{H}\right)\left(B+\frac{{\left(\frac{1}{2}\max\left(\sigma_{X}^2,\sigma_{C}^2\right)\max\left(\gamma,1\right)\right)}^2a^2}{2\lambda}+\left(\alpha_{X}L_{X}+\beta_{X}L_{C}\right)\mathbb{E}_{\mu}{(|X|)}^2\right.\nonumber\\
	&\left.\hspace{2cm}+\left(\alpha_{C}L_{X}+ 
	\beta_{C}L_{C}\right)\mathbb{E}_{\mu}{(|C|)}^2-\frac{\lambda}{2} 
	H\right),\label{eq:dyn_exp_H}
\end{align}
where for this last inequality we used Young's inequality 
\[
\frac{1}{2}\max\left(\sigma_{X}^2,\sigma_{C}^2\right)\max\left(\gamma,1\right)a\sqrt{H}\leq 
\frac{\lambda}{2}H+\frac{{\left(\frac{1}{2}\max\left(\sigma_{X}^2,\sigma_{C}^2\right)\max\left(\gamma,1\right)\right)}^2a^2}{2\lambda}.
\]
Notice that~\eqref{eq:dyn_exp_H} ensures that this new Lyapunov function also tends to 
bring back particles into a compact set, and at an even greater rate. This new rate 
$H\exp(\sqrt{H})$ however comes at a cost: the initial condition must have a finite exponential 
moment, and no longer just have a finite second moment. First, by 
Lemma~\ref{lem:borne_unif_moment_2}, 
$\mathbb{E}{(\bar{X}_{t})}^2+\mathbb{E}{(\bar{C}_{t})}^2\leq\mathcal{C}_{init,2}$. 
Furthermore, 
the 
function $h\mapsto \exp\left(a\sqrt{h}\right)\left(B-\frac{\lambda}{4}h\right)$ is bounded from 
above for $h\geq0$. We therefore obtain from~\eqref{eq:dyn_exp_H} the existence of 
$\tilde{B}$ such that
\begin{align}
	\mathcal{L}_{\bar{\mu}_{t}}\left(\tilde{H}\left(\bar{Z}^i_{t}\right)\right)\leq& 
	\tilde{B}-\frac{\lambda}{4}\left(H\left(\bar{Z}^i_{t}\right)\exp\left(a\sqrt{H\left(\bar{Z}^i_{t}\right)}\right)\right)\label{eq:dyn_{t}ilde_H_non_lin}\\
	\frac{d}{dt}\mathbb{E} \tilde{H}\left(\bar{Z}^i_{t}\right)\leq& 
	\tilde{B}-\frac{\lambda}{4}\mathbb{E}\left(H\left(\bar{Z}^i_{t}\right)\exp\left(a\sqrt{H\left(\bar{Z}^i_{t}\right)}\right)\right)\label{eq:dyn_esp_exp_H}\\
	\text{and}\quad\frac{d}{dt}\mathbb{E} 
	\tilde{H}\left(\bar{Z}^i_{t}\right)\leq&\tilde{B}-\frac{\lambda }{4}\mathbb{E} 
	\tilde{H}\left(\bar{Z}^i_{t}\right),\label{eq:Gronwall_exp_H}
\end{align}
where for this last inequality, we used~\eqref{eq:control_{t}ilde_H}. 
While~\eqref{eq:dyn_{t}ilde_H_non_lin} and~\eqref{eq:dyn_esp_exp_H} will be useful in 
ensuring a sufficient restoring force, Equation~\eqref{eq:Gronwall_exp_H} gives us a uniform 
in 
time bound on $\mathbb{E} \tilde{H}\left(\bar{Z}^i_{t}\right)$, provided we have an initial 
bound. 
These inequalities are to be understood in the sense of SDEs, 
where~\eqref{eq:Gronwall_exp_H} should for instance be rigorously written
\begin{align*}
	\mathbb{E} \tilde{H}\left(\bar{Z}^i_{t}\right)\leq \mathbb{E} 
	\tilde{H}\left(\bar{Z}^i_0\right)+\int_0^{t} \left(\tilde{B}-\frac{\lambda }{4}\mathbb{E} 
	\tilde{H}\left(\bar{Z}^i_s\right)\right)ds.
\end{align*}
Now, for the system of particles, we have, using~\eqref{eq:dyn_exp_H}, $\forall i, \forall 
x_i,v_i\in\mathbb{R}^d$,
\begin{align*}
	\mathcal{L}^{N}\tilde{H}\left(z_i\right)\leq& 
	\exp\left(a\sqrt{H\left(z_i\right)}\right)\left(\tilde{B}+\left(\alpha_{X}L_{X}+\beta_{X}L_{C}\right){\left(\frac{1}{N}\sum_{j=1}^{N}|x_j|\right)}^2\right.\\
	&\left.\hspace{3cm}+\left(\alpha_{C}L_{X}+ 
	\beta_{C}L_{C}\right){\left(\frac{1}{N}\sum_{j=1}^{N}|c_j|\right)}^2-\frac{\lambda}{2} 
	H\left(z_i\right)\right).
\end{align*}
Summing over $i\in\{1,\ldots,N\}$, we may calculate
\begin{align*}
	\left(\alpha_{X}L_{X}+\beta_{X}L_{C}\right)&\sum_{j=1}^{N}{\left(\frac{\sum_{j=1}^{N}|x_j|}{N}\right)}^2\sum_{i=1}^{N}\frac{\exp\left(a\sqrt{H\left(z_i\right)}\right)}{N}-\frac{\lambda}{16}\sum_{i=1}^{N}\frac{H\left(z_i\right)\exp\left(a\sqrt{H\left(z_i\right)}\right)}{N}\nonumber\\
	\leq&\frac{\lambda}{16}\left(\sum_{i,j=1}^{N}\frac{H\left(z_i\right)}{N}\frac{\exp\left(a\sqrt{H\left(z_j\right)}\right)}{N}-\sum_{i=1}^{N}\frac{H\left(z_i\right)\exp\left(a\sqrt{H\left(z_i\right)}\right)}{N}\right)\nonumber\\
	\leq&0.%\label{annuler_somme_non_line}.
\end{align*}
Here, we used Lemma~\ref{prop:H}, the fact that $\forall x,y\geq0, 
xe^{\sqrt{y}}+ye^{\sqrt{x}}-xe^{\sqrt{x}}-ye^{\sqrt{y}}=(e^{\sqrt{x}}-e^{\sqrt{y}})(y-x)\leq0$ 
and assumed 
\[
\left(\alpha_{X}L_{X}+\beta_{X}L_{C}\right)\leq\frac{\gamma\lambda}{64}.
\]
Likewise,
\begin{multline*}
	\left(\alpha_{C}L_{X}+\beta_{C}L_{C}\right) 
	\sum_{j=1}^{N}{\left(\frac{\sum_{j=1}^{N}|c_j|}{N}\right)}^2\sum_{i=1}^{N}\frac{\exp\left(a\sqrt{H\left(z_i\right)}\right)}{N}\\
	-\frac{\lambda}{16}\sum_{i=1}^{N}\frac{H\left(z_i\right)\exp\left(a\sqrt{H\left(z_i\right)}\right)}{N}\leq0,%
	% \label{annuler_somme_non_line_{C}},
\end{multline*}
provided
\[
\left(\alpha_{C}L_{X}+\beta_{C}L_{C}\right)\leq\frac{\lambda}{64}.
\]
There is therefore a constant, which for the sake of clarity we will also denote $\tilde{B}$ (as we may take the maximum of the previous constants), such that we get
\begin{align}
	\mathcal{L}^{i,N}\tilde{H}{(Z^{i,N}_{t})}\leq&\tilde{B}+\left(\alpha_{X}L_{X}+\beta_{X}L_{C}\right){\left(\frac{\sum_{j=1}^{N}|X^{j,N}_{t}|}{N}\right)}^2\exp\left(a\sqrt{H\left(Z^{i,N}_{t}\right)}\right)\nonumber\\
	&+\left(\alpha_{C}L_{X}+\beta_{C}L_{C}\right){\left(\frac{\sum_{j=1}^{N}|C^{j,N}_{t}|}{N}\right)}^2\exp\left(a\sqrt{H\left(Z^{i,N}_{t}\right)}\right)\nonumber\\
	&-\frac{\lambda}{4} H\left(Z^{i,N}_{t}\right)\exp\left(a\sqrt{H\left(Z^{i,N}_{t}\right)}\right), 
	\label{eq:dyn_part_{t}ilde_H}
\end{align}
\begin{align}
	\mathcal{L}^{N}\left(\frac{1}{N}\sum_{i=1}^{N}\tilde{H}{(Z^{i,N}_{t})}\right)\leq& 
	\tilde{B}-\frac{\lambda}{4}\left(\frac{1}{N}\sum_{i=1}^{N}H{(Z^{i,N}_{t})}\exp\left(a\sqrt{H\left(Z^{i,N}_{t}\right)}\right)\right)
	 , \label{eq:gen_part_sum}
\end{align}
and
\begin{align}
	\mathcal{L}^{N}&\left(\frac{1}{N}\sum_{i=1}^{N}\tilde{H}{(Z^{i,N}_{t})}\right)\leq 
	\tilde{B}-\frac{\lambda}{4}\left(\frac{1}{N}\sum_{i=1}^{N}\tilde{H}{(Z^{i,N}_{t})}\right) . 
	\label{eq:gron_part_exp}
\end{align}
Once again,~\eqref{eq:dyn_part_{t}ilde_H} and~\eqref{eq:gen_part_sum} will be useful in 
ensuring a sufficient restoring force, and~\eqref{eq:gron_part_exp} yields a uniform in time 
bound on the expectation of $\tilde{H}{(Z^{i,N}_{t})}$, since by exchangeability of the 
particles,  
$\mathbb{E}\left(\frac{1}{N}\sum_{j=1}^{N}\tilde{H}(Z^{j,N}_{t})\right) = 
\mathbb{E}\left(\tilde{H}{(Z^{i,N}_{t})}\right)$.

%
%Subsection
%

\subsection{Parameters}\label{subsec:hyp_{C}onstants}

We start by fixing the values of some parameters. The somewhat intricate expressions in this section are dictated by the computations arising in the proofs later on. They are somewhat roughly chosen and far from optimal as we only wish to convey the fact that every constant is explicit. On first reading, the exact choice of parameters can and should be skipped, as they are only meant to satisfy Lemma~\ref{lem:parameters}, which is the crucial Lemma of this subsection.

Recall $\alpha_{X}$, $\beta_{X}$, $\alpha_{C}$ and $\beta_{C}$ given in 
Lemma~\ref{lem:Lya_limit}. $a > 0$ is fixed from the last Subsection and the definition of 
$\tilde{H}$, and $\lambda$ and $\tilde{B}$ are obtained from the same Subsection.

Given any $\eta>4$ and $\tilde{\delta}>0$, consider the following set of parameters
\begin{align*}
	\delta=&(1+\tilde{\delta})\frac{1+L_{X,\max}}{1-L_{C,\max}},\ 
	R_0=\sqrt{\frac{1280\tilde{B}}{\lambda\min(\gamma,1)}},\ R=\sqrt{1+\delta^2}R_0,\\
	\mathcal{C}_{f,1}=&16\left(\left(\gamma+a\left(\beta+\frac{\alpha}{\delta}\right)\sqrt{2\max\left(\gamma,1\right)}\right)\frac{\e^{a^2/2}-1}{a^2}+\sqrt{2\max\left(\gamma,1\right)}\left(\sqrt{\gamma}+\frac{1}{\delta}\right)(\e-2)\right) , \\
	\mathcal{C}_{f,2}=&4\left(\gamma+\left(a\left(\beta+\frac{\alpha}{\delta}\right)+2a^2\left(\sqrt{\gamma}+\frac{1}{\delta}\right)\right)\sqrt{2\max\left(\gamma,1\right)}\right) , \\
	c=&\min\left\{\frac{2\tilde{B}}{\eta}, \frac{\lambda}{160}\frac{\eta-4}{\eta},\right.\\
	&\left.\hspace{0.2cm} \frac{\min\left(\frac{\sigma_{X}}{\sqrt{\pi}R},1-L_{C, 
	\max}-\frac{1+L_{X, \max}}{\delta}\right)}{2(1+\eta)}\right.\\
	&\left.\hspace{2cm} 
	\times\exp\left(-\frac{1}{4\sigma_{X}^2}\left(1+\delta\gamma+L_{X,\max}+ \delta 
	L_{C,\max}+\left(\mathcal{C}_{f,1}+\mathcal{C}_{f,2}\right)\sigma_{X}^2\right)R^2\right)\right\},\\
	\epsilon=&\frac{\eta c}{2\tilde{B}}, 
	\phi_{\min}=\exp\left(-\frac{1}{4\sigma_{X}^2}(1+\delta\gamma+L_{X,\max}+\delta 
	L_{C,\max}+(\epsilon\mathcal{C}_{f,1}+\mathcal{C}_{f,2})\sigma_{X}^2)R^2\right) , \\
	\mathcal{C}_1=&\frac{1}{\min\left(\delta, 
	1\right)}\frac{2}{\phi_{\min}}\max\left(\frac{16(1+\delta^2)}{\epsilon\min\left(\gamma,1\right)},1\right),
	\mathcal{C}_2=\frac{1}{\min\left(\delta^2, 
	1\right)}\frac{2}{\phi_{\min}}\max\left(\frac{16(1+\delta^2)}{\epsilon\min\left(\gamma,1\right)},1\right)
	 , \\
	\mathcal{C}_{z}=&\frac{2}{\phi_{\min}}\max\left(1,\frac{4}{\epsilon}\max\left(\sqrt{\frac{1}{\gamma}},1\right)\right).
\end{align*}
We define $f$ as follows
\begin{align}
	f(r)=&\int_{0}^{r\wedge R}\phi(s)g(s)ds, \label{eq:def_f} \\
	\phi(r)=&\exp\left(-\frac{1}{4\sigma_{X}^2}\left(1+\delta\gamma+L_{X}+ \delta 
	L_{C}+\left(\epsilon \mathcal{C}_{f,1}+\mathcal{C}_{f,2}\right)\sigma_{X}^2\right)r^2\right), 
	\nonumber 
	\end{align}
\begin{align*}
	\Phi(s)=&\int_0^s\phi(u)du, \nonumber \\
	g(r)=&1-\frac{c+2\epsilon \tilde{B}}{\sigma_{X}^2}\int_0^r\Phi(s){\phi(s)}^{-1}ds. \nonumber
\end{align*}
Assume furthermore that $L_{X}$ and $L_{C}$, the Lipschitz constants, satisfy
\begin{align}\label{eq:conditions_LXLC1}
	L_{X}\leq\min\left(\frac{\lambda}{128\mathcal{C}_{z}},
	\frac{\lambda a}{512\epsilon\mathcal{C}_{z}}, \frac{c}{2\mathcal{C}_1}\right)
	\quad\text{ and }&\quad L_{C}\leq\min\left(\frac{\lambda}{128\delta\mathcal{C}_{z}},
	\frac{\lambda a}{512\epsilon\delta\mathcal{C}_{z}}, \frac{c}{2\delta\mathcal{C}_1}\right) , \\
	\alpha_{X}L_{X}+\beta_{X}L_{C}\leq \frac{\gamma\lambda}{128}\quad\text{ and }
	&\quad \alpha_{C}L_{X}+ \beta_{C}L_{C}\leq \frac{\lambda}{128}, 
	\label{eq:conditions_LXLC2}\\
	\frac{L_{X}}{8}+L_{C}\left(2+\frac{1}{8}\right)<1-\frac{\lambda}{2}. 
	\label{eq:conditions_LXLC3}
\end{align}
Notice how the bounds on $L_{X}$ and $L_{C}$ depend on $c$. This is one of the reasons 
why 
we use the \textit{a priori} bounds $L_{X}\in[0,L_{X,\max}]$ and $L_{C}\in[0,L_{C,\max}]$ 
given 
in 
the assumptions of Theorem~\ref{thm:unif}: they allow us to bound $c$ and $\delta$ 
independently of $L_{C}$ and $L_{X}$. We are thus able to begin by choosing acceptable 
values 
for those parameters, before then giving upper bounds on $L_{X}$ and $L_{C}$. The 
condition 
of 
taking $L_{X}$ and $L_{C}$ small enough (the condition $L_{X} <c^K$ and $L_{C}<c^K$ for 
a 
well chosen
$c^K$, given in Theorem~\ref{thm:unif}) is necessary to satisfy the conditions 
of~\eqref{eq:conditions_LXLC1},~\eqref{eq:conditions_LXLC2} 
and~\eqref{eq:conditions_LXLC3}.

We quickly mention that the constants $\mathcal{C}_1$, $\mathcal{C}_2$ and 
$\mathcal{C}_{z}$ above come from Lemma~\ref{lem:rho_1_2} later. We gather some 
properties required in the calculations of the proof of Theorem~\ref{thm:unif} in the following 
lemma. Again, these properties are the ones motivating the choice of parameters

\begin{lemma}\label{lem:parameters}
	The set of parameters given in Subsection~\ref{subsec:hyp_{C}onstants} satisfy
	\begin{itemize}
		\item $f$ is $\mathcal{C}^2$ on $(0,R)$ such that $f'_+\left(0\right)=1$ and $f'_-\left(R\right)>0$, and constant on $[R,\infty)$. Moreover, $f$ is non-negative, non-decreasing, and concave, and for all $s\geq 0$,
		\begin{equation*}
			\min\left(s,R\right)f'_-\left(R\right)\leq f\left(s\right)\leq\min\left(s,f\left(R\right)\right)\leq\min\left(s,R\right).
		\end{equation*} 
		\item For all $r\in[0,R]$, $\phi(r)\geq\phi_{\min}$ and $g(r)\geq\frac{1}{2}$.
		\item We have the conditions
		\begin{align*}\frac
			2 f'(R) \geq &\exp\left(-\frac{1}{4\sigma_{X}^2}\left(1+\delta\gamma+L_{X}+ \delta 
			L_{C}+\left(\epsilon 
			\mathcal{C}_{f,1}+\mathcal{C}_{f,2}\right)\sigma_{X}^2\right)R^2\right),\\
			2c+4\epsilon \tilde{B}\leq 
			&\left(1-L_{C}-\frac{1+L_{X}}{\delta}\right)\min_{r\in(0,R]}\frac{f'(r)r}{f(r)},\\
			c\leq & \frac{\lambda}{160}\frac{\frac{80\epsilon \tilde{B}}{\lambda}}{1+\frac{80\epsilon 
			\tilde{B}}{\lambda}},\quad
			\frac{1+L_{X}}{1-L_{C}}<\delta\quad\text{ and }\quad\epsilon\leq 1.
		\end{align*}
	\end{itemize}
\end{lemma}

The proof of this lemma is done in Appendix~\ref{subsec:choix_parametres_section}.

%
%Subsection
%

\subsection{Control of the usual distances}
As explained previously, we consider a modified semi-metric. For $z=(x,c) \in \mathbb{R}^2$ and $z'=(x',c') \in \mathbb{R}^2$, define
\begin{equation}\label{eq:def_r}
	r(z,z')= r(x,c,x',c')=|x-x'|+\delta|c-c'|,
\end{equation}
where $\delta$ is given in Subsection~\ref{subsec:hyp_{C}onstants}, and let 
$\rho({(z_j, z'_j)}_{1\leq j \leq N})$ be defined as follows 
\begin{align}
	\label{eq:def_rho}
	\rho\left({(z_j,z'_j)}_{1\leq j \leq N} \right) &=\frac{1}{N}\sum_{i=1}^{N} 
	f\left(r\left(z_i,z'_i\right)\right) G^i\left({(z_j,z'_j)}_j\right),
\end{align}
where for each $i\in\{1,\ldots,N\}$, 
\begin{align}\label{eq:def_G}
	G^i\left({(z_j,z'_j)}_j\right)&=1+\epsilon \tilde{H}\left(z_i\right)+\epsilon 
	\tilde{H}\left(z'_i\right)+\frac{\epsilon}{N}\sum_{j=1}^{N}\tilde{H}\left(z_j\right)+\frac{\epsilon}{N}\sum_{j=1}^{N}\tilde{H}\left(z'_j\right).
\end{align}
An immediate corollary of the definition and properties of $H$ is that $\rho$ is a quantity on $\mathbb{R}^{4N}$ which controls the usual $L^1$ and $L^2$ distances.

%Lemma

\begin{lemma}\label{lem:rho_1_2}
	The constants $\mathcal{C}_1,\mathcal{C}_2, \mathcal{C}_{z}>0$, given in 
	Subsection~\ref{subsec:hyp_{C}onstants}, are such that for all $z=(x,c) \in\mathbb{R}^{2}$ 
	and $z'=(x',c') \in\mathbb{R}^{2}$
	\begin{itemize}
		\item[(i)] $\|z-z'\|_1 \leq \mathcal{C}_1f\left(r\left(z, z'\right)\right)\left(1+\epsilon 
		\tilde{H}(z)+\epsilon \tilde{H}(z')\right)$, 
		\item[(ii)] $\|z-z'\|_2^2 \leq \mathcal{C}_2f\left(r\left(z,z'\right)\right)\left(1+\epsilon 
		\tilde{H}(z)+\epsilon \tilde{H}(z')\right)$,
		\item[(iii)] $\|z-z'\|_1 \leq\mathcal{C}_{z} f(r(z, z'))\left(1+\epsilon \sqrt{H(z)}+\epsilon 
		\sqrt{H(z')}\right)$. 
	\end{itemize}
\end{lemma}

The proof of this lemma is postponed to 
Appendix~\ref{subsec:preuve_lemme_{C}ontrole_distance}.

%
%
%Section
%
%

\section{Proof of Theorem~\ref{thm:unif} in the case 
$\sigma_{X}>0$}\label{sec:proof_{t}hm_unif}

Let $\xi>0$ be a parameter destined for vanishing, and let 
$\varphi_{\text{sc}}:\mathbb{R}^+\mapsto\mathbb{R}^+$ and 
$\varphi_{\text{rc}}:\mathbb{R}^+\mapsto\mathbb{R}^+$ be two Lipschitz continuous 
functions such that
\begin{align}
	\forall x,\quad \varphi_{\text{sc}}^2(x)+\varphi_{\text{rc}}^2(x)=&1 ,\label{eq:Levy} \\
	\varphi_{\text{rc}}(x)=&1\text{ if }\xi\leq x\leq R ,\nonumber \\
	\varphi_{\text{rc}}(x)=&0\text{ if }x\leq \frac{\xi}{2}\text{ or }x\geq R+\xi.\nonumber
\end{align}
Intuitively, $\varphi_{\text{rc}}$ represents the region of space in which we consider a reflection coupling, and $\varphi_{\text{sc}}$ the one in which we consider a synchronous coupling. In reality, we would like to consider $\varphi_{\text{sc}}$ and $\varphi_{\text{rc}}$ indicator functions of the regions of space. However, we need to consider a Lipschitz approximation of indicator functions to ensure continuity (to apply Itô's calculus) and the strong existence and uniqueness of the stochastic processes. We thus simultaneously construct the following solutions, for $1 \leq i \leq N$
\begin{equation*}
	\left\{
	\begin{array}{ll}
		dX^{i,N}_{t}=&(X^{i,N}_{t}-{{(X^{i,N}_{t})}}^3-C^{i,N}_{t}-\alpha)dt+\frac{1}{N}\sum_{j=1}^{N}K_{X}(Z^{i,N}_{t}-Z^{j,N}_{t})dt\\
		&+\sigma_{X} 
		\varphi_{\text{sc}}\left(|X^{i,N}_{t}-\bar{X}^i_{t}|\right)dB^{i,sc,X}_{t}+\sigma_{X} 
		\varphi_{\text{rc}}\left(|X^{i,N}_{t}-\bar{X}^i_{t}|\right)dB^{i,rc,X}_{t}\\
		dC^{i,N}_{t}=&(\gamma 
		X^{i,N}_{t}-C^{i,N}_{t}+\beta)dt+\frac{1}{N}\sum_{j=1}^{N}K_{C}(Z^{i,N}_{t}-Z^{j,N}_{t})dt+\sigma_{C}
		 dB^{i,C}_{t},
	\end{array}
	\right.
\end{equation*}
and
\begin{equation*}
	\left\{
	\begin{array}{ll}
		d\bar{X}^i_{t}=&(\bar{X}^i_{t}-{(\bar{X}^i_{t})}^3-\bar{C}^i_{t}-\alpha)dt+K_{X}\ast\bar{\mu}_{t}(\bar{Z}^i_{t})dt\\
		&+\sigma_{X} 
		\varphi_{\text{sc}}\left(|X^{i,N}_{t}-\bar{X}^i_{t}|\right)dB^{i,sc,X}_{t}-\sigma_{X} 
		\varphi_{\text{rc}}\left(|X^{i,N}_{t}-\bar{X}^i_{t}|\right)dB^{i,rc,X}_{t}\\
		d\bar{C}^i_{t}=&(\gamma 
		\bar{X}^i_{t}-\bar{C}^i_{t}+\beta)dt+K_{C}\ast\bar{\mu}_{t}(\bar{Z}^i_{t})dt+\sigma_{C}dB^{i,C}_{t},
	\end{array}
	\right.
\end{equation*}
where ${(B^{i,sc,X})}_i$ and ${(B^{i,rc,X})}_i$ are independent Brownian motions (also 
independent of ${(B^{i,C})}_i$). 
Notice that we consider a symmetric coupling on the dynamics of $C$. By Levy's 
characterization of Brownian motion, using~\eqref{eq:Levy}, we thus construct a solution 
of~\eqref{eq:FN_MF} and $N$  independent copies of a solution of~\eqref{eq:FN_limit}.

%
%Subsection
%

\subsection{Main proof and results}
\begin{proposition}\label{prop:majorfG}
	We denote $r^i_{t}= r(Z^{i,N}_{t}, \bar{Z}^i_{t})$ and $G^i_{t} = G^i 
	({(Z^{j,N}_{t})}_j,{(\bar{Z}^j_{t})}_j)$.
	For all $c \in \mathbb{R}$, for each $i \in \{1, \dots, N\}$, we have
	\begin{equation}\label{eq:time_evol_dist}
		d(e^{ct}f(r^i_{t})G^i_{t})\leq e^{ct}K^i_{t}dt+dM^i_{t},
	\end{equation}
	where $M^i_{t}$ is a continuous local martingale and $K^i_{t}$ can be written as 
	\begin{equation}\label{eq:K}
		K^i_{t} = \tilde{K}^i_{t} + I^{1,i}_{t} + I^{2,i}_{t} + I^{3,i}_{t} .
	\end{equation}
	We define $\tilde{K}^i_{t}$, $I^{1,i}_{t}$, $I^{2,i}_{t}$ and $I^{3,i}_{t}$ as follows
	\begin{align}
		\tilde{K}^i_{t}=&G^i_{t}\Big[2cf(r^i_{t})+\frac{1}{2}f''(r^i_{t})\left(2\sigma_{X}^2 
		\varphi_{\text{rc}}{\left(|X^{i,N}_{t}-\bar{X}^i_{t}|\right)}^2\right) \nonumber\\
		&+f'(r^i_{t})\Big((1+\gamma\delta+L_{X}+\delta 
		L_{C})|X^{i,N}_{t}-\bar{X}^i_{t}|-|{{(X^{i,N}_{t})}}^3-{(\bar{X}^i_{t})}^3| \nonumber\\
		&+(1+L_{X}+ \delta L_{C}-\delta)|C^{i,N}_{t}-\bar{C}^i_{t}|+\left(\epsilon 
		\mathcal{C}_{f,1}+\mathcal{C}_{f,2}\right)\sigma_{X}^2\varphi_{\text{rc}}{\left(|X^{i,N}_{t}-\bar{X}^i_{t}|\right)}^2r^i_{t}\Big)\Big]
		 \nonumber\\
		&+\epsilon f(r^i_{t})\left(4\tilde{B}-\frac{\lambda}{16} 
		\tilde{H}(\bar{Z}^i_{t})-\frac{\lambda}{16} 
		\tilde{H}{(Z^{i,N}_{t})}-\frac{\lambda}{16N}\sum_{j=1}^{N} 
		\tilde{H}(\bar{Z}^j_{t})-\frac{\lambda}{16N}\sum_{j=1}^{N} 
		\tilde{H}(Z^{j,N}_{t})\right)\label{eq:tildeK},\\
		I^{1,i}_{t} =& 
		G^i_{t}f'(r^i_{t})\left(\left|\frac{1}{N}\sum_{j=1}^{N}K_{X}(\bar{Z}^i_{t}-\bar{Z}^j_{t})-K_{X}\ast\bar{\mu}_{t}(\bar{Z}^i_{t})\right|\right)\nonumber\\
		&\hspace{2cm}+\delta 
		G^i_{t}f'(r^i_{t})\left(\left|\frac{1}{N}\sum_{j=1}^{N}K_{C}(\bar{Z}^i_{t}-\bar{Z}^j_{t})-K_{C}\ast\bar{\mu}_{t}(\bar{Z}^i_{t})\right|\right)
		 \label{eq:I1},\\
		I^{2,i}_{t}= &G^i_{t}f'(r^i_{t})\left(\frac{L_{X}}{N}\left(\sum_{j=1}^{N} 
		\|Z^{j,N}_{t}-\bar{Z}^j_{t}\|_1 
		\right)\right)+\delta G^i_{t}f'(r^i_{t})\left(\frac{L_{C}}{N}\left(\sum_{j=1}^{N} 
		\|Z^{j,N}_{t}-\bar{Z}^j_{t}\|_1\right)\right) \nonumber\\
		&-cf(r^i_{t})G^i_{t}-\epsilon 
		f(r^i_{t})\left[\frac{\lambda}{16}H(\bar{Z}^i_{t})\exp\left(a\sqrt{H(\bar{Z}^i_{t})}\right)+\frac{\lambda}{16}H{(Z^{i,N}_{t})}\exp\left(a\sqrt{H{(Z^{i,N}_{t})}}\right)\right]
		 \nonumber\\
		&-\epsilon f(r^i_{t})\left[ 
		\frac{\lambda}{16N}\sum_{j=1}^{N}H(\bar{Z}^j_{t})\exp\left(a\sqrt{H(\bar{Z}^j_{t})}\right) 
		\right.\nonumber\\
		&\quad \quad \quad \quad \quad \left. 
		+\frac{\lambda}{16N}\sum_{j=1}^{N}H(Z^{j,N}_{t})\exp\left(a\sqrt{H(Z^{j,N}_{t})}\right)\right],
		\label{eq:I2}\\
		 I^{3,i}_{t}=& \epsilon 
		 f(r^i_{t})\left(\left(\alpha_{X}L_{X}+\beta_{X}L_{C}\right){\left(\frac{\sum_{j=1}^{N}|X^{j,N}_{t}|}{N}\right)}^2\exp\left(a\sqrt{H{(Z^{i,N}_{t})}}\right)\right.
		  \nonumber\\
		&\left. +\left(\alpha_{C}L_{X}+ 
		\beta_{C}L_{C}\right){\left(\frac{\sum_{j=1}^{N}|C^{j,N}_{t}|}{N}\right)}^2\exp\left(a\sqrt{H{(Z^{i,N}_{t})}}\right)\right.
		 \nonumber\\
		&\left. 
		-\frac{\lambda}{16}H{(Z^{i,N}_{t})}\exp\left(a\sqrt{H{(Z^{i,N}_{t})}}\right)-\frac{\lambda}{16N}\sum_{j=1}^{N}H(Z^{j,N}_{t})\exp\left(a\sqrt{H(Z^{j,N}_{t})}\right)\right).\label{eq:I3}
	\end{align}
\end{proposition}
We need to control $\mathbb{E}(G^i_{t})$. This control is a consequence of Lyapunov's 
properties on $\tilde{H}$ and the initial assumption of the Theorem~\ref{thm:unif}. A proof is 
given in Appendix~\ref{subsec:preuve_lem_majorderivertildeH}.
\begin{lemma}\label{lem:controlG}
	There exist $\mathcal{C}_{G,1}$ and $\mathcal{C}_{G,2}$, such that for each $i \leq N$, for 
	all $t>0$, we have 
	\begin{align*}
		\mathbb{E}(G^i_{t})&\leq \mathcal{C}_{G,1}\quad\text{ and 
		}\quad\mathbb{E}[{(G^i_{t})}^2] 
		\leq \mathcal{C}_{G,2}. 
	\end{align*}
\end{lemma}

The decomposition given in the Proposition~\ref{prop:majorfG} is true for all $c \in 
\mathbb{R}$. To control exactly the behavior of each term, we will now consider $c$, defined 
in Subsection~\ref{subsec:hyp_{C}onstants}.

Each term given in Proposition~\ref{prop:majorfG} will be controlled differently. The following 
lemmas summarize it. The first term, $\tilde{K}^i_{t}$, contains the various behaviors we have 
previously identified: we deal with it either through a synchronous coupling (when the 
deterministic drift is contracting) or through a reflection coupling (notice the second derivative 
$f''$ which will provide contraction provided $f$ is sufficiently concave). Finally, notice the 
effect of the Lyapunov function $\tilde{H}$ which yields a restoring force.

\begin{lemma}\label{lem:majorTildeK}
	With the parameters and functions given in Subsection~\ref{subsec:hyp_{C}onstants}, for 
	each $i \leq N$, for all $t >0$, 
	\begin{align}
		\mathbb{E}\tilde{K}^i_{t} \leq \xi\left(2+\delta\gamma+L_{X}+\delta 
		L_{C}-L_{C}-\frac{1+L_{X}}{\delta}\right)\mathbb{E}G^i_{t} .
	\end{align}
\end{lemma}

The interaction term $\frac{1}{N}\sum_j K_{X}(Z^{j,N}_{t} - Z^{i,N}_{t}) - K_{X} \ast 
\bar{\mu}_{t}(\bar{Z}^i_{t})$ can be decomposed into the following two parts: ${\frac{1}{N} 
\sum_{j}K_{X}(\bar{Z}^{j}_{t} - \bar{Z}^{i}_{t})- K_{X} \ast \bar{\mu}_{t}(\bar{Z}^i_{t})}$ and 
$\frac{1}{N} \sum_{j} [K_{X}(Z^{j,N}_{t} - Z^{i,N}_{t})$ $- K_{X}(\bar{Z}^{j}_{t} - 
\bar{Z}^{i}_{t})]$. The first part, which is in $I^{1,i}_{t}$, is dealt with using some form of 
the law 
of large numbers in a similar way to what has been done in the proof of 
Theorem~\ref{thm:non_unif}.
\begin{lemma}\label{lem:I1}
	With the parameters and functions given in Subsection~\ref{subsec:hyp_{C}onstants}, for 
	each $i \leq N$, for all $t >0$, 
	\begin{align}
		\mathbb{E}(I^{1,i}_{t}) \leq 4\sqrt{\frac{\mathcal{C}_{init,2}\mathcal{C}_{G,2}}{N}} (L_{X} 
		+ 
		L_{C}), 
	\end{align}
	where $\mathcal{C}_{G,2}$ is defined in Lemma~\ref{lem:controlG} and $\mathcal{C}_{init,2}$ is defined in Lemma~\ref{lem:borne_unif_moment_2}. 
\end{lemma}

$I^{2,i}_{t}$ contains the leftovers of this decomposition and some of the additional terms of 
the Lyapunov function.
\begin{lemma}\label{lem:I2}
	With the parameters and functions given in Subsection~\ref{subsec:hyp_{C}onstants}, for 
	all $t>0$, 
	\begin{align}
		\frac{1}{N} \sum_{i=1}^{N} I^{2,i}_{t} \leq 0. 
	\end{align}
\end{lemma}

Finally, $I^{3,i}_{t}$ deals with the non-linearity appearing in the dynamics of the Lyapunov 
function, and will be non-positive for values of $L_{X}$ and $L_{C}$ sufficiently small. It is 
also 
here we justify adding the last two terms in~\eqref{eq:def_G}.

\begin{lemma}\label{lem:I3}
	With the parameters and functions given in Subsection~\ref{subsec:hyp_{C}onstants}, for 
	each $i \leq N$, for all $t>0$, 
	\begin{align}
		I^{3,i}_{t} \leq 0.
	\end{align}
\end{lemma}

\begin{proof}[Proof of Theorem~\ref{thm:unif}]
	With these four Lemmas, we can calculate
	\begin{align*}
		\frac{1}{N}&\sum_{i=1}^{N}\mathbb{E}K^i_{t} = \frac{1}{N}\sum_{i=1}^{N}\mathbb{E} 
		\tilde{K}^i_{t} + \frac{1}{N}\sum_{i=1}^{N}\mathbb{E} I^{1,i}_{t} + 
		\frac{1}{N}\sum_{i=1}^{N}\mathbb{E} I^{2,i}_{t} + \frac{1}{N}\sum_{i=1}^{N}\mathbb{E} 
		I^{3,i}_{t}\\
		%\leq & \frac{1}{N}\sum_{i=1}^{N} \xi\left(2+\delta\gamma+L_{X}+\delta 
		%L_{C}-L_{C}-\frac{1+L_{X}}{\delta}\right)\mathbb{E}G^i_{t} + \frac{1}{N}\sum_{i=1}^{N} 
		%4\sqrt{\frac{\mathcal{C}_{init,2}\mathcal{C}_{G,2}}{N}} (L_{X} + L_{C})\\
		\leq& \xi\left(2+\delta\gamma+L_{X}+\delta L_{C}-L_{C}-\frac{1+L_{X}}{\delta}\right) 
		\frac{1}{N}\sum_{i=1}^{N} \mathbb{E}G^i_{t} + 
		4\sqrt{\frac{\mathcal{C}_{init,2}\mathcal{C}_{G,2}}{N}} (L_{X} + L_{C})
	\end{align*}
	Since by Lemma~\ref{lem:controlG}, we have$\frac{1}{N}\sum_{i=1}^{N} \mathbb{E}G^i_{t} 
	\leq \mathcal{C}_{G,1}$, we obtain
	\begin{align*}
		\frac{1}{N}\sum_{i=1}^{N}\mathbb{E}K^i_{t} \leq& \xi A+(L_{X}+L_{C})\frac{B}{\sqrt{N}}
	\end{align*}
	where $A$ and $B$ are constants. 
	
	For all initial couplings such that $\mathbb{E} \rho\left({(Z^{j,N}_0,\bar{Z}^j_0)}_{1\leq j \leq 
	N} \right) < \infty$, by taking the expectation of~\eqref{eq:time_evol_dist} along a sequence 
	of increasing localizing stopping times, we have thanks to Fatou's lemma 
	\begin{align*}
		e^{c t}\mathbb{E}\left( \rho\left({(Z^{j,N}_{t},\bar{Z}^j_{t})}_{1\leq j \leq N} \right) 
		\right)\leq&\mathbb{E}\left( \rho\left({(Z^{j,N}_0,\bar{Z}^j_0)}_{1\leq j \leq N} \right) 
		\right)+\xi A\int_{0}^{t}e^{c s}ds\\
		&+(L_{X}+L_{C})\frac{B}{\sqrt{N}}\int_{0}^{t}e^{c s}ds .\\
		%\leq & \mathbb{E}\left(\rho\left({(Z^{j,N}_0,\bar{Z}^j_0)}_{1\leq j \leq N} \right) \right) + 
		%\xi 
		%A \frac{e^{c t}-1}{c } + (L_{X}+L_{C})\frac{B}{\sqrt{N}} \frac{e^{c t}-1}{c }.
	\end{align*}
	We obtain
	\begin{align*}
		\mathbb{E}\left(\rho\left({(Z^{j,N}_{t},\bar{Z}^j_{t})}_{1\leq j \leq N} \right) \right)\leq & 
		\mathbb{E}\left(\rho\left({(Z^{j,N}_0,\bar{Z}^j_0)}_{1\leq j \leq N} \right)\right)e^{-c t} + 
		\frac{\xi A }{c }\left(1-e^{-c t}\right) \\
		&+\frac{(L_{X}+L_{C})B}{c }\frac{1}{\sqrt{N}}\left(1-e^{-c t}\right).
	\end{align*}
	By using the exchangeability of the particles, we have $\mathbb{E}\left( 
	\rho\left({(Z^{j,N}_{t},\bar{Z}^j_{t})}_{1\leq j \leq N} \right) 
	\right)=\mathbb{E}\left(\frac{1}{N}\sum_{i=1}^{N} f(r^i_{t})G^i_{t} 
	\right)=\mathbb{E}\left(\frac{1}{k}\sum_{i=1}^k f(r^i_{t})G^i_{t} \right)$ for all 
	$k\in\mathbb{N}$. Then
	\[
	\mathbb{E}\left(\sum_{i=1}^k f(r^i_{t})G^i_{t} \right) = k \mathbb{E}\left( 
	\rho\left({(Z^{j,N}_{t},\bar{Z}^j_{t})}_{1\leq j \leq N} \right) \right).
	\]
	Let $\mu_0$ be a measure on $\mathbb{R}^2$, $\mu^{k,N}_{t}$ the marginal distribution at 
	time $t$ of the first $k$ neurons 
	$\left((X^{1,N}_{t},C^{1,N}_{t}),\ldots,(X^{k,N}_{t},C^{k,N}_{t})\right)$ of an $N$-particle 
	system~\eqref{eq:FN_MF} with initial distribution $\mu_0^{\otimes N}$, and $\bar \mu_{t}$ 
	a solution of~\eqref{eq:FN_limit} with initial distribution $\mu_0$. This implies 
	$\mathbb{E}\left(\rho\left({(Z^{j,N}_0,\bar{Z}^j_0)}_{1\leq j \leq N} \right)\right)=0$. By 
	Lemma~\ref{lem:rho_1_2}, we obtain for the $L^1$ Wasserstein distance
	\begin{align*}
		\mathcal{W}_1(\mu^{k,N}_{t}, \bar{\mu}^{\otimes k}_{t}) = & \inf \left\lbrace \mathbb{E}[ \| 
		Z^{(k)}- \bar{Z}^{(k)} \|_1], \mathbb{P}_{Z^{(k)}} = \mu^{k,N}_{t}, 
		\mathbb{P}_{\bar{Z}^{(k)}} = \bar{\mu}^{\otimes k}_{t} \right\rbrace \\
		%=& \inf \left\lbrace \mathbb{E}\left[ \sum_{i=1}^k \| Z^{i,N}_{t}- \bar{Z}^i_{t} \|_1\right], 
		%\mathbb{P}_{{(Z^{i,N}_{t})}_i} = \mu^{k,N}_{t}, \mathbb{P}_{{(\bar{Z}^i_{t})}_i} = 
		%\bar{\mu}^{\otimes k}_{t} \right\rbrace\\
		\leq & \inf \left\lbrace \mathcal{C}_1 \mathbb{E}\left[ \sum_{i=1}^k f(r^i_{t})G^i_{t}\right], 
		\mathbb{P}_{{(Z^{i,N}_{t})}_i} = \mu^{k,N}_{t}, \mathbb{P}_{{(\bar{Z}^i_{t})}_i} = 
		\bar{\mu}^{\otimes k}_{t} \right\rbrace \\
		\leq & \inf \left\lbrace k \mathcal{C}_1 \mathbb{E}\left( 
		\rho\left({(Z^{j,N}_{t},\bar{Z}^j_{t})}_{1\leq j \leq N} \right) \right), 
		\mathbb{P}_{{(Z^{i,N}_{t})}_i} 
		= \mu^{k,N}_{t}, \mathbb{P}_{{(\bar{Z}^i_{t})}_i} = \bar{\mu}^{\otimes k}_{t} \right\rbrace\\
		\leq& \frac{\xi A k\mathcal{C}_1}{c}\left(1-e^{-c t}\right) 
		+\frac{(L_{X}+L_{C})B\mathcal{C}_1}{c }\frac{k}{\sqrt{N}}\left(1-e^{-c t}\right)
	\end{align*}
	By taking the limit as $\xi \rightarrow 0$ uniformly in time, we obtain the desired result. The same lemma and the same type of calculations yield the result for the $L^2$ Wasserstein distance
	\begin{align*}
		\mathcal{W}_2{(\mu^{k,N}_{t}, \bar{\mu}^{\otimes k}_{t})}^2 \leq & 
		\frac{k}{\sqrt{N}}\mathcal{C}_2 \frac{(L_{X}+L_{C})B}{c }.
	\end{align*}
\end{proof}

\subsection{Proof of the decomposition}
\begin{proof}[Proof of Proposition~\ref{prop:majorfG}]
	First, we need to calculate $d(e^{ct} f(r^i_{t}) G^i_{t})$, where we recall 
	\[
	r^i_{t}=|X^{i,N}_{t}-\bar{X}^i_{t}|+\delta|C^{i,N}_{t}-\bar{C}^i_{t}|
	\] 
	and 
	\[
	G^i_{t}=1+\epsilon \tilde{H}(\bar{Z}^i_{t})+\epsilon 
	\tilde{H}{(Z^{i,N}_{t})}+\frac{\epsilon}{N}\sum_{j=1}^{N} 
	\tilde{H}(Z^{j,N}_{t})+\frac{\epsilon}{N}\sum_{j=1}^{N} \tilde{H}(\bar{Z}^j_{t}).
	\]
	We have already calculated $d(X^{i,N}_{t}-\bar{X}^i_{t})$ and $d|X^{i,N}_{t}-\bar{X}^i_{t}|$ 
	in the 
	case of symmetric coupling in Subsection~\ref{subsec:NonUniformTimeChaos} 
	in~\eqref{eq:dX}. Here, we need to use Ito's formula and usual convergence lemmas, see 
	Lemma~\ref{lem:ito_L1} below, to take care of the Brownian term (recall that the coefficient 
	in front of the Brownian vanishes in the vicinity of the singularity of the absolute value). We 
	obtain 
	\begin{align*}
		d|X^{i,N}_{t}-\bar{X}^i_{t}|=&A^X_{t}dt+2\text{sign}(X^{i,N}_{t}-\bar{X}^i_{t})\sigma_{X} 
		\varphi_{\text{rc}}\left(|X^{i,N}_{t}-\bar{X}^i_{t}|\right)dB^{i,rc,X}_{t},
	\end{align*}
	with
	\begin{align*}
		A^X_{t}\leq& 
		|X^{i,N}_{t}-\bar{X}^i_{t}|-\left|{{(X^{i,N}_{t})}}^3-{(\bar{X}^i_{t})}^3\right|+\left|C^{i,N}_{t}-\bar{C}^i_{t}\right|\\
		&+\left|\frac{1}{N}\sum_{j=1}^{N}K_{X}(Z^{i,N}_{t}-Z^{j,N}_{t})-K_{X}\ast\bar{\mu}_{t}(\bar{Z}^i_{t})\right|.
	\end{align*}
	Likewise, as it has already been calculated in~\eqref{eq:dC} in 
	Subsection~\ref{subsec:NonUniformTimeChaos}, 
	\begin{equation*}%\label{dyn_diff_{C}}
		d|C^{i,N}_{t}-\bar{C}^i_{t}|=A^C_{t}dt,
	\end{equation*}
	with
	\begin{align*}
		A^C_{t}\leq\gamma\left|X^{i,N}_{t}-\bar{X}^i_{t}\right|-|C^{i,N}_{t}-\bar{C}^i_{t}|+\left|\frac{1}{N}\sum_{j=1}^{N}K_{C}(Z^{i,N}_{t}-Z^{j,N}_{t})-K_{C}\ast\bar{\mu}_{t}(\bar{Z}^i_{t})\right|.
	\end{align*}
	Now we have
	\begin{align*}
		dr^i_{t}=&\left(A^X_{t}+\delta 
		A^C_{t}\right)dt+2\text{sign}(X^{i,N}_{t}-\bar{X}^i_{t})\sigma_{X} 
		\varphi_{\text{rc}}\left(|X^{i,N}_{t}-\bar{X}^i_{t}|\right)dB^{i,rc,X}_{t}
	\end{align*}
	and we deduce with Ito's formula
	\begin{align*}
		df(r^i_{t})=f'(r^i_{t})dr^i_{t}+\frac{1}{2}f''(r^i_{t}){\left(2\sigma_{X} 
		\varphi_{\text{rc}}\left(|X^{i,N}_{t}-\bar{X}^i_{t}|\right)\right)}^2 dt.
	\end{align*}
	Finally, for $c>0$,
	\begin{align*}
		d(e^{ct}f(r^i_{t}))=ce^{ct}f(r^i_{t})dt+e^{ct}df(r^i_{t}).
	\end{align*}
	Then, by Ito's formula,
	\begin{align*}
		\frac{1}{\epsilon}dG^i_{t}=&\left(\mathcal{L}_{\bar{\mu}_{t}}\tilde{H}(\bar{Z}^i_{t})+\mathcal{L}^{N}
		 \tilde{H}{(Z^{i,N}_{t})}\right)dt\\
		&+\sigma_{X} 
		\varphi_{\text{rc}}\left(|X^{i,N}_{t}-\bar{X}^i_{t}|\right)\left(\partial_{X}\tilde{H}{(Z^{i,N}_{t})}-\partial_{X}\tilde{H}(\bar{Z}^i_{t})\right)dB^{i,rc,X}_{t}\\
		&+\sigma_{X} 
		\varphi_{\text{sc}}\left(|X^{i,N}_{t}-\bar{X}^i_{t}|\right)\left(\partial_{X}\tilde{H}{(Z^{i,N}_{t})}+\partial_{X}\tilde{H}(\bar{Z}^i_{t})\right)dB^{i,sc,X}_{t}\\
		&+\sigma_{C} 
		\left(\partial_{C}\tilde{H}{(Z^{i,N}_{t})}+\partial_{C}\tilde{H}(\bar{Z}^i_{t})\right)dB^{i,C}_{t}\\
		&+\frac{1}{N}\sum_{j=1}^{N}\left(\mathcal{L}_{\bar{\mu}_{t}}\tilde{H}(\bar{Z}^j_{t})+\mathcal{L}^{N}
		 \tilde{H}(Z^{j,N}_{t})\right)dt\\
		&+\frac{\sigma_{X}}{N}\sum_{j=1}^{N} 
		\varphi_{\text{rc}}\left(|X^{j,N}_{t}-\bar{X}^j_{t}|\right)\left(\partial_{X}\tilde{H}(Z^{j,N}_{t})-\partial_{X}\tilde{H}(\bar{Z}^j_{t})\right)dB^{j,rc,X}_{t}\\
		&+\frac{\sigma_{X}}{N}\sum_{j=1}^{N} 
		\varphi_{\text{sc}}\left(|X^{j,N}_{t}-\bar{X}^j_{t}|\right)\left(\partial_{X}\tilde{H}(Z^{j,N}_{t})+\partial_{X}\tilde{H}(\bar{Z}^j_{t})\right)dB^{j,sc,X}_{t}\\
		&+\frac{\sigma_{C}}{N}\sum_{j=1}^{N} 
		\left(\partial_{C}\tilde{H}(Z^{j,N}_{t})+\partial_{C}\tilde{H}(\bar{Z}^j_{t})\right)dB^{j,C}_{t}.
	\end{align*}
	We finally get 
	\begin{multline*}
		d(e^{ct}f(r^i_{t}) G^i_{t}) 
		=G^i_{t}d(e^{ct}f(r^i_{t}))+e^{ct}f(r^i_{t})dG^i_{t}\\
		+2\epsilon\left(1+\frac{1}{N}\right)\sigma_{X}^2\varphi_{\text{rc}}{\left(|X^{i,N}_{t}-\bar{X}^i_{t}|\right)}^2\text{sign}(X^{i,N}_{t}-\bar{X}^i_{t})\\
		\times \left(\partial_{X}\tilde{H}{(Z^{i,N}_{t})}-\partial_{X}\tilde{H}(\bar{Z}^i_{t})\right)
		 e^{ct}f'(r^i_{t})dt.
	\end{multline*}
	Now, we need to use the following Lemma, proven in Appendix~\ref{subsec:preuve_lem_majorderivertildeH}, to have a more tractable expression
	\begin{lemma}\label{lem:majorderivertildeH}
		We have the upper bound
		\begin{align*}
			2\epsilon\left(1+\frac{1}{N}\right)\sigma_{X}^2 
			\varphi_{\text{rc}}{\left(|X^{i,N}_{t}-\bar{X}^i_{t}|\right)}^2& 
			\text{sign}(X^{i,N}_{t}-\bar{X}^i_{t})\left(\partial_{X}\tilde{H}{(Z^{i,N}_{t})}-\partial_{X}\tilde{H}(\bar{Z}^i_{t})\right)
			 \\
			&\leq \left(\epsilon 
			\mathcal{C}_{f,1}+\mathcal{C}_{f,2}\right)\sigma_{X}^2\varphi_{\text{rc}}{\left(|X^{i,N}_{t}-\bar{X}^i_{t}|\right)}^2r^i_{t}
			 G^i_{t}.
		\end{align*}
	\end{lemma}
	Eventually, by denoting the terms in $dB^{i,rc,X}_{t}$, $dB^{i,sc,X}_{t}$, $dB^{i,C}_{t}$, 
	\(\ldots\) 
	as the local martingale $dM^i_{t}$, we obtain
	\begin{align*}
		d(e^{ct}&f(r^i_{t}) G^i_{t})\leq G^i_{t} ce^{ct}f(r^i_{t})dt + e^{ct} G^i_{t} 
		f'(r^i_{t})\left(A^X_{t}+\delta A^C_{t}\right)dt \\
		&+e^{ct} G^i_{t}\frac{1}{2}f''(r^i_{t}){\left(2\sigma_{X} 
		\varphi_{\text{rc}}\left(|X^{i,N}_{t}-\bar{X}^i_{t}|\right)\right)}^2 dt\\ &+e^{ct}f(r^i_{t}) 
		\left( 
		\left(\mathcal{L}_{\bar{\mu}_{t}}\tilde{H}(\bar{Z}^i_{t})+\mathcal{L}^{N} 
		\tilde{H}{(Z^{i,N}_{t})}\right)dt + 
		\frac{1}{N}\sum_{j=1}^{N}\left(\mathcal{L}_{\bar{\mu}_{t}}\tilde{H}(\bar{Z}^j_{t})+\mathcal{L}^{N}
		 \tilde{H}(Z^{j,N}_{t})\right)dt \right) \\
		&+\left(\epsilon 
		\mathcal{C}_{f,1}+\mathcal{C}_{f,2}\right)\sigma_{X}^2\varphi_{\text{rc}}{\left(|X^{i,N}_{t}-\bar{X}^i_{t}|\right)}^2r^i_{t}
		 G^i_{t} e^{ct}f'(r^i_{t})dt + dM^i_{t}.
	\end{align*}
	We use~\eqref{eq:dyn_{t}ilde_H_non_lin} to bound 
	$\mathcal{L}_{\bar{\mu}_{t}}\tilde{H}(\bar{Z}^i_{t})$ and~\eqref{eq:dyn_part_{t}ilde_H} to 
	bound 
	$\mathcal{L}^{N} \tilde{H}{(Z^{i,N}_{t})}$. The interaction terms in $A^X_{t}$ and 
	$A^C_{t}$ are 
	decomposed and we define $I^{1,i}_{t}$ as follows
	\begin{align*}
		I^{1,i}_{t} =& 
		G^i_{t}f'(r^i_{t})\left(\left|\frac{1}{N}\sum_{j=1}^{N}K_{X}(\bar{Z}^i_{t}-\bar{Z}^j_{t})-K_{X}\ast\bar{\mu}_{t}(\bar{Z}^i_{t})\right|\right)\\
		&\hspace{2cm}+\delta 
		G^i_{t}f'(r^i_{t})\left(\left|\frac{1}{N}\sum_{j=1}^{N}K_{C}(\bar{Z}^i_{t}-\bar{Z}^j_{t})-K_{C}\ast\bar{\mu}_{t}(\bar{Z}^i_{t})\right|\right).
	\end{align*}
	The second part of the decomposition is grouped in $I^{2,i}_{t}$, with compensating terms 
	that appear with the use of~\eqref{eq:dyn_{t}ilde_H_non_lin} 
	and~\eqref{eq:dyn_part_{t}ilde_H}, to control the sum
	\begin{align*} 
		I^{2,i}_{t}= 
		&G^i_{t}f'(r^i_{t})\left(\frac{L_{X}}{N}\left(\sum_{j=1}^{N}|X^{j,N}_{t}-\bar{X}^j_{t}|+|C^{j,N}_{t}-\bar{C}^j_{t}|\right)\right)
		 \nonumber\\
		&+\delta 
		G^i_{t}f'(r^i_{t})\left(\frac{L_{C}}{N}\left(\sum_{j=1}^{N}|X^{j,N}_{t}-\bar{X}^j_{t}|+|C^{j,N}_{t}-\bar{C}^j_{t}|\right)\right)
		 \nonumber\\
		&-cf(r^i_{t})G^i_{t}-\epsilon 
		f(r^i_{t})\left[\frac{\lambda}{16}H(\bar{Z}^i_{t})\exp\left(a\sqrt{H(\bar{Z}^i_{t})}\right)+\frac{\lambda}{16}H{(Z^{i,N}_{t})}\exp\left(a\sqrt{H{(Z^{i,N}_{t})}}\right)\right]
		 \nonumber\\
		&-\epsilon f(r^i_{t})\frac{\lambda}{16N}\left[ 
		\sum_{j=1}^{N}H(\bar{Z}^j_{t})\exp\left(a\sqrt{H(\bar{Z}^j_{t})}\right) 
		+\sum_{j=1}^{N}H(Z^{j,N}_{t})\exp\left(a\sqrt{H(Z^{j,N}_{t})}\right)\right].
	\end{align*}
	We gather the expectations terms, obtained with~\eqref{eq:dyn_part_{t}ilde_H}, in 
	$I^{3,i}_{t}$, and we keep a fraction of the Lyapunov function to control it
	\begin{align*}
		I^{3,i}_{t}= \epsilon 
		f(r^i_{t})&\left(\left(\alpha_{X}L_{X}+\beta_{X}L_{C}\right){\left(\frac{\sum_{j=1}^{N}|X^{j,N}_{t}|}{N}\right)}^2\exp\left(a\sqrt{H{(Z^{i,N}_{t})}}\right)\right.
		 \nonumber\\
		&\left. +\left(\alpha_{C}L_{X}+ 
		\beta_{C}L_{C}\right){\left(\frac{\sum_{j=1}^{N}|C^{j,N}_{t}|}{N}\right)}^2\exp\left(a\sqrt{H{(Z^{i,N}_{t})}}\right)\right.
		 \nonumber\\
		&\left. 
		-\frac{\lambda}{16}H{(Z^{i,N}_{t})}\exp\left(a\sqrt{H{(Z^{i,N}_{t})}}\right)-\frac{\lambda}{16N}\sum_{j=1}^{N}H(Z^{j,N}_{t})\exp\left(a\sqrt{H(Z^{j,N}_{t})}\right)\right).
	\end{align*}
	Finally, we define $\tilde{K}^i_{t}$ with the leftovers. It will, in particular, give the constraints 
	on $f$ which explain its choice.
	\begin{align*} 
		\tilde{K}^i_{t}=&G^i_{t}\Big[2cf(r^i_{t})+\frac{1}{2}f''(r^i_{t})\left(2\sigma_{X}^2 
		\varphi_{\text{rc}}{\left(|X^{i,N}_{t}-\bar{X}^i_{t}|\right)}^2\right) \nonumber\\
		&+f'(r^i_{t})\Big((1+\gamma\delta+L_{X}+\delta 
		L_{C})|X^{i,N}_{t}-\bar{X}^i_{t}|-|{{(X^{i,N}_{t})}}^3-{(\bar{X}^i_{t})}^3| \nonumber\\
		&+(1+L_{X}+ \delta L_{C}-\delta)|C^{i,N}_{t}-\bar{C}^i_{t}| +\left(\epsilon 
		\mathcal{C}_{f,1}+\mathcal{C}_{f,2}\right)\sigma_{X}^2\varphi_{\text{rc}}{\left(|X^{i,N}_{t}-\bar{X}^i_{t}|\right)}^2r^i_{t}\Big)\Big]
		 \nonumber\\
		&+\epsilon f(r^i_{t})\left(4\tilde{B}-\frac{\lambda}{16} 
		\tilde{H}(\bar{Z}^i_{t})-\frac{\lambda}{16} 
		\tilde{H}{(Z^{i,N}_{t})}-\frac{\lambda}{16N}\sum_{j=1}^{N} 
		\tilde{H}(\bar{Z}^j_{t})-\frac{\lambda}{16N}\sum_{j=1}^{N} \tilde{H}(Z^{j,N}_{t})\right).
	\end{align*}
\end{proof}

%
%Subsection
%

\subsection{Controls of $I^{1,i}_{t}$, $I^{2,i}_{t}$ and $I^{3,i}_{t}$}
\begin{proof}[Proof of Lemma~\ref{lem:I3}]
	Since we assume 
	\[
	\frac{4}{\gamma}\left(\alpha_{X}L_{X}+\beta_{X}L_{C}\right)\leq\frac{\lambda}{32}\quad\text{and}
	\quad 4\left(\alpha_{C}L_{X}+ \beta_{C}L_{C}\right)\leq\frac{\lambda}{32} ,
	\]
	 and since
	\begin{align*}
		H(Z^{j,N}_{t})\exp\left(a\sqrt{H{(Z^{i,N}_{t})}}\right)\leq 
		H{(Z^{i,N}_{t})}\exp\left(a\sqrt{H{(Z^{i,N}_{t})}}\right)+H(Z^{j,N}_{t})\exp\left(a\sqrt{H(Z^{j,N}_{t})}\right)
	\end{align*}
	we obtain
	\begin{align*}
		(\alpha_{X}&L_{X}+\beta_{X}L_{C}){\left(\frac{\sum_{j=1}^{N}|X^{j,N}_{t}|}{N}\right)}^2\exp\left(a\sqrt{H{(Z^{i,N}_{t})}}\right)\nonumber\\
		&+\left(\alpha_{C}L_{X}+ 
		\beta_{C}L_{C}\right){\left(\frac{\sum_{j=1}^{N}|C^{j,N}_{t}|}{N}\right)}^2\exp\left(a\sqrt{H{(Z^{i,N}_{t})}}\right)\nonumber\\
		&-\frac{\lambda}{16N}\left(NH{(Z^{i,N}_{t})}\exp\left(a\sqrt{H{(Z^{i,N}_{t})}}\right)+\sum_{j=1}^{N}H(Z^{j,N}_{t})\exp\left(a\sqrt{H(Z^{j,N}_{t})}\right)\right)\leq
		 0.
	\end{align*} 
	Then, for each $i \leq N$, and for all $t>0$, $I^{3,i}_{t} \leq 0$. 
\end{proof}

\begin{proof}[Proof of Lemma~\ref{lem:I2}] 
	We prove the non-positivity of $\frac{1}{N} \sum_{i=1}^{N} I^{2,i}_{t}$. First, since 
	$f'\left(r^i_{t}\right)\leq1$, we have 
	\begin{align*}
		\frac{1}{N}\sum_{i=1}^{N}&\left(\frac{1}{N} f'\left(r^i_{t}\right)G^i_{t} \sum_{j=1}^{N} 
		\|Z^{j,N}_{t}-\bar{Z}^j_{t} \|_1\right)\\
		%\leq &\frac{1}{N^2} \sum_{i,j=1}^{N} \|Z^{j,N}_{t}-\bar{Z}^j_{t} \|_1 G^i_{t}\\
		\leq &\frac{1}{N}\sum_{i=1}^{N} \|Z^{i,N}_{t}-\bar{Z}^i_{t} \|_1 +\frac{2\epsilon}{N^2} 
		\sum_{i,j=1}^{N} \|Z^{i,N}_{t}-\bar{Z}^i_{t} \|_1 \left(\tilde{H}\left(\bar{Z}^j_{t}\right)+ 
		\tilde{H}(Z^{j,N}_{t})\right), 
	\end{align*}
	and, using Lemma~\ref{lem:rho_1_2} (i)
	\[
	\frac{1}{N}\sum_{i=1}^{N} \|Z^{i,N}_{t}-\bar{Z}^i_{t} \|_1 \leq 
	\frac{\mathcal{C}_1}{N}\sum_{i=1}^{N}f(r^i_{t}) G^i_{t} , 
	\]
	and with Lemma~\ref{lem:rho_1_2} (iii)
	\begin{align*}
		\sum_{i,j=1}^{N}  \|Z^{i,N}_{t}-\bar{Z}^i_{t} \|_1 & \left(\tilde{H}\left(\bar{Z}^j_{t}\right)+ 
		\tilde{H}(Z^{j,N}_{t})\right)
		% \leq &
		%\mathcal{C}_{z} \sum_{i,j=1}^{N} f(r^i_{t}) \left(1 + \epsilon \sqrt{H{(Z^{i,N}_{t})}} + 
		%\epsilon 
		%\sqrt{H(\bar{Z}^i_{t})}\right) \left(\tilde{H}\left(\bar{Z}^j_{t}\right)+ 
		%\tilde{H}(Z^{j,N}_{t})\right)\\
		\leq \mathcal{C}_{z} \sum_{i,j=1}^{N} f(r^i_{t}) \left(\tilde{H}\left(\bar{Z}^j_{t}\right)+ 
		\tilde{H}(Z^{j,N}_{t})\right) \\
		&+ \epsilon \mathcal{C}_{z} \sum_{i,j=1}^{N} f(r^i_{t}) \left(\sqrt{H{(Z^{i,N}_{t})}} + 
		\sqrt{H(\bar{Z}^i_{t})}\right) \left(\tilde{H}\left(\bar{Z}^j_{t}\right)+ 
		\tilde{H}(Z^{j,N}_{t})\right). 
	\end{align*}
	Using~\eqref{eq:control_{t}ilde_H} from Lemma~\ref{lem:control_{t}ilde_H}, we obtain for 
	the 
	first sum 
	\begin{align*}
		\mathcal{C}_{z} \sum_{i,j=1}^{N} f(r^i_{t})&\left(\tilde{H}\left(\bar{Z}^j_{t}\right)+ 
		\tilde{H}(Z^{j,N}_{t})\right)\\
		&\leq \mathcal{C}_{z} \sum_{i,j=1}^{N} 
		f(r^i_{t})\left(H(\bar{Z}^j_{t})\exp\left(a\sqrt{H(\bar{Z}^j_{t})}\right)+H(Z^{j,N}_{t})\exp\left(a\sqrt{H(Z^{j,N}_{t})}\right)\right).
	\end{align*}
	With~\eqref{eq:control_{t}ilde_H_2} from the same Lemma, we obtain for the second sum
	\begin{align*}
		\epsilon \mathcal{C}_{z} & \sum_{i,j=1}^{N} f(r^i_{t}) \left(\sqrt{H{(Z^{i,N}_{t})}} + 
		\sqrt{H(\bar{Z}^i_{t})}\right) \left(\tilde{H}\left(\bar{Z}^j_{t}\right)+ 
		\tilde{H}(Z^{j,N}_{t})\right)\\
		\leq & \epsilon \mathcal{C}_{z} \frac{2}{a} \sum_{i,j=1}^{N} f(r^i_{t}) 
		\left(\sqrt{H{(Z^{i,N}_{t})}} 
		+ 
		\sqrt{H(\bar{Z}^i_{t})}\right) \\
		&\hspace{2cm}\times\left(\sqrt{H\left(\bar{Z}^j_{t}\right)} \exp\left(a 
		\sqrt{H\left(\bar{Z}^j_{t}\right)} \right)+ \sqrt{H\left(Z^{j,N}_{t}\right)} \exp\left(a 
		\sqrt{H\left(Z^{j,N}_{t}\right)}\right)\right).
	\end{align*}
	Since for all $(y_1, y_2, y_3, y_4) \in {(\mathbb{R}^+)}^4$, we have
	\[
	(y_1 + y_2) \left( y_3 e^{a y_3 } + y_4 e^{a y_4} \right) 
	\leq 2 \left(y_1^2 e^{a y_1 } + y_2^2 e^{a y_2 } + y_3^2 e^{a y_3 } + y_4^2 e^{a y_4 
	}\right), 
	\]
	we obtain for this last sum
	\begin{align*}
		\frac{2 \epsilon \mathcal{C}_{z}}{a} &\sum_{i,j=1}^{N} f(r^i_{t}) \left(\sqrt{H{(Z^{i,N}_{t})}} 
		+ 
		\sqrt{H(\bar{Z}^i_{t})}\right)\\ &\hspace{2cm}\times\left(\sqrt{H\left(\bar{Z}^j_{t}\right)} 
		\exp\left(a \sqrt{H\left(\bar{Z}^j_{t}\right)} \right)+ \sqrt{H\left(Z^{j,N}_{t}\right)} 
		\exp\left(a 
		\sqrt{H\left(Z^{j,N}_{t}\right)}\right)\right)\\
		%\leq & \frac{4 \epsilon \mathcal{C}_{z}}{a} \sum_{i,j=1}^{N} 
		%f(r^i_{t})\left(H(\bar{Z}^i_{t})\exp\left(a\sqrt{H(\bar{Z}^i_{t})}\right)+ 
		%H{(Z^{i,N}_{t})}\exp\left(a\sqrt{H{(Z^{i,N}_{t})}}\right)\right) \\ 
		%&+ \frac{4 \epsilon \mathcal{C}_{z}}{a} \sum_{i,j=1}^{N} f(r^i_{t}) \left( 
		%H(\bar{Z}^j_{t})\exp\left(a\sqrt{H(\bar{Z}^j_{t})}\right)+H(Z^{j,N}_{t})\exp\left(a\sqrt{H(Z^{j,N}_{t})}\right)\right)
		% \\
		\leq & \frac{4 \epsilon \mathcal{C}_{z}}{a} N \sum_{i=1}^{N} 
		f(r^i_{t})\left(H(\bar{Z}^i_{t})\exp\left(a\sqrt{H(\bar{Z}^i_{t})}\right)+ 
		H{(Z^{i,N}_{t})}\exp\left(a\sqrt{H{(Z^{i,N}_{t})}}\right)\right) \\ 
		&+ \frac{4 \epsilon \mathcal{C}_{z}}{a} \sum_{i,j=1}^{N} f(r^i_{t}) \left( 
		H(\bar{Z}^j_{t})\exp\left(a\sqrt{H(\bar{Z}^j_{t})}\right)+H(Z^{j,N}_{t})\exp\left(a\sqrt{H(Z^{j,N}_{t})}\right)\right).
	\end{align*}
	Then, by reconsidering the first expression 
	\begin{align*}
		\frac{1}{N} & \sum_{i=1}^{N}\left(\frac{1}{N} f'\left(r^i_{t}\right)G^i_{t} \sum_{j=1}^{N} 
		\|Z^{j,N}_{t}-\bar{Z}^j_{t} \|_1\right) \\
		\leq &\frac{\mathcal{C}_1}{N}\sum_{i=1}^{N}f(r^i_{t}) G^i_{t} \\
		&+ \frac{2\epsilon}{N^2} \mathcal{C}_{z} \sum_{i,j=1}^{N} 
		f(r^i_{t})\left(H(\bar{Z}^j_{t})\exp\left(a\sqrt{H(\bar{Z}^j_{t})}\right)+H(Z^{j,N}_{t})\exp\left(a\sqrt{H(Z^{j,N}_{t})}\right)\right)
		 \\
		&+ \frac{2\epsilon}{N^2} \frac{4 \epsilon \mathcal{C}_{z}}{a} N \sum_{i=1}^{N} 
		f(r^i_{t})\left(H(\bar{Z}^i_{t})\exp\left(a\sqrt{H(\bar{Z}^i_{t})}\right)+ 
		H{(Z^{i,N}_{t})}\exp\left(a\sqrt{H{(Z^{i,N}_{t})}}\right)\right) \\
		&+ \frac{2\epsilon}{N^2} \frac{4 \epsilon \mathcal{C}_{z}}{a} \sum_{i,j=1}^{N} f(r^i_{t}) 
		\left( 
		H(\bar{Z}^j_{t})\exp\left(a\sqrt{H(\bar{Z}^j_{t})}\right)+H(Z^{j,N}_{t})\exp\left(a\sqrt{H(Z^{j,N}_{t})}\right)\right)
	\end{align*}
	This way, by~\eqref{eq:conditions_LXLC1} since
	\[
	L_{X}\mathcal{C}_1\leq \frac{c}{2},\quad2\mathcal{C}_{z}L_{X}\leq\frac{\lambda}{64} \quad 
	\text{ 
	and }\quad L_{X}\epsilon\frac{8\mathcal{C}_{z}}{a}\leq\frac{\lambda}{64} ,
	\]
	we get
	\begin{align*}
		\frac{1}{N} & \sum_{i=1}^{N}\left(\frac{1}{N} f'\left(r^i_{t}\right)G^i_{t} \sum_{j=1}^{N} 
		\|Z^{j,N}_{t}-\bar{Z}^j_{t} \|_1\right) \\
		\leq &\frac{1}{N} \frac{c }{2 L_{X}}\sum_{i=1}^{N}f(r^i_{t}) G^i_{t} \\
		&+ \frac{\epsilon}{N^2} \frac{\lambda}{64 L_{X}} \sum_{i,j=1}^{N} 
		f(r^i_{t})\left(H(\bar{Z}^j_{t})\exp\left(a\sqrt{H(\bar{Z}^j_{t})}\right)+H(Z^{j,N}_{t})\exp\left(a\sqrt{H(Z^{j,N}_{t})}\right)\right)
		 \\
		&+ \frac{\epsilon}{N} \frac{\lambda}{64 L_{X}} \sum_{i=1}^{N} 
		f(r^i_{t})\left(H(\bar{Z}^i_{t})\exp\left(a\sqrt{H(\bar{Z}^i_{t})}\right)+ 
		H{(Z^{i,N}_{t})}\exp\left(a\sqrt{H{(Z^{i,N}_{t})}}\right)\right) \\
		&+ \frac{\epsilon}{N^2} \frac{\lambda}{64 L_{X}} \sum_{i,j=1}^{N} f(r^i_{t}) \left( 
		H(\bar{Z}^j_{t})\exp\left(a\sqrt{H(\bar{Z}^j_{t})}\right)+H(Z^{j,N}_{t})\exp\left(a\sqrt{H(Z^{j,N}_{t})}\right)\right)
		 , 
	\end{align*}
	and we finally obtain ``half'' the result
	\begin{align*}
		&\frac{1}{N} \sum_{i=1}^{N} G^i_{t} f'(r^i_{t})\left(\frac{L_{X}}{N}\left(\sum_{j=1}^{N} 
		\|Z^{j,N}_{t}-\bar{Z}^j_{t} \|_1 \right)\right) -\frac{c }{2} \frac{1}{N} \sum_{i=1}^{N} 
		f(r^i_{t})G^i_{t}\\
		&-\frac{\epsilon}{2} \frac{\lambda}{16 N} \sum_{i=1}^{N} 
		f(r^i_{t})\left[H(\bar{Z}^i_{t})\exp\left(a\sqrt{H(\bar{Z}^i_{t})}\right)+ 
		H{(Z^{i,N}_{t})}\exp\left(a\sqrt{H{(Z^{i,N}_{t})}}\right)\right] \nonumber\\
		&-\frac{\epsilon}{2}  \frac{\lambda}{16N^2} \sum_{i=1}^{N} f(r^i_{t})\left[ 
		\sum_{j=1}^{N}H(\bar{Z}^j_{t})\exp\left(a\sqrt{H(\bar{Z}^j_{t})}\right) 
		+\sum_{j=1}^{N}H(Z^{j,N}_{t})\exp\left(a\sqrt{H(Z^{j,N}_{t})}\right)\right] \leq 0. 
	\end{align*}
	Likewise, by~\eqref{eq:conditions_LXLC1}, since 
	\[
	\delta L_{C}\mathcal{C}_1\leq \frac{c }{2},\quad2\mathcal{C}_{z}\delta 
	L_{C}\leq\frac{\lambda}{64}\quad\text{ and }\quad\delta 
	L_{C}\epsilon\frac{8\mathcal{C}_{z}}{a}\leq\frac{\lambda}{64} ,
	\]
	we obtain the second ``half''
	\begin{align*}
		&\frac{1}{N} \sum_{i=1}^{N} \delta 
		G^i_{t}f'(r^i_{t})\left(\frac{L_{C}}{N}\left(\sum_{j=1}^{N}\|Z^{j,N}_{t}-\bar{Z}^j_{t} \|_1 
		\right)\right)  -\frac{c }{2} \frac{1}{N} \sum_{i=1}^{N} f(r^i_{t})G^i_{t}\\
		&-\frac{\epsilon}{2} \frac{\lambda}{16N} \sum_{i=1}^{N} 
		f(r^i_{t})\left[H(\bar{Z}^i_{t})\exp\left(a\sqrt{H(\bar{Z}^i_{t})}\right)+H{(Z^{i,N}_{t})}\exp\left(a\sqrt{H{(Z^{i,N}_{t})}}\right)\right]
		 \nonumber\\
		&- \frac{\epsilon}{2}\frac{\lambda}{16N^2} \sum_{i=1}^{N} f(r^i_{t})\left[ 
		\sum_{j=1}^{N}H(\bar{Z}^j_{t})\exp\left(a\sqrt{H(\bar{Z}^j_{t})}\right) 
		+\sum_{j=1}^{N}H(Z^{j,N}_{t})\exp\left(a\sqrt{H(Z^{j,N}_{t})}\right)\right] \leq 0.
	\end{align*}
	Eventually, we have proved $\sum_{i=1}^{N} I^{2,i}_{t} \leq 0$.
\end{proof}

\begin{proof}[Proof of Lemma~\ref{lem:I1}]
	Since $f'(r) \leq 1$, we have by Cauchy-Schwarz inequality
	\begin{align*}
		\mathbb{E}&\left(G^i_{t}f'(r^i_{t})\left(\left|\frac{1}{N}\sum_{j=1}^{N}K_{X}(\bar{Z}^i_{t}-\bar{Z}^j_{t})-K_{X}\ast\bar{\mu}_{t}(\bar{Z}^i_{t})\right|\right)\right)\\
		\leq&\mathbb{E}{\left(|G^i_{t}|^2\right)}^{1/2}\mathbb{E}{\left(\left|\frac{1}{N}\sum_{j=1}^{N}K_{X}(\bar{Z}^i_{t}-\bar{Z}^j_{t})-K_{X}\ast\bar{\mu}_{t}(\bar{Z}^i_{t})\right|^2\right)}^{1/2}.
	\end{align*}
	By Lemma~\ref{lem:controlG}, we have for each $i \leq N$, for all $t \geq 0$, 
	$\mathbb{E}[{(G^i_{t})}^2] \leq \mathcal{C}_{G, 2}$. 
	
	Moreover, we notice that the ${(\bar{Z}^j_{t})}_j$ are i.i.d with law $\bar{\mu}_{t}$. Let's 
	denote 
	$\bar{Z}_{t}$ a generic random variable of law $\bar{\mu}_{t}$ independent of 
	$\bar{Z}^i_{t}$. 
	The calculus of the right term of the product has already been done in 
	Subsection~\ref{subsec:NonUniformTimeChaos}, and we 
	have~\eqref{eq:controle_interactionmoyenne} 
	\begin{align*}
		\mathbb{E}\left(\left|\frac{1}{N}\sum_{j=1}^{N}K_{X}(\bar{Z}^i_{t}-\bar{Z}^j_{t})-K_{X}\ast\bar{\mu}_{t}(\bar{Z}^i_{t})\right|^2\right)
		 \leq &\frac{8 L_{X}^2}{N} \mathbb{E}(\| \bar{Z}_{t}\|_1^2) . 
	\end{align*}
	A similar calculation yields
	\begin{align*}
		\mathbb{E}\left(\left|\frac{1}{N}\sum_{j=1}^{N}K_{C}(\bar{Z}^i_{t}-\bar{Z}^j_{t})-K_{C}\ast\bar{\mu}_{t}(\bar{Z}^i_{t})\right|^2\right)
		 \leq \frac{8 L_{C}^2}{N} \mathbb{E}(\| \bar{Z}_{t}\|_1^2) . 
	\end{align*}
	By Lemma~\ref{lem:borne_unif_moment_2}, $\mathbb{E}\left(|\bar{X}_{t}|^2+ 
	|\bar{C}_{t}|^2\right) \leq \mathcal{C}_{init,2}$. In particular, 
	\begin{align*}
		\mathbb{E}\left(\|\bar{Z}_{t}\|_1^2\right) =
		 \mathbb{E} \left( {\left[|\bar{X}_{t}|+ |\bar{C}_{t}|\right]}^2 \right) \leq 2 
		 \mathbb{E}\left(|\bar{X}_{t}|^2+ |\bar{C}_{t}|^2\right) 
		\leq 2 \mathcal{C}_{init,2}.
	\end{align*}
	Thus
	\begin{align*}
		\mathbb{E}&\left(G^i_{t}f'(r^i_{t})\left(\left|\frac{1}{N}\sum_{j=1}^{N}K_{X}(\bar{Z}^i_{t}-\bar{Z}^j_{t})-K_{X}\ast\bar{\mu}_{t}(\bar{Z}^i_{t})\right|\right)\right)\leq
		 L_{X}\mathcal{C}_{G,2}^{1/2} \sqrt{2 \mathcal{C}_{init,2}} \sqrt{\frac{8}{N}} , 
	\end{align*}
	and likewise
	\begin{align*}
		\mathbb{E}\left(G^i_{t}f'(r^i_{t})\left(\left|\frac{1}{N}\sum_{j=1}^{N}K_{C}(\bar{Z}^i_{t}-\bar{Z}^j_{t})-K_{C}\ast\bar{\mu}_{t}(\bar{Z}^i_{t})\right|\right)\right)\leq
		 L_{C}\mathcal{C}_{G,2}^{1/2} \sqrt{2 \mathcal{C}_{init,2}} \sqrt{\frac{8}{N}}.
	\end{align*}
\end{proof}

%
%Subsection
%

\subsection{Contraction in various regions of space}

The goal of this section is to prove Lemma~\ref{lem:majorTildeK}, i.e show that for each $i \leq 
N$, for all $t >0$, we have the following control
\begin{align*}
	\mathbb{E}\tilde{K}^i_{t} \leq \xi\left(2+\delta\gamma+L_{X}+\delta 
	L_{C}-L_{C}-\frac{1+L_{X}}{\delta}\right)\mathbb{E}G^i_{t} .
\end{align*}
Recall
\begin{align*}
	\tilde{K}^i_{t}=&G^i_{t}\Big[2c f(r^i_{t})+\frac{1}{2}f''(r^i_{t})\left(2\sigma_{X}^2 
	\varphi_{\text{rc}}{\left(|X^{i,N}_{t}-\bar{X}^i_{t}|\right)}^2\right)\\
	&+f'(r^i_{t})\Big((1+\gamma\delta+L_{X}+\delta 
	L_{C})|X^{i,N}_{t}-\bar{X}^i_{t}|-|{(X^{i,N}_{t})}^3-{(\bar{X}^i_{t})}^3|\\
	&+(1+L_{X}+ \delta L_{C}-\delta)|C^{i,N}_{t}-\bar{C}^i_{t}| +\left(\epsilon 
	\mathcal{C}_{f,1}+\mathcal{C}_{f,2}\right)\sigma_{X}^2\varphi_{\text{rc}}{\left(|X^{i,N}_{t}-\bar{X}^i_{t}|\right)}^2r^i_{t}\Big)\Big]\\
	&+\epsilon f(r^i_{t})\left(4\tilde{B}-\frac{\lambda}{16} 
	\tilde{H}(\bar{Z}^i_{t})-\frac{\lambda}{16} 
	\tilde{H}{(Z^{i,N}_{t})}-\frac{\lambda}{16N}\sum_{j=1}^{N} 
	\tilde{H}(\bar{Z}^j_{t})-\frac{\lambda}{16N}\sum_{j=1}^{N} \tilde{H}(Z^{j,N}_{t})\right),
\end{align*}
which is a quantity that contains every term we have not yet dealt with. To prove 
Lemma~\ref{lem:majorTildeK}, we divide for each $i\in\{1,\ldots,N\}$ the space into three 
regions 
\begin{align*}
	\mbox{Reg}^i_1=&\left\{(\bar{Z}^i_{t},Z^{i,N}_{t})\text{ s.t. }|\bar{X}^i_{t}-X^{i,N}_{t}|\geq\xi 
	\text{ and } r^i_{t}\leq R\right\},\\
	\mbox{Reg}^i_2=&\left\{(\bar{Z}^i_{t},Z^{i,N}_{t})\text{ s.t. }|\bar{X}^i_{t}-X^{i,N}_{t}|<\xi 
	\text{ and } r^i_{t}\leq R_1\right\},\\
	\mbox{Reg}^i_3=&\left\{(\bar{Z}^i_{t},Z^{i,N}_{t})\text{ s.t. }r^i_{t}> 
	R\right\},
\end{align*}
where $R$ was given in Lemma~\ref{prop:H}, and consider
\begin{align*}
	\frac{1}{N}\sum_{i=1}^{N}\mathbb{E}\tilde{K}^i_{t}=\frac{1}{N}\sum_{i=1}^{N}\left(\mathbb{E}\left(\tilde{K}^i_{t}\mathds{1}_{\mbox{Reg}^i_1}\right)+\mathbb{E}\left(\tilde{K}^i_{t}\mathds{1}_{\mbox{Reg}^i_2}\right)+\mathbb{E}\left(\tilde{K}^i_{t}\mathds{1}_{\mbox{Reg}^i_3}\right)\right).
\end{align*}

%
%Subsubsection
%

\subsubsection{Region 1: $ \xi\leq |X^{i,N}_{t}-\bar{X}^i_{t}|$ and $r^i_{t}\leq R$.}
In this region of space, since $\varphi_{\text{rc}}(|X^{i,N}_{t}-\bar{X}^i_{t}|)=1$, we have
\begin{align*}
	\tilde{K}^i_{t} \mathds{1}_{\mbox{Reg}^i_1} =& 
	\mathds{1}_{\mbox{Reg}^i_1}\left(G^i_{t}\Big[2c f(r^i_{t})+2 \sigma_{X}^2f''(r^i_{t}) + 
	f'(r^i_{t}) 
	\left(\epsilon \mathcal{C}_{f,1}+\mathcal{C}_{f,2}\right)\sigma_{X}^2 r^i_{t} \right.\\
	&+f'(r^i_{t})(1+\gamma\delta+L_{X}+\delta L_{C})|X^{i,N}_{t}-\bar{X}^i_{t}|\Big] \\
	&- G^i_{t} f'(r^i_{t}) (\delta -1-L_{X}- \delta L_{C})|C^{i,N}_{t}-\bar{C}^i_{t}| - G^i_{t} 
	f'(r^i_{t}) 
	|{(X^{i,N}_{t})}^3-{(\bar{X}^i_{t})}^3|\\
	&+\epsilon f(r^i_{t}) 4\tilde{B}\\
	&\left. -\epsilon f(r^i_{t})\left(\frac{\lambda}{16} \tilde{H}(\bar{Z}^i_{t})+\frac{\lambda}{16} 
	\tilde{H}{(Z^{i,N}_{t})}+\frac{\lambda}{16N}\sum_{j=1}^{N} 
	\tilde{H}(\bar{Z}^j_{t})+\frac{\lambda}{16N}\sum_{j=1}^{N} \tilde{H}(Z^{j,N}_{t})\right)\right),
\end{align*}
and since $\tilde{H}(z)\geq 0$, $|X^{i,N}_{t}-\bar{X}^i_{t}|\leq r^i_{t}$, $\delta > 
\frac{1+L_{X}}{1-L_{C}}$ (by the choice given in Subsection~\ref{subsec:hyp_{C}onstants}) 
and 
$1 \leq G^i_{t}$ we have 
\begin{align*}
	\tilde{K}^i_{t}\mathds{1}_{\mbox{Reg}^i_1}\leq& 
	\mathds{1}_{\mbox{Reg}^i_1}G^i_{t}\left[(2c 
	+4\epsilon \tilde{B})f(r^i_{t})+2\sigma_{X}^2f''(r^i_{t})\right.\\
	&\left.+f'(r^i_{t})\left(1+\delta\gamma+L_{X}+ \delta L_{C}+ \left(\epsilon 
	\mathcal{C}_{f,1}+\mathcal{C}_{f,2}\right)\sigma_{X}^2\right)r^i_{t}\right].
\end{align*}
Using the definition $f$ given in~\eqref{eq:def_f} we get
\begin{align*}
	2\sigma_{X}^2&f''(r^i_{t})+f'(r^i_{t})\left(1+\delta \gamma+ L_{X}+ L_{C}+ \left(\epsilon 
	\mathcal{C}_{f,1}+\mathcal{C}_{f,2}\right)\sigma_{X}^2 \right) r^i_{t}\\
	=&~ 2\sigma_{X}^2\phi'(r^i_{t})g(r^i_{t})+2\sigma_{X}^2\phi(r^i_{t})g'(r^i_{t})\\
	&+\phi(r^i_{t})g(r^i_{t})\left(1+\delta \gamma+ L_{X}+ \delta L_{C}+\left(\epsilon 
	\mathcal{C}_{f,1}+\mathcal{C}_{f,2}\right)\sigma_{X}^2\right) r^i_{t}\\
	=&~ 2\sigma_{X}^2\phi(r^i_{t})g'(r^i_{t})=-(2c +4\epsilon \tilde{B})\Phi(r^i_{t}).
\end{align*}
Thus
\begin{align}
	(2c +4\epsilon 
	\tilde{B})f(r^i_{t})+2\sigma_{X}^2f''(r^i_{t})+f'(r^i_{t})&\left(1+\delta\gamma+L_{X}+ 
	\delta L_{C}+\left(\epsilon 
	\mathcal{C}_{f,1}+\mathcal{C}_{f,2}\right)\sigma_{X}^2\right)r^i_{t} 
	\nonumber \\
	&=(2c +4\epsilon \tilde{B})f(r^i_{t})-(2c +4\epsilon \tilde{B})\Phi(r^i_{t}) \label{eq:ineqf}\\
	&\leq0. \nonumber
\end{align}
Eventually, in this region of space
\begin{align*}
	\tilde{K}^i_{t}\mathds{1}_{\mbox{Reg}^i_1}\leq0.
\end{align*}

%
%Subsubsection
%

\subsubsection{Region 2: $|X^{i,N}_{t}-\bar{X}^i_{t}|<\xi$ and $r^i_{t}\leq R$.}
In this region, we can write $\tilde{K}^i_{t}$ as
\begin{align*}
	\tilde{K}^i_{t}&\mathds{1}_{\mbox{Reg}^i_2} =\mathds{1}_{\mbox{Reg}^i_2}G^i_{t}\Big[2c 
	f(r^i_{t})+ \varphi_{\text{rc}}{\left(|X^{i,N}_{t}-\bar{X}^i_{t}|\right)}^2 \left[2\sigma_{X}^2 
	f''(r^i_{t}) 
	+\left(\epsilon \mathcal{C}_{f,1}+\mathcal{C}_{f,2}\right)\sigma_{X}^2 r^i_{t} f'(r^i_{t}) 
	\right]\\
	&+f'(r^i_{t})\left((1+\gamma\delta+L_{X}+\delta L_{C})|X^{i,N}_{t}-\bar{X}^i_{t}|-(\delta 
	-1-L_{X} 
	- 
	\delta L_{C})|C^{i,N}_{t}-\bar{C}^i_{t}|\right)\Big]\\
	&-\mathds{1}_{\mbox{Reg}^i_2}G^i_{t} f'(r^i_{t})|{(X^{i,N}_{t})}^3-{(\bar{X}^i_{t})}^3|\\ 
	&+\mathds{1}_{\mbox{Reg}^i_2}\epsilon f(r^i_{t}) 4\tilde{B}\\
	&-\epsilon f(r^i_{t})\mathds{1}_{\mbox{Reg}^i_2}\left(\frac{\lambda}{16} 
	\tilde{H}(\bar{Z}^i_{t})+\frac{\lambda}{16} 
	\tilde{H}{(Z^{i,N}_{t})}+\frac{\lambda}{16N}\sum_{j=1}^{N} 
	\tilde{H}(\bar{Z}^j_{t})+\frac{\lambda}{16N}\sum_{j=1}^{N} \tilde{H}(Z^{j,N}_{t})\right). 
\end{align*}
Since $r^i_{t} = |X^{i,N}_{t}-\bar{X}^i_{t}|+ \delta |C^{i,N}_{t}-\bar{C}^i_{t}|$ and 
$|X^{i,N}_{t}-\bar{X}^i_{t}|<\xi$, we have $|C^{i,N}_{t}-\bar{C}^i_{t}| \geq (r^i_{t}-\xi)/\delta$. 
Since 
$\delta > \dfrac{1+L_{X}}{1-L_{C}}$, we obtain
\begin{align*}
	\tilde{K}^i_{t}\mathds{1}_{\mbox{Reg}^i_2}\leq&
	\mathds{1}_{\mbox{Reg}^i_2}G^i_{t}\Big[2c f(r^i_{t})+ 
	\varphi_{\text{rc}}{\left(|X^{i,N}_{t}-\bar{X}^i_{t}|\right)}^2 \left[2\sigma_{X}^2 f''(r^i_{t}) 
	+\left(\epsilon \mathcal{C}_{f,1}+\mathcal{C}_{f,2}\right)\sigma_{X}^2 r^i_{t} f'(r^i_{t}) 
	\right]\\
	&+f'(r^i_{t})\left((1+\gamma\delta+L_{X}+\delta L_{C}) \xi -(\delta - \delta L_{C} -1-L_{X} 
	)\frac{r^i_{t}-\xi}{\delta}\right)\Big]\\
	&+\epsilon f(r^i_{t}) \mathds{1}_{\mbox{Reg}^i_2}4\tilde{B}\\
	\leq& \varphi_{\text{rc}}{\left(|X^{i,N}_{t}-\bar{X}^i_{t}|\right)}^2 G^i_{t} 
	\mathds{1}_{\mbox{Reg}^i_2}\left[2\sigma_{X}^2 f''(r^i_{t}) +\left(\epsilon 
	\mathcal{C}_{f,1}+\mathcal{C}_{f,2}\right)\sigma_{X}^2 r^i_{t} f'(r^i_{t}) \right]\\
	&+ \mathds{1}_{\mbox{Reg}^i_2}G^i_{t} f'(r^i_{t})\xi \left[ 1+\gamma\delta+L_{X}+\delta 
	L_{C} 
	+ 
	1 - L_{C} - \frac{1+ L_{X}}{\delta} \right]\\
	& +\mathds{1}_{\mbox{Reg}^i_2} G^i_{t} \left( (2c +4\epsilon \tilde{B})f(r^i_{t}) - r^i_{t} 
	f'(r^i_{t}) 
	\left(1 - L_{C} - \frac{1+ L_{X}}{\delta}\right) \right).
\end{align*}
By~\eqref{eq:ineqf}, 
\begin{align*}
	2\sigma_{X}^2 f''(r^i_{t}) +\left(\epsilon 
	\mathcal{C}_{f,1}+\mathcal{C}_{f,2}\right)\sigma_{X}^2 
	r^i_{t} f'(r^i_{t}) &= -(2c +4\epsilon \tilde{B})\Phi(r^i_{t}) - f'(r^i_{t}) r^i_{t} \left(1+\delta 
	\gamma+ 
	L_{X}+ L_{C} \right)\\
	& \leq 0, 
\end{align*}
and by Lemma~\ref{lem:parameters}
\begin{align*}
	2c +4\epsilon 
	\tilde{B}\leq\left(1-L_{C}-\frac{1+L_{X}}{\delta}\right)\min_{r\in(0,R]}\frac{f'(r)r}{f(r)},
\end{align*}
we obtain
\begin{align*}
	\tilde{K}^i_{t}\mathds{1}_{\mbox{Reg}^i_2}
	\leq& \mathds{1}_{\mbox{Reg}^i_2}G^i_{t} f'(r^i_{t})\xi \left[ 1+\gamma\delta+L_{X}+\delta 
	L_{C} 
	+ 1 - L_{C} - \frac{1+ L_{X}}{\delta} \right].
\end{align*}
Finally, since $f'(r) \leq 1$, 
\begin{align*}
	\mathbb{E}\tilde{K}^i_{t}\mathds{1}_{\mbox{Reg}^i_2}\leq 
	\xi\left(2+\delta\gamma+L_{X}+\delta 
	L_{C}-L_{C}-\frac{1+L_{X}}{\delta}\right)\mathbb{E}G^i_{t}.
\end{align*}

%
%Subsubsection
%

\subsubsection{Region 3: $r^i_{t}\geq R$.}

In this region of space $f'=f''=0$ and $f$ is constant, and we therefore have
\begin{align*}
	\tilde{K}^i_{t} \mathds{1}_{\mbox{Reg}^i_3} =& f(r^i_{t})\mathds{1}_{\mbox{Reg}^i_3} 
	\left[2c G^i_{t} + 4 \epsilon \tilde{B} \right.\\
	& \left.-\frac{\lambda \epsilon}{16} \left( \tilde{H}(\bar{Z}^i_{t})+ 
	\tilde{H}{(Z^{i,N}_{t})}+\frac{1}{N}\sum_{j=1}^{N} 
	\tilde{H}(\bar{Z}^j_{t})+\frac{1}{N}\sum_{j=1}^{N} \tilde{H}(Z^{j,N}_{t}) \right) \right].
\end{align*}
Since $G^i_{t} = 1+\epsilon \tilde{H}(\bar{Z}^i_{t})+\epsilon 
\tilde{H}{(Z^{i,N}_{t})}+\frac{\epsilon}{N}\sum_{j=1}^{N} 
\tilde{H}(Z^{j,N}_{t})+\frac{\epsilon}{N}\sum_{j=1}^{N} \tilde{H}(\bar{Z}^j_{t})$ by 
definition~\eqref{eq:def_G}, we can write
\begin{align*}
	\tilde{K}^i_{t} \mathds{1}_{\mbox{Reg}^i_3} 
	=& f(r^i_{t})\mathds{1}_{\mbox{Reg}^i_3} \left[2c + 4 \epsilon \tilde{B} \right.\\
	&\left.+\epsilon\left(2c -\frac{\lambda }{16}\right) \left( \tilde{H}(\bar{Z}^i_{t})+ 
	\tilde{H}{(Z^{i,N}_{t})}+\frac{1}{N}\sum_{j=1}^{N} 
	\tilde{H}(\bar{Z}^j_{t})+\frac{1}{N}\sum_{j=1}^{N} \tilde{H}(Z^{j,N}_{t}) \right) \right].
\end{align*}
Since $c \leq \lambda/32$ by the choice given in Subsection~\ref{subsec:hyp_{C}onstants}, 
we obtain
\begin{align*}
	\tilde{K}^i_{t}\mathds{1}_{\mbox{Reg}^i_3}\leq f(r^i_{t})\mathds{1}_{\mbox{Reg}^i_3} 
	\left[2c 
	+4\epsilon \tilde{B} -\epsilon\left(\frac{\lambda}{16}-2c 
	\right)(H(\bar{Z}^i_{t})+H{(Z^{i,N}_{t})})\right]. 
\end{align*}
We have chosen $R$ such that, for $z,z'$ satisfying $r\geq R$, we have $H(z)+H(z')\geq80\frac{\tilde{B}}{\lambda}$ by Lemma~\ref{prop:H} (iv). Therefore
\begin{align*}
	\tilde{K}^i_{t}\mathds{1}_{\mbox{Reg}^i_3}\leq& f(r^i_{t})\mathds{1}_{\mbox{Reg}^i_3} 
	\left(2c +4\epsilon \tilde{B} -\epsilon\left(\frac{\lambda}{16}-2c 
	\right)80\frac{\tilde{B}}{\lambda}\right)\\
	=&f(r^i_{t})\mathds{1}_{\mbox{Reg}^i_3} \left(2c \left(1+80\frac{\epsilon 
	\tilde{B}}{\lambda}\right)-\epsilon \tilde{B}\right)
\end{align*}
Lemma~\ref{lem:parameters} and more specifically the inequality
\begin{align*}
	c \leq\frac{1}{2}\frac{\epsilon \tilde{B}}{1+80\frac{\epsilon \tilde{B}}{\lambda}}=\frac{\lambda}{160}\frac{\frac{80\epsilon \tilde{B}}{\lambda}}{1+\frac{80\epsilon \tilde{B}}{\lambda}}
\end{align*}
yields the desired result: $\tilde{K}^i_{t}\mathds{1}_{\mbox{Reg}^i_3} \leq 0$.

%
%
%Appendice
%
%

\appendix

%
%
%Section
%
%

\section{Various technical lemmas}

	%
	%Subsection
	%
	
	\subsection{On Itô's formula for the $L^1$ norm}\label{subsec:ito_L1}
	
	Let us here detail the calculations leading to the use of Itô's formula to derive the dynamics 
	of the $L^1$ norm of the processes. At first glance it should not be possible, as the 
	absolute value is not a twice continuously differentiable function. However, in our case, we 
	consider a diffusion coefficient which is zero around the point of discontinuity of the 
	function. The following lemma is based on the calculations done in Lemma~7 
	of~\cite{DEGZ20}, and relies on an approximation of the absolute value function and usual 
	convergence lemmas. We here give a quite general result.
	
	%Lemma
	
	\begin{lemma}\label{lem:ito_L1}
		Let $(X_{t},X'_{t},C_{t},C'_{t})$ be continuous processes and 
		$F,G:\mathbb{R}^+\times\mathbb{R}^4\mapsto\mathbb{R}$ be two continuous functions. 
		Assume furthermore that there is $R_G>0$ such that $G(t,x,x',c,c')=0$ if $|x-x'|<R_G$ 
		and that $G$ is bounded. Consider the dynamics
		\begin{align*}
			d(X_{t}-X'_{t})=F(t,X_{t},X'_{t},C_{t},C'_{t})dt+G(t,X_{t},X'_{t},C_{t},C'_{t})dB_{t},
		\end{align*}
		where $B$ is a Brownian motion. Then almost surely for all $t\geq0$
		\begin{align}\nonumber
			d|X_{t}-X'_{t}|=&\text{sign}(X_{t}-X'_{t})F(t,X_{t},X'_{t},C_{t},C'_{t})dt\\
			&+\text{sign}(X_{t}-X'_{t})G(t,X_{t},X'_{t},C_{t},C'_{t})dB_{t},\label{eq:ito_L1}
		\end{align}
		where 
		\begin{align*}
			\text{sign}(x)=\left\{\begin{array}{ll}1 & \text{ if }x>0\\
				0& \text{ if }x=0\\
				-1 & \text{ if }x<0.\end{array}\right.
		\end{align*}
	\end{lemma}
	
	%Proof
	
	\begin{proof}
		By the standard Itô's formula for twice continuously differentiable functions, we have
		\begin{align*}
			d{(X_{t}-X'_{t})}^2=&2(X_{t}-X'_{t})F(t,X_{t},X'_{t},C_{t},C'_{t})dt+2(X_{t}-X'_{t})G(t,X_{t},X'_{t},C_{t},C'_{t})dB_{t}\\
			& +G^2(t,X_{t},X'_{t},C_{t},C'_{t})dt.
		\end{align*}
		Consider, for $\eta>0$ (which in the end will go to 0), the function 
		$\psi_\eta(r)={(r+\eta)}^{1/2}$ which is smooth on $[0,\infty[$ and satisfies
		\begin{align*}
			\forall r>0,& \lim_{\eta\rightarrow0}\psi_\eta(r)=r^{1/2}, \lim_{\eta\rightarrow0}2\psi'_\eta(r)=r^{-1/2}, \lim_{\eta\rightarrow0}4\psi''_\eta(r)=-r^{-3/2}\\
			\text{and thus }&\lim_{\eta\rightarrow0}2r\psi''_\eta(r)+\psi'_\eta(r)=0\text{ and }\forall 
			r\in\mathbb{R}, \lim_{\eta\rightarrow0}2r\psi'_\eta(r^2)=\text{sign}(r).
		\end{align*}
		Then
		\begin{align*}
			d\psi_\eta\left({(X_{t}-X'_{t})}^2\right)=&2(X_{t}-X'_{t})\psi'_\eta\left({(X_{t}-X'_{t})}^2\right)F(t,X_{t},X'_{t},C_{t},C'_{t})dt\\
			&+2(X_{t}-X'_{t})\psi'_\eta\left({(X_{t}-X'_{t})}^2\right)G(t,X_{t},X'_{t},C_{t},C'_{t})dB_{t}\\
			&+\psi'_\eta\left({(X_{t}-X'_{t})}^2\right)G^2(t,X_{t},X'_{t},C_{t},C'_{t})dt\\
			&+2{(X_{t}-X'_{t})}^2\psi''_\eta\left({(X_{t}-X'_{t})}^2\right)G^2(t,X_{t},X'_{t},C_{t},C'_{t})dt,
		\end{align*}
		which is just another way of writing that for all $t\geq0$
		\begin{align*}
			\psi_\eta&\left({(X_{t}-X'_{t})}^2\right)=\psi_\eta\left({(X_0-X'_0)}^2\right)\\
			&+\int_0^{t}2(X_s-X'_s)\psi'_\eta\left({(X_s-X'_s)}^2\right)F(s,X_s,X'_s,C_s,C'_s)ds\\
			&+\int_0^{t}2(X_s-X'_s)\psi'_\eta\left({(X_s-X'_s)}^2\right)G(s,X_s,X'_s,C_s,C'_s)dB_s\\
			&+\int_0^{t}\left(\psi'_\eta\left({(X_s-X'_s)}^2\right)+2{(X_s-X'_s)}^2\psi''_\eta\left({(X_s-X'_s)}^2\right)\right)G^2(s,X_s,X'_s,C_s,C'_s)ds
		\end{align*}
		We now compute the limit of each term. First
		\begin{align*}
			\psi_\eta\left({(X_{t}-X'_{t})}^2\right)\xrightarrow[\eta\rightarrow0]{}|X_{t}-X'_{t}|\quad\text{
					 			and 
					 			}\quad\psi_\eta\left({(X_0-X'_0)}^2\right)\xrightarrow[\eta\rightarrow0]{}|X_0-X'_0|.
		\end{align*}
		Then
		\begin{align*}
			\left|2(X_s-X'_s)\psi'_\eta\left({(X_s-X'_s)}^2\right)\right|=\frac{\left|X_s-X'_s\right|}{\sqrt{{(X_s-X'_s)}^2+\eta}}\leq1.
		\end{align*}
		Thus for all $t\geq0$, by dominated convergence (recall $F$ is a continuous function, thus integrable on $[0,t]$), we obtain almost surely
		\begin{align*}
			\int_0^{t}2(X_s-X'_s)\psi'_\eta\left({(X_s-X'_s)}^2\right)&F(s,X_s,X'_s,C_s,C'_s)ds\\
			&\xrightarrow[\eta\rightarrow0]{}\int_0^{t}\text{sign}(X_s-X'_s)F(s,X_s,X'_s,C_s,C'_s)ds,
		\end{align*}
		and by Theorem 2.12 Chapter 4 of~\cite{RY99}, almost surely we have
		\begin{align*}
			\int_0^{t}2(X_s-X'_s)\psi'_\eta\left({(X_s-X'_s)}^2\right)&G(s,X_s,X'_s,C_s,C'_s)dB_s\\
			&\xrightarrow[\eta\rightarrow0]{}\int_0^{t}\text{sign}(X_s-X'_s)G(s,X_s,X'_s,C_s,C'_s)dB_s,
		\end{align*}
		Finally, since $G(s,X_s,X'_s,C_s,C'_s)=0$ if  $|X_s-X'_s|<R_G$ and
		\begin{align*}
			\psi'_\eta\left({(X_s-X'_s)}^2\right)+2{(X_s-X'_s)}^2&\psi''_\eta\left({(X_s-X'_s)}^2\right)\\
			=&\frac{1}{2}\left(\frac{1}{\sqrt{{(X_s-X'_s)}^2+\eta}}-\frac{{(X_s-X'_s)}^2}{{\left({(X_s-X'_s)}^2+\eta\right)}^{3/2}}\right)\\
			=&\frac{1}{2}\frac{\eta}{{\left({(X_s-X'_s)}^2+\eta\right)}^{3/2}}\\
			\leq&\frac{1}{2}\frac{\eta}{|X_s-X'_s|^3},
		\end{align*}
		by dominated convergence we almost surely have
		\begin{align*}
			\int_0^{t}\left(\psi'_\eta\left({(X_s-X'_s)}^2\right)+2{(X_s-X'_s)}^2\psi''_\eta\left({(X_s-X'_s)}^2\right)\right)G^2(s,X_s,X'_s,C_s,C'_s)ds\xrightarrow[\eta\rightarrow0]{}0
		\end{align*}
		Thus for all $t\geq0$ we almost surely have~\eqref{eq:ito_L1}, and continuity allows us to 
		conclude that we almost surely have for all $t\geq0$~\eqref{eq:ito_L1}.
	\end{proof}

%
%Subsection
%

\subsection{On Lemma~\ref{prop:H}}\label{subsec:preuve_lemme_prop_H}

\begin{lemma}\label{lem:prop-H}
	For all $z=(x,c),z'=(x',c')\in \mathbb{R}^d$, denoting $r(z,z')=|x-x'|+\delta|c-c'|$
	\begin{equation}\label{eq:r_min_rho}
		{r(z,z')}^2\leq \frac{16(1+\delta^2)}{\min\left(\gamma,1\right)}\left(H(z)+H(z')\right)\,,
	\end{equation}
	so that, in particular, for any constant $B>0$, if $r(z,z')\geq R=\sqrt{\dfrac{1280(1+\delta^2)B}{\lambda\min(\gamma,1)}}$, then 
	\begin{equation*}
		\lambda H(z)+\lambda H(z')\geq80B.
	\end{equation*}
\end{lemma}

\begin{proof}
	We have $H(z)\geq\frac{\gamma}{4}x^2+\frac{c^2}{4}\geq \frac{1}{4}\min\left(\gamma,1\right)\left(x^2+c^2\right)$. Thus
	\begin{align*}
		{r(z,z')}^2=&{\left(|x-x'|+\delta|c-c'|\right)}^2\\
		%\leq& 2|x-x'|^2+2\delta^2|c-c'|^2\\
		%\leq& 4x^2+4x'^2+4\delta^2c^2+4\delta^2c'^2\\
		\leq& 4(1+\delta^2)(x^2+c^2)+4(1+\delta^2)(x'^2+c'^2)\\
		\leq& 16\frac{(1+\delta^2)}{\min\left(\gamma,1\right)}\left(H(z)+H(z')\right)
	\end{align*}
\end{proof}

%
%Subsection
%

\subsection{Proof of Lyapunov's property of $H$ and its consequences}\label{subsec:preuve_lya}
\paragraph{Lyapunov's property}

\begin{proof}[Proof of Lemma~\ref{lem:Lya_limit}]
	We write the proof for~\eqref{eq:dyn_H}, as it also yields~\eqref{eq:dyn_H_part} by 
	considering $\mu$ to be the empirical measure. We notice
	\begin{align*}
		\partial_{C}H=c+\alpha\quad\text{ and }\quad\partial_{X}H=\gamma x+\beta,
	\end{align*}
	so
	\begin{align*}
		\mathcal{L}_{\mu}H(z)=&\partial_{x}H(z)(x-x^3)+\partial_{x}H(z)K_{X}\ast\mu(z)-c\partial_{c}H(z)\\
		&+\partial_{c}H(z)K_{C}\ast\mu(z)+\frac{\sigma_{X}
		 ^2\gamma}{2}+\frac{\sigma_{C} ^2}{2}\\
		=&(\gamma x+\beta)(x-x^3) -c(c+\alpha) +(\gamma x+\beta)K_{X}\ast\mu(z) \\
		& +(c+\alpha)K_{C}\ast\mu(z) + \frac{\sigma_{X} ^2\gamma}{2}+\frac{\sigma_{C} ^2}{2} .
	\end{align*}
	First, we focus on the interaction terms. We have
	\begin{align*}
		|K_{X}\ast\mu(z)|\leq&\int_{\mathbb{R}^2}|K_{X}(z - z')| \mu(dz')\\
		\leq&\int_{\mathbb{R}^2}L_{X}(\|z\|_1+\|z'\|_1)\mu(dz').
	\end{align*}
	Hence,
	\begin{align*}
		(\gamma x+\beta)&K_{X}\ast\mu(z) \\
		\leq&L_{X}(\gamma |x|+\beta)(|x|+|c|+\mathbb{E}_{\mu}(|X|)+\mathbb{E}_{\mu}(|C|))\\
		\leq&L_{X}\left(\gamma |x|^2+\gamma 
		|x||c|+\gamma|x|\mathbb{E}_{\mu}(|X|)+\gamma|x|\mathbb{E}_{\mu}(|C|)+\beta|x|+\beta|c|+\beta\mathbb{E}_{\mu}(|X|)\right.\\
		& \left.+\beta\mathbb{E}_{\mu}(|C|)\right),
	\end{align*}
	and using Young's inequality $ab\leq\frac{\alpha}{2}a^2+\frac{1}{2\alpha}b^2$ ($\alpha = 
	16$ when we separate the $x$ and $c$ terms, and $\alpha =1$ otherwise on the various 
	terms) we get
	\begin{align*}
		(\gamma x+\beta)&K_{X}\ast\mu(z)\\
		\leq&L_{X}\Big(\gamma 
		|x|^2+8\gamma^2|x|^2+\frac{|c|^2}{32}+\frac{\gamma}{2}|x|^2+\frac{\gamma}{2}\mathbb{E}_{\mu}{(|X|)}^2+8\gamma^2|x|^2+\frac{\mathbb{E}_{\mu}{(|C|)}^2}{32}+\frac{\beta^2}{2}\\
		&+\frac{|x|^2}{2}+8\beta^2+\frac{|c|^2}{32}+\frac{\beta^2}{2}+\frac{1}{2}\mathbb{E}_{\mu}{(|X|)}^2+8\beta^2+\frac{\mathbb{E}_{\mu}{(|C|)}^2}{32}\Big)\\
		=&L_{X}\left(17\beta^2+|x|^2\left(\frac{1}{2}+\frac{3}{2}\gamma+16\gamma^2\right)+\frac{|c|^2}{16}+\mathbb{E}_{\mu}{(|X|)}^2\left(\frac{\gamma}{2}+\frac{1}{2}\right)+\frac{\mathbb{E}_{\mu}{(|C|)}^2}{16}\right).
	\end{align*}
	Likewise
	\begin{align*}
		(c+\alpha)K_{C}\ast\mu(z)\leq&L_{C}\left(17\alpha^2+\frac{17}{2}|x|^2+|c|^2\left(\frac{3}{2}+\frac{3}{32}\right)+\frac{17}{2}
		 \mathbb{E}_{\mu}{(|X|)}^2\right. \\
		& \left. +\mathbb{E}_{\mu}{(|C|)}^2\left(\frac{1}{2}+\frac{1}{32}\right)\right).
	\end{align*}
	The idea is to bound $\lambda H(z)+ \mathcal{L}_{\mu}H(z)$, by distinguishing 3 types of 
	terms: we isolate terms in $\mathbb{E}_{\mu}{(|C|)}^2-c^2$, 
	$\mathbb{E}_{\mu}{(|X|)}^2-x^2$, 
	and we group polynomial terms. Then, we notice the polynomial is upper bounded by a 
	constant $A$. Thus
	\begin{align*}
		\lambda H(z)&+ \mathcal{L}_{\mu}H(z) - \frac{\sigma_{X} ^2\gamma}{2} - 
		\frac{\sigma_{C} 
		^2}{2} \\
		=& \lambda \left(\frac{1}{2}\gamma x^2+\beta x+\frac{1}{2}c^2+\alpha c+H_0\right) + (\gamma x+\beta)(x-x^3) -c(c+\alpha) \\
		&+(\gamma x+\beta)K_{X}\ast\mu(z) +(c+\alpha)K_{C}\ast\mu(z)\\ 
		\leq&\left(\lambda H_0+17\beta^2 L_{X}+17\alpha^2 L_{C}\right)-\gamma x^4-\beta 
		x^3+(1+\lambda)\beta x\\
		&+\left((1+\frac{\lambda}{2})\gamma+L_{X}\left(1+2\gamma+16\gamma^2\right)+17L_{C}\right)x^2\\
		&+\left(\frac{L_{C}}{8}+L_{C}\left(2+\frac{1}{8}\right)-\left(1-\frac{\lambda}{2}\right)\right)c^2-(1-\lambda)
		 \alpha c\\
		&+\left(\frac{L_{X}}{16}+\frac{L_{C}}{2}+\frac{L_{C}}{32}\right)\left(\mathbb{E}_{\mu}{(|C|)}^2-c^2\right)\\
		&+\left(\frac{\gamma}{2}L_{X}+\frac{1}{2}L_{X}+\frac{17}{2}L_{C}\right)\left(\mathbb{E}_{\mu}{(|X|)}^2-x^2\right).
	\end{align*}
	Provided that
	\[
	\frac{L_{X}}{8}+L_{C}\left(2+\frac{1}{8}\right)<1-\frac{\lambda}{2},
	\]
	there is 
	$A\geq0$ such that
	\begin{align*}
		-\gamma x^4-&\beta x^3+(1+\lambda)\beta 
		x+\left((1+\frac{\lambda}{2})\gamma+L_{X}\left(1+2\gamma+16\gamma^2\right)+17L_{C}\right)x^2\\
		&+\left(\frac{L_{C}}{8}+L_{C}\left(2+\frac{1}{8}\right)-\left(1-\frac{\lambda}{2}\right)\right)c^2-(1-\lambda)
		 \alpha c\leq A. 
	\end{align*}
	Hence the result
	\begin{multline*}
		\mathcal{L}_{\mu}H(\bar{z}) \leq B+\left(\alpha_{X} L_{X}+\beta_{X} L_{C}\right) 
		\left(\mathbb{E}_{\mu}{(|X|)}^2-\bar{x}^2\right)\\+\left(\alpha_{C} L_{X}
		+ \beta_{C} 
		L_{C}\right)\left(\mathbb{E}_{\mu}{(|C|)}^2-\bar{c}^2\right)
		-\lambda H(\bar{z}).
	\end{multline*}
\end{proof}

\paragraph{First consequences}
\begin{proof}[Proof of Proposition~\ref{prop:Lya}]
	Inequality~\eqref{eq:gronwall_particles} simply relies on the sum of~\eqref{eq:dyn_H_part} 
	for each $i$ and the fact that $\mathcal{L}^{j,N}\left( H\left(Z^{i,N}_{t}\right) 
	\right)= 0$ for $i \neq j$
	\begin{align*}
		\frac{1}{N}\sum_{i=1}^{N} & \mathcal{L}^{N}\left(H\left(Z^{i,N}_{t}\right)\right)= 
		\frac{1}{N}\sum_{i=1}^{N} \mathcal{L}^{i,N}\left( H\left(Z^{i,N}_{t}\right)\right)\\
		\leq & \frac{1}{N}\sum_{i=1}^{N} \left[B + \left(\alpha_{X} L_{X}+\beta_{X} L_{C}\right) 
		\left({\left(\frac{1}{N} \sum_{k=1}^{N} |X^{k,N}_{t}| \right)}^2 - {(X^{i,N}_{t})}^2 \right) 
		\right. \\
		& + \left. \left(\alpha_{C} L_{X}+ \beta_{C} L_{C}\right) 
		\left({\left(\frac{1}{N}\sum_{k=1}^{N}|C^{k,N}_{t}|\right)}^2 -{(C^{i,N}_{t})}^2 \right) 
		-\lambda 
		H\left(Z^{i,N}_{t}\right)\right]\\
		\leq & B -\lambda \frac{1}{N}\sum_{i=1}^{N} H\left(Z^{i,N}_{t}\right) .
	\end{align*}
	The last inequality uses the fact that ${\left(\frac{1}{N} \sum_{i=1}^{N} |y_i| \right)}^2 - 
	\frac{1}{N} \sum_{i=1}^{N} {(y_i)}^2 \leq 0$ for all ${(y_i)}_{1 \leq i \leq N} \in 
	\mathbb{R}^{N}$. 
\end{proof}

\paragraph{Bounds on the second moments of processes}
We can now prove the uniform in bounds on the second moments of $X^{i,N}_{t}$, 
$C^{i,N}_{t}$, 
$\bar{X}^i_{t}$ and $\bar{C}^i_{t}$ from~\eqref{eq:dyn_H} and~\eqref{eq:dyn_H_part}. Let's 
notice that ${(X^{i,N,\kappa}_{t}, C^{i,N, \kappa}_{t})}_i$ coincides with 
${(X^{i,N}_{t}, C^{i,N}_{t})}_i$ 
before the time $T_\kappa$ defined in Subsection~\ref{subsec:Existence}. Since our interest 
is 
in ${(X^{i,N}_{t}, C^{i,N}_{t})}_i$, we chose to give the proof of the 
Proposition~\ref{prop:carreintegrable_{X}C}. The proof of the 
Lemma~\ref{lem:carreintegrable_{X}CK} is 
very 
similar. 
\begin{proof}[Proof of Proposition~\ref{prop:carreintegrable_{X}C}]
	$K_{X}$ and $K_{C}$ are Lipschitz with constants $L_{X}$ and $L_{C}$ respectively. We 
	do 
	not assume any bounds on these constants. We assume for each $i \leq N$, 
	$\mathbb{E}(|X^{i,N}_0|^2) < +\infty$ and $\mathbb{E}(|C^{i,N}_0|^2) < +\infty$. We have
	\begin{align*}
		d\left(\frac{e^{\lambda t}}{N}\sum_{i=1}^{N}H\left(Z^{i,N}_{t}\right)\right)=\lambda 
		\frac{e^{\lambda t}}{N}\sum_{i=1}^{N}H\left(Z^{i,N}_{t}\right) dt+e^{\lambda 
		t}\mathcal{L}^{N}\left(\frac{1}{N}\sum_{i=1}^{N}H\left(Z^{i,N}_{t}\right)\right)dt+dM_{t},
	\end{align*}
	where $M_{t}$ is a local martingale. Using~\eqref{eq:dyn_H_part}
	\begin{align*}
		d\left(\frac{e^{\lambda t}}{N}\sum_{i=1}^{N}H\left(Z^{i,N}_{t}\right)\right)=A_{t} 
		dt+dM_{t},
	\end{align*}
	where $A_{t}\leq Be^{\lambda t}$. Let $\tau_n$ be an increasing sequence of localizing 
	stopping times converging to $\infty$ for $M_{t}$
	\begin{align*}
		\mathbb{E}\left(\frac{e^{\lambda t\wedge 
		\tau_n}}{N}\sum_{i=1}^{N}H\left(Z^{i,N}_{t\wedge \tau_n}\right)\right)\leq& 
		\mathbb{E}\left(\frac{1}{N}\sum_{i=1}^{N}H\left(Z^{i,N}_{0}\right)\right)+\mathbb{E}\left(\int_0^{t\wedge
		 \tau_n}Be^{\lambda s}ds\right)\\
		\leq&\mathbb{E}\left(\frac{1}{N}\sum_{i=1}^{N}H\left(Z^{i,N}_{0}\right)\right)+B\frac{\mathbb{E}\left(e^{\lambda
		 t\wedge \tau_n}\right)-1}{\lambda}\\
		\leq&\mathbb{E}\left(\frac{1}{N}\sum_{i=1}^{N}H\left(Z^{i,N}_{0}\right)\right)+B\max\left(\frac{e^{\lambda
		 t}-1}{\lambda},\frac{1}{|\lambda|}\right),
	\end{align*}
	where the maximum on this last inequality depends on the sign of $\lambda$. By Fatou's lemma, we obtain
	\begin{align*}
		e^{\lambda 
		t}\mathbb{E}\left(\frac{1}{N}\sum_{i=1}^{N}H\left(Z^{i,N}_{t}\right)\right)=&\mathbb{E}\left(\liminf_{n\to\infty}\frac{e^{\lambda
		 t\wedge \tau_n}}{N}\sum_{i=1}^{N}H\left(Z^{i,N}_{t\wedge \tau_n}\right)\right)\\
		\leq&\liminf_{n\to\infty}\mathbb{E}\left(\frac{e^{\lambda t\wedge 
		\tau_n}}{N}\sum_{i=1}^{N}H\left(Z^{i,N}_{t\wedge \tau_n}\right)\right)\\
		\leq&\mathbb{E}\left(\frac{1}{N}\sum_{i=1}^{N}H\left(Z^{i,N}_{0}\right)\right)+B\max\left(\frac{e^{\lambda
		 t}-1}{\lambda},\frac{1}{|\lambda|}\right).
	\end{align*}
	Hence the various bounds on $\mathbb{E}\left(|X^{i,N}_{t}|^2\right)$ and 
	$\mathbb{E}\left(|C^{i,N}_{t}|^2\right)$, since by Lemma~\ref{prop:H} (i) we have 
	\begin{align*}
	&{\mathbb{E}H\left(Z^{i,N}_{t}\right)\geq 
	\frac{\gamma}{4}\mathbb{E}\left(|X^{i,N}_{t}|^2\right)+\frac{1}{4}\mathbb{E}\left(|C^{i,N}_{t}|^2\right)}\\
	 \text{ and } & \mathbb{E}H\left(Z^{i,N}_0\right) \leq \gamma 
	 \mathbb{E}\left(|X^{i,N}_0|^2\right) + \mathbb{E}\left(|C^{i,N}_0|^2\right) + \frac{3}{2} H_0.
	 \end{align*}
	These bounds are uniform in time provided $\lambda>0$, i.e 
	$\frac{L_{X}}{8}+L_{C}\left(2+\frac{3}{32}\right)<1$.
\end{proof}

\begin{proof}[Proof of Proposition~\ref{prop:carreintegrable_barXC} and 
Lemma~\ref{lem:borne_unif_moment_2}]
	The proof is done in exactly the same way as the proof of 
	Proposition~\ref{prop:carreintegrable_{X}C} above using~\eqref{eq:dyn_H}.
	
\end{proof}

%
%Subsection
%

\subsection{Proof of Lemma~\ref{lem:parameters}}\label{subsec:choix_parametres_section}

We now prove that there are constants $c, \epsilon$ and $\delta$ such that
\begin{align}
	c +2\epsilon \tilde{B}\leq 
	&\frac{\sigma_{X}^2}{2}{\left(\int_0^R\Phi(s){\phi(s)}^{-1}ds\right)}^{-1}\label{eq:cond_{C}_1}\\
	2c +4\epsilon 
	\tilde{B}\leq&\left(1-L_{C}-\frac{1+L_{X}}{\delta}\right)\min_{r\in(0,R]}\frac{f'(r)r}{f(r)}\label{eq:cond_{C}_2}\\
	c \leq&\frac{\lambda}{160}\frac{\frac{80\epsilon \tilde{B}}{\lambda}}{1+\frac{80\epsilon 
	\tilde{B}}{\lambda}}\label{eq:cond_{C}_3}\\
	\delta>&\frac{1+L_{X}}{1-L_{C}}\label{eq:cond_delta}
\end{align}

\begin{itemize}
	\item Since for all $u\geq0$, $0<\phi\left(u\right)\leq1$, we have 
	$0<\Phi\left(s\right)=\int_0^{s}\phi\left(u\right)du\leq s$, i.e $s/\Phi\left(s\right)\geq1$.
	Therefore
	\begin{align*}
		\inf_{r\in(0,R]}\frac{r\phi\left(r\right)}{\Phi\left(r\right)}\geq\inf_{r\in(0,R]}\phi\left(r\right)=\phi\left(R\right).
	\end{align*}
	It is thus sufficient for~\eqref{eq:cond_{C}_2} to have
	\begin{equation*}
		2c +4\epsilon 
		\tilde{B}\leq\frac{1}{2}\left(1-L_{C}-\frac{1+L_{X}}{\delta}\right)\phi\left(R\right).
	\end{equation*}
	
	\item We have
	\begin{align*}
		\phi\left(r\right)\leq\exp\left(-\frac{1}{4\sigma_{X}^2}r^2\right).
	\end{align*}
	So
	\begin{align*}
		\Phi\left(r\right)\leq\int_0^{\infty}\exp\left(-\frac{r^2}{4\sigma_{X}^2}\right)dr=\sigma_{X}\sqrt{\pi}.
	\end{align*}
	Then
	\begin{align*}
		\int_0^{R}\frac{\Phi\left(r\right)}{\phi\left(r\right)}dr\leq\sigma_{X}\sqrt{\pi}R\frac{1}{\phi\left(R\right)}.
	\end{align*}
	It is thus sufficient for~\eqref{eq:cond_{C}_1} that 
	\begin{equation*}
		c +2\epsilon \tilde{B}\leq\frac{\sigma_{X}}{2\sqrt{\pi}}\frac{\phi\left(R\right)}{R}.
	\end{equation*}
	
	\item The various conditions involving $c $ invite us to consider $2\epsilon \tilde{B}=\eta c $. Then
	\begin{align*}
		c \leq\frac{\lambda}{160}\frac{\frac{80\epsilon \tilde{B}}{\lambda}}{1+\frac{80\epsilon \tilde{B}}{\lambda}}&\iff c \leq \frac{\lambda }{160}\frac{40\eta c }{\lambda+40\eta c }\\
		&\iff 1\leq\lambda\frac{\eta}{4\lambda+160\eta c }\text{ (since $c \geq0$)}\\
		&\iff c \leq\frac{\lambda}{160}\frac{\eta-4}{\eta}.
	\end{align*}
	
	\item We choose to write
	\begin{equation*}
		\delta=(1+\tilde{\delta})\frac{1+L_{X,\max}}{1-L_{C,\max}}>\frac{1+L_{X}}{1-L_{C}}
	\end{equation*}
	
	\item Let us assume, for simplicity, that $\epsilon\leq1$. It is sufficient for this later condition to have
	\begin{align*}
		c \leq \frac{2\tilde{B}}{\eta}.
	\end{align*}
	
	\item The appearance of $\phi\left(R\right)$ suggests we should try to minimize it.
	We recall
	\begin{align*}
		\phi(r)=&\exp\left(-\frac{1}{4\sigma_{X}^2}\left(1+\delta\gamma+L_{X}+ \delta 
		L_{C}+\left(\epsilon 
		\mathcal{C}_{f,1}+\mathcal{C}_{f,2}\right)\sigma_{X}^2\right)r^2\right)\\
		\geq&\exp\left(-\frac{1}{4\sigma_{X}^2}\left(1+\delta\gamma+L_{X}+ \delta 
		L_{C}+\left(\mathcal{C}_{f,1}+\mathcal{C}_{f,2}\right)\sigma_{X}^2\right)r^2\right).
	\end{align*}
	It is therefore sufficient for~\eqref{eq:cond_{C}_1} to have
	\begin{align*}
		c 
		\leq\frac{1}{1+\eta}\frac{\sigma_{X}}{2\sqrt{\pi}}\frac{1}{R}\exp\left(-\frac{1}{4\sigma_{X}^2}\left(1+\delta\gamma+L_{X}+
		 \delta L_{C}+\left(\mathcal{C}_{f,1}+\mathcal{C}_{f,2}\right)\sigma_{X}^2\right)R^2\right),
	\end{align*}
	and for~\eqref{eq:cond_{C}_2} to have
	\begin{multline*}
		c 
		\leq\frac{1}{2(1+\eta)}\left(1-L_{C}-\frac{1+L_{X}}{\delta}\right)\\
		\times\exp\left(-\frac{R^2}{4\sigma_{X}^2}\left(1+\delta\gamma+L_{X}+
		 \delta L_{C}+\left(\mathcal{C}_{f,1}+\mathcal{C}_{f,2}\right)\sigma_{X}^2\right)\right).
	\end{multline*}
	\item Finally, we bound $L_{X}$ and $L_{C}$ by either $0$ or $L_{X,\max}$ and 
	$L_{C,\max}$, 
	to obtain bounds on $c $ independent of $L_{X}$ and $L_{C}$.

\end{itemize}

%
%Subsection
%

\subsection{Proof of 
Lemma~\ref{lem:rho_1_2}}\label{subsec:preuve_lemme_{C}ontrole_distance}

Let $z, z'\in \mathbb{R}^2$.
\paragraph{Proof of control of the $L^1$ distance:}
We have
\begin{align*}
	\|z-z'\|_1 = |x-x'|+|c-c'|\leq\frac{1}{\min\left(\delta, 1\right)}\left(|x-x'|+\delta|c-c'|\right)=\frac{1}{\min\left(\delta, 1\right)}r(z,z').
\end{align*}
If $r(z,z')\leq 1\leq R$, we have, using Lemma~\ref{lem:parameters}
\[
	r(z,z')\leq\frac{f(r)}{f'_-(R)}\leq\frac{f(r)}{\phi(R)g(R)}\left(1+\epsilon \tilde{H}(z)+\epsilon \tilde{H}(z')\right).
\]
If $r(z,z')\geq 1$, we have, using~\eqref{eq:r_min_rho}
\begin{align*}
	r(z,z')\leq& {r(z,z')}^2\\
	\leq&\frac{16(1+\delta^2)}{\epsilon\min\left(\gamma,1\right)}\left(\epsilon H(z)+\epsilon H(z')\right)\\
	\leq &\frac{16(1+\delta^2)}{\epsilon\min\left(\gamma,1\right)}\frac{f(r)}{f(1)}\left(1+\epsilon H(z)+\epsilon H(z')\right)\\
	\leq&\frac{16(1+\delta^2)}{\epsilon\min\left(\gamma,1\right)}\frac{f(r)}{\phi(R)g(R)}\left(1+\epsilon \tilde{H}(z)+\epsilon \tilde{H}(z')\right).
\end{align*}
Thus
\begin{align*}
	\|z-z'\|_1 \leq \frac{1}{\min\left(\delta, 1\right)}\frac{1}{\phi(R)g(R)}\max\left(\frac{16(1+\delta^2)}{\epsilon\min\left(\gamma,1\right)},1\right)f(r(z,z'))\left(1+\epsilon \tilde{H}(z)+\epsilon \tilde{H}(z')\right).
\end{align*}

\paragraph{Proof of control of the $L^2$ distance:}

We have
\begin{align*}
	{r(z,z')}^2={\left(|x-x'|+\delta|c-c'|\right)}^2\geq & |x-x'|^2+\delta^2|c-c'|^2\\
	\geq& \min\left(1,\delta^2\right)\left(|x-x'|^2+|c-c'|^2\right).
\end{align*}
If $r(z,z')\geq 1$, we have, using~\eqref{eq:r_min_rho}
\begin{align*}
	{r(z,z')}^2%\leq&\frac{16(1+\delta^2)}{\epsilon\min\left(\gamma,1\right)}\left(\epsilon 
	%H(z)+\epsilon H(z')\right)\\
	%\leq &\frac{16(1+\delta^2)}{\epsilon\min\left(\gamma,1\right)}\frac{f(r)}{f(1)}\left(1+\epsilon H(z)+\epsilon H(z')\right)\\
	\leq&\frac{16(1+\delta^2)}{\epsilon\min\left(\gamma,1\right)}\frac{f(r)}{\phi(R)g(R)}\left(1+\epsilon \tilde{H}(z)+\epsilon \tilde{H}(z')\right).
\end{align*}
If $r(z,z')\leq 1\leq R$, we have, using Lemma~\ref{lem:parameters}
\begin{align*}
	{r(z,z')}^2\leq r(z,z')\leq\frac{f(r)}{f'_-(R)}\leq\frac{f(r)}{\phi(R)g(R)}\left(1+\epsilon 
	\tilde{H}(z)+\epsilon \tilde{H}(z')\right).
\end{align*}
Thus
\begin{align*}
	\|z-z'\|_2^2\leq \frac{1}{\min\left(\delta^2, 1\right)}\frac{1}{\phi(R)g(R)}\max\left(\frac{16(1+\delta^2)}{\epsilon\min\left(\gamma,1\right)},1\right)f(r(z,z'))\left(1+\epsilon \tilde{H}(z)+\epsilon \tilde{H}(z')\right).
\end{align*}

\paragraph{Proof of the second control of the $L^1$ distance:} We have, if $r(z,z')\leq 1\leq R$
\begin{align*}
	r(z,z')\leq\frac{f(r)}{f'_-(R)}\leq\frac{f(r)}{\phi(R)g(R)}\left(1+\epsilon \sqrt{H(z)}+\epsilon \sqrt{H(z')}\right).
\end{align*}
and, if $r(z,z')\geq 1$, recall Lemma~\ref{prop:H}.
\begin{align*}
	\|z-z'\|_1 \leq& \sqrt{\frac{4}{\gamma}H(z)}+\sqrt{\frac{4}{\gamma}H(z')}+\sqrt{4H(z)}+\sqrt{4H(z')}\\
	\leq& 4\max\left(\sqrt{\frac{1}{\gamma}},1\right)\left(\sqrt{H(z)}+\sqrt{H(z')}\right)\\
	\leq&\frac{4}{\epsilon}\max\left(\sqrt{\frac{1}{\gamma}},1\right)\frac{f(r)}{\phi(R)g(R)}\left(1+\epsilon\sqrt{H(z)}+\epsilon\sqrt{H(z')}\right),
\end{align*}
and thus
\begin{align*}
	\|z-z'\|_1 \leq\frac{1}{\phi(R)g(R)}\max\left(1,\frac{4}{\epsilon}\max\left(\sqrt{\frac{1}{\gamma}},1\right)\right)f(r(z,z'))\left(1+\epsilon\sqrt{H(z)}+\epsilon\sqrt{H(z')}\right).
\end{align*}

\paragraph{Independence with respect to $L_{X}$ and $L_{C}$}
The \textit{a priori} bounds $L_{X}\in[0,L_{X,\max}]$ and $L_{C}\in[0,L_{C,\max}]$ allow us to 
bound $\phi(R)$ independently of $L_{C}$ and $L_{X}$ by $\phi_{\min}$ (and we also use 
$g(R)\geq\frac{1}{2}$), thus giving us constant $\mathcal{C}_1$, $\mathcal{C}_2$ and 
$\mathcal{C}_{z}$ independent of $L_{C}$ and $L_{X}$.

%
%
%Subsection
%
%

\subsection{Proof of Lemmas~\ref{lem:controlG} and~\ref{lem:majorderivertildeH}} 
~\label{subsec:preuve_lem_majorderivertildeH}

\begin{proof}[Proof of Lemma~\ref{lem:controlG}]
	Let's prove there exists a uniform in time bound on $\mathbb{E}(G^i_{t})$ and 
	$\mathbb{E}[{(G^i_{t})}^2]$. First, let's recall the definition of $G$ from\eqref{eq:def_G}
	\[
	G^i_{t}=1+\epsilon \tilde{H}(\bar{Z}^i_{t})+\epsilon 
	\tilde{H}{(Z^{i,N}_{t})}+\frac{\epsilon}{N}\sum_{j=1}^{N} 
	\tilde{H}(Z^{j,N}_{t})+\frac{\epsilon}{N}\sum_{j=1}^{N} \tilde{H}(\bar{Z}^j_{t}).
	\]
	The idea is to bound the different expectations in terms of the expectations at time $t=0$.
	 Since $\mathbb{E}(e^{\tilde{a}(|X_0|+|C_0|)})$ is finite, we know that for each 
	 $k \in \mathbb{N}$, $\mathbb{E}(|X_0|^k)$ and $\mathbb{E}(|C_0|^k)$ are also finite. 
	 We deduce that for each $k \in \mathbb{N}$, for each $j\leq N$, 
	 $\mathbb{E}[{H(\bar{Z}^j_0)}^k]$ and $\mathbb{E}[{H(Z^{j,N}_0)}^k]$ are finite. 
	
	In fact, to bound uniformly in time the first moment, we only have to bound 
	$\mathbb{E}(\tilde{H}(Z^{j,N}_{t}))$ and $\mathbb{E}(\tilde{H}(\bar{Z}^j_{t}))$ for each $j 
	\leq 
	N$. Let's begin with $\bar{Z}^j$. By~\eqref{eq:Gronwall_exp_H}, we have 
	\begin{align*}
		\frac{d}{dt}\mathbb{E}\left[ 
		\tilde{H}\left(\bar{Z}^j_{t}\right)\right]\leq&\tilde{B}-\frac{\lambda }{4}\mathbb{E}\left[ 
		\tilde{H}\left(\bar{Z}^j_{t}\right)\right].
	\end{align*}
	By using Itô's formula on $e^{\lambda t/4}\tilde{H}\left(\bar{Z}^j_{t}\right)$ and the bound 
	above, we obtain 
	\begin{align*}
		\mathbb{E} \left[ \tilde{H}\left(\bar{Z}^j_{t}\right)\right] \leq& \frac{4 \tilde{B}}{\lambda} + 
		e^{-\frac{\lambda}{4}t} \left(\mathbb{E}\left[ \tilde{H}\left(\bar{Z}^j_0\right)\right] - 
		\frac{4 \tilde{B}}{\lambda} \right)\\
		\leq& \max\left(\mathbb{E}\left[ \tilde{H}\left(\bar{Z}^j_0\right)\right],\frac{4 \tilde{B}}{\lambda} \right). 
	\end{align*}
	By~\eqref{eq:control_{t}ilde_H}, in Lemma~\ref{lem:control_{t}ilde_H}, we deduce the 
	following inequality and we apply Cauchy-Schwarz inequality
	\begin{align} 
		\mathbb{E}\left[ \tilde{H}\left(\bar{Z}^j_0\right)\right] \leq & \mathbb{E}\left[ H\left(\bar{Z}^j_0\right) \exp \left(a \sqrt{H\left(\bar{Z}^j_0\right)}\right) \right]\nonumber \\
		\leq & \mathbb{E}{\left[ H{\left(\bar{Z}^j_0\right)}^2 \right]}^{1/2} \mathbb{E}{\left[\exp 
		\left(2a \sqrt{H\left(\bar{Z}^j_0\right)}\right) \right]}^{1/2}. \label{eq:CS_{t}ildeH}
	\end{align}
	We already know $\mathbb{E}\left[ H{\left(\bar{Z}^j_0\right)}^2 \right]$ is bounded. Now, it 
	is enough to prove that there exist $C$ such that for all $z \in \mathbb{R}^2$
	\begin{align*}
		\exp \left(2a \sqrt{H\left(z\right)}\right) \leq C \times e^{\tilde{a}(|x|+|c|)}.
	\end{align*}
	In fact, from the definition of $H$ in~\eqref{eq:def_H}, we have
	\begin{align*}
		2 \sqrt{H(z)} =& \sqrt{2} \sqrt{\gamma {\left(x + \frac{\beta}{\gamma}\right)}^2 + 
		{(c+\alpha)}^2 + H_0} \\
		\leq & \sqrt{2 \gamma} \left| x + \frac{\beta}{\gamma} \right| + \sqrt{2} \left|c+\alpha\right| + \sqrt{H_0}\\
		\leq & \sqrt{2 \gamma} |x| + \sqrt{2} |c| + \frac{1}{a}\ln{C}, 
	\end{align*}
	where $C$ is a constant independent of $z$. Finally, since $\max{(a \sqrt{2\gamma}, a\sqrt{2})} \leq \tilde{a}$, we have
	\begin{align*}
		\exp \left(2a \sqrt{H\left(z\right)}\right) \leq C \times e^{\tilde{a}(|x|+|c|)}.
	\end{align*}
	Then, $\mathbb{E}\left[\exp \left(2a \sqrt{H\left(\bar{Z}^j_0\right)}\right) \right]$ is bounded 
	and we deduce $\mathbb{E}(\tilde{H}(\bar{Z}^j_{t}))$ is bounded for each $j \leq N$ and all 
	$t \geq 0$. 
	
	The same calculations can be done for $Z^{j,N}_{t}$. By~\eqref{eq:gron_part_exp}, we have
	\begin{align*}
		\mathcal{L}^{N}&\left(\frac{1}{N}\sum_{i=1}^{N}\tilde{H}{(Z^{i,N}_{t})}\right)\leq 
		\tilde{B}-\frac{\lambda}{4}\left(\frac{1}{N}\sum_{i=1}^{N}\tilde{H}{(Z^{i,N}_{t})}\right). 
	\end{align*}
	In particular, 
	\begin{align*}
		\frac{d}{dt} 
		\left[\mathbb{E}\left(\frac{1}{N}\sum_{i=1}^{N}\tilde{H}{(Z^{i,N}_{t})}\right)\right] 
		\leq & \mathbb{E} \left[ 
		\mathcal{L}^{N}\left(\frac{1}{N}\sum_{i=1}^{N}\tilde{H}{(Z^{i,N}_{t})}\right) \right] \leq 
		\tilde{B}-\frac{\lambda}{4}\mathbb{E}\left(\frac{1}{N}\sum_{i=1}^{N}\tilde{H}{(Z^{i,N}_{t})}\right),
	\end{align*}
	and we can use the same method as above. 
	
	Finally, we have proved that for each $j \leq N$, $\mathbb{E}(\tilde{H}(Z^{j,N}_{t}))$ and 
	$\mathbb{E}(\tilde{H}(\bar{Z}^j_{t}))$ are bounded uniformly in time. Thus, 
	$\mathbb{E}(G^i_{t})$ is bounded uniformly in time (and in $N$).

	To bound the second moment of $G^i_{t}$, we have to bound each type of the following 
	expectations $\mathbb{E}[\tilde{H}(Z^{j_1,N}_{t}) \tilde{H}(Z^{j_2,N}_{t})]$, 
	$\mathbb{E}[\tilde{H}(Z^{j_1,N}_{t}) \tilde{H}(\bar{Z}^{j_2}_{t})]$, 
	$\mathbb{E}[\tilde{H}(\bar{Z}^{j_1}_{t}) \tilde{H}(\bar{Z}^{j_2}_{t})]$, 
	$\mathbb{E}[\tilde{H}{(Z^{j,N}_{t})}^2]$ and $\mathbb{E}[\tilde{H}{(\bar{Z}^{j}_{t})}^2]$. 
	By 
	Cauchy-Schwarz inequality, it is in fact enough to bound 
	$\mathbb{E}[\tilde{H}{(Z^{j,N}_{t})}^2]$ and $\mathbb{E}[\tilde{H}{(\bar{Z}^{j}_{t})}^2]$. 
	
	First, by the definition of $\tilde{H}$ in~\eqref{eq:deftildeH}, 
	\begin{align*}
		{\tilde{H}(z)}^2 =& 
		{\left(\frac{2}{a^2}\exp\left(a\sqrt{H(z)}\right)\left(a\sqrt{H(z)}-1\right)+\frac{2}{a^2}\right)}^2\\
		\leq& 2 \frac{2^2}{a^4}\exp\left(2a\sqrt{H(z)}\right){\left(a\sqrt{H(z)}-1\right)}^2+2 
		\frac{2^2}{a^4}\\
		\leq& \frac{8}{a^4} \exp\left(2a\sqrt{H(z)}\right) \left(2a^2 H(z) + 2\right) + \frac{8}{a^4}. 
	\end{align*}
	As for the first moment, the study of $Z^{j,N}_{t}$ is very similar to the one of 
	$\bar{Z}^j_{t}$. Here, we only focus on the second one. 
	
Using Cauchy-Schwarz inequality, bounds on $\mathbb{E}\left[H{(\bar{Z}^{j}_{t})}^2 \right]$ 
and $\mathbb{E}\left[ 
\exp\left(4a\sqrt{H(\bar{Z}^{j}_{t})}\right) \right]$	are sufficient  to bound 
$\mathbb{E}[\tilde{H}{(\bar{Z}^{j}_{t})}^2]$. The latter has already been bounded 
	uniformly in time, and the former can be obtained by the same calculations as previously, 
	replacing $a$ by $4a$ (and thus assuming $\tilde{a}\geq 4\sqrt{2}a\max(\sqrt{\gamma},1)$, 
	which we do).
	
	Finally, we deduce $\mathbb{E}\left({(G^i_{t})}^2\right)$ is bounded uniformly in time. 
\end{proof}

\begin{proof}[Proof of Lemma~\ref{lem:majorderivertildeH}]
	Using $\partial_{X}H(z)=\gamma x+\beta$, we have
	\begin{multline*}
		\left|\partial_{X}\tilde{H}{(Z^{i,N}_{t})}-\partial_{X}\tilde{H}(\bar{Z}^i_{t})\right|\\
		= \left|\left(\gamma
		 X^{i,N}_{t}+\beta \right)\exp\left(a\sqrt{H{(Z^{i,N}_{t})}}\right)-\left(\gamma 
		\bar{X}^i_{t}+\beta \right)\exp\left(a\sqrt{H(\bar{Z}^i_{t})}\right)\right|\\
		\leq\left|\gamma X^{i,N}_{t}-\gamma 
		\bar{X}^i_{t}\right|\left(\exp\left(a\sqrt{H{(Z^{i,N}_{t})}}\right)+\exp\left(a\sqrt{H(\bar{Z}^i_{t})}\right)\right)\\
		+\left|\gamma\bar{X}^i_{t}+\beta\right|\left|\exp\left(a\sqrt{H{(Z^{i,N}_{t})}}\right)-\exp\left(a\sqrt{H(\bar{Z}^i_{t})}\right)\right|.
	\end{multline*}
	Since $\left|X^{i,N}_{t}- \bar{X}^i_{t}\right| \leq r^i_{t}$,
	\begin{align*}
		\left|\gamma X^{i,N}_{t}-\gamma 
		\bar{X}^i_{t}\right|&\left(\exp\left(a\sqrt{H{(Z^{i,N}_{t})}}\right)+\exp\left(a\sqrt{H(\bar{Z}^i_{t})}\right)\right)
		 \\
		\leq& \gamma 
		r^i_{t}\left(\exp\left(a\sqrt{H{(Z^{i,N}_{t})}}\right)+\exp\left(a\sqrt{H(\bar{Z}^i_{t})}\right)\right)
				 	\end{align*}
	By Lemma~\ref{prop:H} (ii), we have $H(z) \geq \frac{1}{2} 
	\min\left(\frac{1}{\gamma},1\right){(\gamma x+\beta)}^2$. By the mean value theorem, for all 
	$y_1 \leq y_2$ in $\mathbb{R}$, there exists $y_3 \in [y_1, y_2]$ such that $e^{ay_1} - 
	e^{ay_2} = a (y_1-y_2) e^{ay_3}$. In particular, we have the following control ${|e^{ay_1} - 
	e^{ay_2}| \leq a |y_1-y_2| (e^{ay_1} + e^{ay_2})}$. Thus
	\begin{align*}
		\left|\gamma\bar{X}^i_{t}+\beta\right| 
		&\left|\exp\left(a\sqrt{H{(Z^{i,N}_{t})}}\right)-\exp\left(a\sqrt{H(\bar{Z}^i_{t})}\right)\right|\\
		\leq a & 
		\sqrt{\frac{2H(\bar{Z}^i_{t})}{\min\left(\frac{1}{\gamma},1\right)}}\left|\sqrt{H{(Z^{i,N}_{t})}}-\sqrt{H(\bar{Z}^i_{t})}\right|\left(\exp\left(a\sqrt{H{(Z^{i,N}_{t})}}\right)+\exp\left(a\sqrt{H(\bar{Z}^i_{t})}\right)\right)\\
		\leq a &\sqrt{2\max\left(\gamma, 
		1\right)}\left|H{(Z^{i,N}_{t})}-H(\bar{Z}^i_{t})\right|\left(\exp\left(a\sqrt{H{(Z^{i,N}_{t})}}\right)+\exp\left(a\sqrt{H(\bar{Z}^i_{t})}\right)\right).
	\end{align*}
	Then by the definition of $H$ we get
	\begin{align*}
		\left|H{(Z^{i,N}_{t})}\right.&\left.-H(\bar{Z}^i_{t})\right| \\
		=& \left| \frac{1}{2}\gamma \left({(X^{i,N}_{t})}^2 - {(\bar{X}^i_{t})}^2\right) +\beta 
		(X^{i,N}_{t} - 
		\bar{X}^i_{t})+\frac{1}{2}\left({(C^{i,N}_{t})}^2 - {(\bar{C}^i_{t})}^2\right) + \alpha 
		(C^{i,N}_{t} 
		-\bar{C}^i_{t}) \right|\\
		\leq&\frac{1}{2}\gamma\left|X^{i,N}_{t}-\bar{X}^i_{t}\right|\left|X^{i,N}_{t}+\bar{X}^i_{t}\right|+\beta\left|X^{i,N}_{t}-\bar{X}^i_{t}\right|+\frac{1}{2}\left|C^{i,N}_{t}-\bar{C}^i_{t}\right|\left|C^{i,N}_{t}+\bar{C}^i_{t}\right|\\
		&+\alpha\left|C^{i,N}_{t}-\bar{C}^i_{t}\right|.
	\end{align*}
	Now, by Lemma~\ref{prop:H} (i), we have $H(z) \geq \frac{\gamma}{4} x^2 + \frac{1}{4} 
	c^2$ and since $\left|X^{i,N}_{t}- \bar{X}^i_{t}\right| \leq r^i_{t}$ and $\left|C^{i,N}_{t}- 
	\bar{C}^i_{t}\right| \leq r^i_{t}/\delta$, we get
	\begin{align*}
		\left|X^{i,N}_{t}-\bar{X}^i_{t}\right|\left(\frac{1}{2}\gamma\left|X^{i,N}_{t}+\bar{X}^i_{t}\right|+\beta\right)
		 \leq r^i_{t} \left(\sqrt{\gamma} \left(\sqrt{H{(Z^{i,N}_{t})}}+\sqrt{H(\bar{Z}^i_{t})}\right) 
		+\beta\right)
	\end{align*}
	and
	\begin{align*}
		\left|C^{i,N}_{t}-\bar{C}^i_{t}\right|\left(\frac{1}{2}\left|C^{i,N}_{t}+\bar{C}^i_{t}\right|+\alpha\right)
		 \leq \frac{r^i_{t}}{\delta} \left(\sqrt{H{(Z^{i,N}_{t})}}+\sqrt{H(\bar{Z}^i_{t})}+\alpha\right).
	\end{align*}
	Thus
	\begin{align*}
		\left|H{(Z^{i,N}_{t})}-H(\bar{Z}^i_{t})\right| \leq 
		&\left(\beta+\frac{\alpha}{\delta}\right)r^i_{t}+\left(\sqrt{\gamma}+\frac{1}{\delta}\right)r^i_{t}\left(\sqrt{H{(Z^{i,N}_{t})}}+\sqrt{H(\bar{Z}^i_{t})}\right).
	\end{align*}
	Finally, 
	\begin{align*}
		|\partial_{X} & \tilde{H}{(Z^{i,N}_{t})}-\partial_{X}\tilde{H}(\bar{Z}^i_{t})|\\
		\leq &\gamma r^i_{t} 
		\left(\exp\left(a\sqrt{H{(Z^{i,N}_{t})}}\right)+\exp\left(a\sqrt{H(\bar{Z}^i_{t})}\right)\right)\\
		&+ a \sqrt{2\max\left(\gamma, 1\right)} \left(\beta+\frac{\alpha}{\delta}\right)r^i_{t} 
		\left(\exp\left(a\sqrt{H{(Z^{i,N}_{t})}}\right)+\exp\left(a\sqrt{H(\bar{Z}^i_{t})}\right)\right) \\
		&+a \sqrt{2\max\left(\gamma, 1\right)} 
		\left(\sqrt{\gamma}+\frac{1}{\delta}\right)r^i_{t}\left(\sqrt{H{(Z^{i,N}_{t})}} 
		+\sqrt{H(\bar{Z}^i_{t})}\right) \\
		& \hspace{4cm}\times  
		\left(\exp\left(a\sqrt{H{(Z^{i,N}_{t})}}\right)+\exp\left(a\sqrt{H(\bar{Z}^i_{t})}\right)\right)\\
		\leq & r^i_{t} \left( \gamma + a \sqrt{2\max\left(\gamma, 1\right)} 
		\left(\beta+\frac{\alpha}{\delta}\right) \right) 
		\left(\exp\left(a\sqrt{H{(Z^{i,N}_{t})}}\right)+\exp\left(a\sqrt{H(\bar{Z}^i_{t})}\right)\right)\\
		&+a r^i_{t} \sqrt{2\max\left(\gamma, 1\right)} \left(\sqrt{\gamma}+\frac{1}{\delta}\right)\\
		&\hspace{2 cm} \times\left(2 \sqrt{H{(Z^{i,N}_{t})}} \exp\left(a\sqrt{H{(Z^{i,N}_{t})}}\right) 
		+2 
		\sqrt{H(\bar{Z}^i_{t})} \exp\left(a\sqrt{H(\bar{Z}^i_{t})}\right)\right).
	\end{align*}
	Now, we can finally use Lemma~\ref{lem:control_{t}ilde_H}, and more 
	precisely~\eqref{eq:control_{t}ilde_H} and~\eqref{eq:control_{t}ilde_H_2}, we obtain
	\begin{align*}
		|\partial_{X} & \tilde{H}{(Z^{i,N}_{t})}-\partial_{X}\tilde{H}(\bar{Z}^i_{t})|\\
		\leq & r^i_{t} \left( \gamma + a \sqrt{2\max\left(\gamma, 1\right)} 
		\left(\beta+\frac{\alpha}{\delta}\right) \right) \left( \tilde{H}{(Z^{i,N}_{t})} + 
		\tilde{H}(\bar{Z}^i_{t}) + \frac{4}{a^2} \left(\e^{a^2/2} -1 \right) \right)\\
		&+a r^i_{t} \sqrt{2\max\left(\gamma, 1\right)} 
		\left(\sqrt{\gamma}+\frac{1}{\delta}\right)\left(2 a \tilde{H}{(Z^{i,N}_{t})} + \frac{2}{a} 
		(\e-2) 
		+ 2 a \tilde{H}(\bar{Z}^i_{t}) + \frac{2}{a} (\e-2) \right)\\
		\leq & r^i_{t} \left( \tilde{H}{(Z^{i,N}_{t})} + \tilde{H}(\bar{Z}^i_{t}) \right) \left[ \gamma + a 
		\sqrt{2\max\left(\gamma, 1\right)} \left(\beta+\frac{\alpha}{\delta}\right) + 2 a^2 
		\sqrt{2\max\left(\gamma, 1\right)} \left(\sqrt{\gamma}+\frac{1}{\delta}\right) \right]\\
		&+ r^i_{t} \left[ \left( \gamma + a \sqrt{2\max\left(\gamma, 1\right)} 
		\left(\beta+\frac{\alpha}{\delta}\right) \right) \frac{4}{a^2} \left(\e^{a^2/2} -1 \right) \right.\\
		&\quad\quad\quad+\left. 4 
		\sqrt{2\max\left(\gamma, 1\right)} \left(\sqrt{\gamma}+\frac{1}{\delta}\right) (\e-2) \right].
	\end{align*}
	We denote by $\mathcal{C}_{f,1}$ and $\mathcal{C}_{f,2}$ (given in Lemma~\ref{lem:parameters}) the following constants
	\begin{align*}
		\mathcal{C}_{f,1} &= 4 \left[\left( \gamma + a \sqrt{2\max\left(\gamma, 1\right)} 
		\left(\beta+\frac{\alpha}{\delta}\right) \right) \frac{4}{a^2} \left(\e^{a^2/2} -1 \right) \right.\\
&\quad 	\quad \quad 	
		+ \left. 4 \sqrt{2\max\left(\gamma, 1\right)} \left(\sqrt{\gamma}+\frac{1}{\delta}\right) (\e-2) 
		\right]\\
		\mathcal{C}_{f,2} &= 4 \left[\gamma + a \sqrt{2\max\left(\gamma, 1\right)} \left(\beta+\frac{\alpha}{\delta}\right) + 2 a^2 \sqrt{2\max\left(\gamma, 1\right)} \left(\sqrt{\gamma}+\frac{1}{\delta}\right)\right]
	\end{align*}
	By the definition of $G^i_{t}$ and since $G^i_{t} \geq 1$, we obtain
	\begin{align*}
		|\partial_{X} & \tilde{H}{(Z^{i,N}_{t})}-\partial_{X}\tilde{H}(\bar{Z}^i_{t})| 
		\leq r^i_{t} \frac{G^i_{t}}{\epsilon} \frac{\mathcal{C}_{f,2}}{4} + r^i_{t} G^i_{t} 
		\frac{\mathcal{C}_{f,1}}{4}, 
	\end{align*}
	and eventually
	\begin{align*}
		2\epsilon\left(1+\frac{1}{N}\right)\sigma_{X}^2\varphi_{\text{rc}}{\left(|X^{i,N}_{t}-\bar{X}^i_{t}|\right)}^2&\left|\partial_{X}\tilde{H}{(Z^{i,N}_{t})}-\partial_{X}\tilde{H}(\bar{Z}^i_{t})\right|\\
		&\leq \left(\epsilon 
		\mathcal{C}_{f,1}+\mathcal{C}_{f,2}\right)\sigma_{X}^2\varphi_{\text{rc}}{\left(|X^{i,N}_{t}-\bar{X}^i_{t}|\right)}^2r^i_{t}
		 G^i_{t}.
	\end{align*}
	
\end{proof}

%
%
%Section
%
%

\section{Proof of Theorem~\ref{thm:unif} in the case $\sigma_{X}=0$ and 
$\sigma_{C}>0$}\label{app:sigma_{X}_0}

We quickly explain in this section how we may also deal with the case $\sigma_{X}=0$ and 
$\sigma_{C}>0$. Recall how the choice of the coupling method was motivated by the 
observation in~\eqref{eq:dC} that the difference of potentials 
$\left|C^{i,N}_{t}-\bar{C}^i_{t}\right|$ was naturally contracting when 
$\left|X^{i,N}_{t}-\bar{X}^i_{t}\right|$ was close to $0$. This lead us to use a reflection 
coupling 
on the Brownian motions acting on the potential $X$, to bring the difference close to $0$, and 
it was thus necessary for $\sigma_{X}$ to be positive ($\sigma_{C}$ however did not matter). 
In 
the case $\sigma_{X}=0$, we then have to assume $\sigma_{C}>0$, and we do a change of 
variable, motivated by the following observation. We have, when $\sigma_{X}=0$
\begin{align*}
	d(X^{i,N}_{t}-\bar{X}^i_{t})=&\left((X^{i,N}_{t}-\bar{X}^i_{t})-({(X^{i,N}_{t})}^3-{(\bar{X}^i_{t})}^3)-(C^{i,N}_{t}-\bar{C}^i_{t})\right)dt\\
	&+\left(\frac{1}{N}\sum_{j=1}^{N}K_{X}(Z^{i,N}_{t}-Z^{j,N}_{t})-K_{X}\ast\bar{\rho}(\bar{Z}^i_{t})\right)dt\\
	=&\left(2(X^{i,N}_{t}-\bar{X}^i_{t})-(C^{i,N}_{t}-\bar{C}^i_{t})-(X^{i,N}_{t}-\bar{X}^i_{t})-({(X^{i,N}_{t})}^3-{(\bar{X}^i_{t})}^3)\right)dt\\
	&+\left(\frac{1}{N}\sum_{j=1}^{N}K_{X}(Z^{i,N}_{t}-Z^{j,N}_{t})-K_{X}\ast\bar{\rho}(\bar{Z}^i_{t})\right)dt.
\end{align*}
Thus
\begin{align*}
	d|X^{i,N}_{t}-\bar{X}^i_{t}|=&\text{sign}(X^{i,N}_{t}-\bar{X}^i_{t})\left(2(X^{i,N}_{t}-\bar{X}^i_{t})-(C^{i,N}_{t}-\bar{C}^i_{t})\right)dt\\
	&-\left(|{(X^{i,N}_{t})}^3-{(\bar{X}^i_{t})}^3|-|X^{i,N}_{t}-\bar{X}^i_{t}|\right)dt\\
	&+\text{sign}(X^{i,N}_{t}-\bar{X}^i_{t})\left(\frac{1}{N}\sum_{j=1}^{N}K_{X}(Z^{i,N}_{t}-Z^{j,N}_{t})-K_{X}\ast\bar{\rho}(\bar{Z}^i_{t})\right)dt.
\end{align*}
The quantity $|X^{i,N}_{t}-\bar{X}^i_{t}|$ is therefore naturally contracting when 
$|2(X^{i,N}_{t}-\bar{X}^i_{t})-(C^{i,N}_{t}-\bar{C}^i_{t})|$ is close to $0$. Thanks to the 
presence of a Brownian motion in the stochastic differential equations defining the potential 
$C$, we can now use a reflection coupling to have 
$|2(X^{i,N}_{t}-\bar{X}^i_{t})-(C^{i,N}_{t}-\bar{C}^i_{t})|$ go to $0$. Consider the following 
coupling
\begin{equation}\label{eq:FN_MF_enC}
	\left\{
	\begin{array}{ll}
		dX^{i,N}_{t}=(X^{i,N}_{t}-{(X^{i,N}_{t})}^3-C^{i,N}_{t}-\alpha)dt+\dfrac{1}{N}\sum_{j=1}^{N}K_{X}(Z^{i,N}_{t}-Z^{j,N}_{t})dt\\
		dC^{i,N}_{t}=(\gamma 
		X^{i,N}_{t}-C^{i,N}_{t}+\beta)dt+\dfrac{1}{N}\sum_{j=1}^{N}K_{C}(Z^{i,N}_{t}-Z^{j,N}_{t})dt\\
		\hspace{4 cm}+ \sigma_{C} 
		\varphi_{\text{sc}}\left(|2(X^{i,N}_{t}-\bar{X}^i_{t})-(C^{i,N}_{t}-\bar{C}^i_{t})|\right)dB^{i,sc,C}_{t}\\
		\hspace{4 cm} +\sigma_{C} 
		\varphi_{\text{rc}}\left(|2(X^{i,N}_{t}-\bar{X}^i_{t})-(C^{i,N}_{t}-\bar{C}^i_{t})|\right)dB^{i,rc,C}_{t}
		 ,
	\end{array}
	\right.
\end{equation}
and
\begin{equation}\label{eq:FN_limit_enC}
	\left\{
	\begin{array}{ll}
		d\bar{X}^i_{t}=(\bar{X}^i_{t}-{(\bar{X}^i_{t})}^3-\bar{C}^i_{t}-\alpha)dt+K_{X}\ast\bar{\rho}(\bar{Z}^i_{t})dt\\
		d\bar{C}^i_{t}=(\gamma 
		\bar{X}^i_{t}-\bar{C}^i_{t}+\beta)dt+K_{C}\ast\bar{\rho}(\bar{Z}^i_{t})dt 
		\\
		\hspace{4 cm}+ \sigma_{C} 
		\varphi_{\text{sc}}\left(|2(X^{i,N}_{t}-\bar{X}^i_{t})-(C^{i,N}_{t}-\bar{C}^i_{t})|\right)dB^{i,sc,C}_{t}\\
		\hspace{4 cm} -\sigma_{C} 
		\varphi_{\text{rc}}\left(|2(X^{i,N}_{t}-\bar{X}^i_{t})-(C^{i,N}_{t}-\bar{C}^i_{t})|\right)dB^{i,rc,C}_{t}\\
		\bar{\rho}=\mathcal{L}((\bar{X}^1_{t},\bar{C}^1_{t})),
	\end{array}
	\right.
\end{equation}
and for $\delta >0$, the following modified distance 
\[
r^i_{t}=\delta 
|X^{i,N}_{t}-\bar{X}^i_{t}|+|2 (X^{i,N}_{t}-\bar{X}^i_{t}) - (C^{i,N}_{t}-\bar{C}^i_{t})|.
\]
Like 
previously, we consider a modified semi-metric of the form $\frac{1}{N}\sum f(r^i_{t})G^i_{t}$ 
and similar calculations yield
\begin{align*}
	d(e^{ct}f(r^i_{t})G^i_{t})\leq e^{ct}K^i_{t}dt+dM^i_{t},
\end{align*}
where $M^i_{t}$ is a continuous local martingale and
\begin{align*}
	K^i_{t} = \tilde{K}^i_{t} + I^{1,i}_{t} + I^{2,i}_{t} + I^{3,i}_{t}.
\end{align*}
We define $\tilde{K}^i_{t}$, $I^{1,i}_{t}$, $I^{2,i}_{t}$ and $I^{3,i}_{t}$ as follows
\begin{align*}
	\tilde{K}^i_{t} =& G^i_{t} \left[2 c f(r^i_{t}) + 2 f''(r^i_{t})\sigma_{C}^2 
	\varphi_{\text{rc}}{\left(|2(X^{i,N}_{t}-\bar{X}^i_{t})-(C^{i,N}_{t}-\bar{C}^i_{t})|\right)}^2 
	\right.\\
	&\left. + f'(r^i_{t}) \left(\left|2(X^{i,N}_{t}-\bar{X}^i_{t})-(C^{i,N}_{t}-\bar{C}^i_{t})\right| 
	(\delta + 1 ) 
	- |{(X^{i,N}_{t})}^3-{(\bar{X}^i_{t})}^3| ( \delta - 2 ) \right.\right. \\
	& \left.\left. + |X^{i,N}_{t}-\bar{X}^i_{t}| \left( -\delta + \gamma + L_{X}(\delta +2) + L_{C} 
	\right) 
	+ |C^{i,N}_{t}-\bar{C}^i_{t}| (L_{X}(\delta + 2) + L_{C}) \right.\right.\\
	& \left.\left. + \sigma_{C}^2 
	\varphi_{\text{rc}}{\left(|2(X^{i,N}_{t}-\bar{X}^i_{t})-(C^{i,N}_{t}-\bar{C}^i_{t})|\right)}^2\left(\epsilon\mathcal{C}_{f,1}+\mathcal{C}_{f,2}\right)
	 r^i_{t}\right) \right] \\
	&+\epsilon f(r^i_{t})\left(4\tilde{B}-\frac{\lambda}{8} 
	\tilde{H}(\bar{Z}^i_{t})-\frac{\lambda}{8} 
	\tilde{H}{(Z^{i,N}_{t})}-\frac{\lambda}{8N}\sum_{j=1}^{N} 
	\tilde{H}(\bar{Z}^j_{t})-\frac{\lambda}{32N}\sum_{j=1}^{N} \tilde{H}(Z^{j,N}_{t})\right), \\
	I^{1,i}_{t} =& G^i_{t}f'(r^i_{t})\left[(\delta + 2) 
	\left(\left|\frac{1}{N}\sum_{j=1}^{N}K_{X}(\bar{Z}^i_{t}-\bar{Z}^j_{t})-K_{X}\ast\bar{\mu}_{t}(\bar{Z}^i_{t})\right|\right)\right.\\
	&\hspace{5cm}\left.
	+ 
	\left(\left|\frac{1}{N}\sum_{j=1}^{N}K_{C}(\bar{Z}^i_{t}-\bar{Z}^j_{t})-K_{C}\ast\bar{\mu}_{t}(\bar{Z}^i_{t})\right|\right)
	 \right],\\
	I^{2,i}_{t} = &G^i_{t}f'(r^i_{t})\left[(\delta + 
	2)\left(\frac{L_{X}}{N}\left(\sum_{j=1}^{N}|X^{j,N}_{t}-\bar{X}^j_{t}|+|C^{j,N}_{t}-\bar{C}^j_{t}|\right)\right)\right.\\
	&\hspace{5cm}\left.+ 
	\left(\frac{L_{C}}{N}\left(\sum_{j=1}^{N}|X^{j,N}_{t}-\bar{X}^j_{t}|+|C^{j,N}_{t}-\bar{C}^j_{t}|\right)\right)
	 \right] \nonumber\\
	&-cf(r^i_{t})G^i_{t}-\epsilon 
	f(r^i_{t})\left[\frac{\lambda}{32}H(\bar{Z}^i_{t})\exp\left(a\sqrt{H(\bar{Z}^i_{t})}\right)+\frac{\lambda}{32}H{(Z^{i,N}_{t})}\exp\left(a\sqrt{H{(Z^{i,N}_{t})}}\right)\right]
	 \nonumber\\
	&-\epsilon f(r^i_{t})\left[ 
	\frac{\lambda}{32N}\sum_{j=1}^{N}H(\bar{Z}^j_{t})\exp\left(a\sqrt{H(\bar{Z}^j_{t})}\right) 
	\right.\\
	&\quad\quad 
	+\left.\frac{\lambda}{32N}\sum_{j=1}^{N}H(Z^{j,N}_{t})\exp\left(a\sqrt{H(Z^{j,N}_{t})}\right)\right],
	\end{align*}
	\begin{align*}
	I^{3,i}_{t} =& \epsilon 
	f(r^i_{t})\left(\left(\alpha_{X}L_{X}+\beta_{X}L_{C}\right){\left(\frac{\sum_{j=1}^{N}|X^{j,N}_{t}|}{N}\right)}^2\exp\left(a\sqrt{H{(Z^{i,N}_{t})}}\right)\right.
	 \nonumber\\
	&\hspace{1.3cm}\left. +\left(\alpha_{C}L_{X}+ 
	\beta_{C}L_{C}\right){\left(\frac{\sum_{j=1}^{N}|C^{j,N}_{t}|}{N}\right)}^2\exp\left(a\sqrt{H{(Z^{i,N}_{t})}}\right)\right.
	 \nonumber\\
	&\hspace{1.3cm}\left. 
	-\frac{\lambda}{16}H{(Z^{i,N}_{t})}\exp\left(a\sqrt{H{(Z^{i,N}_{t})}}\right)-\frac{\lambda}{16N}\sum_{j=1}^{N}H(Z^{j,N}_{t})\exp\left(a\sqrt{H(Z^{j,N}_{t})}\right)\right).
\end{align*}
We then have the additional constraint of $\delta>2$ (so that the coefficient appearing in front 
of $|{(X^{i,N}_{t})}^3-{(\bar{X}^i_{t})}^3|$ in the expression of $\tilde{K}^i_{t}$ is 
non-positive). 
Otherwise, we deal with the various terms exactly as previously, through the choice of a 
sufficiently concave function $f$ and a law of large numbers, and by considering the regions 
of space
\begin{align*}
	\mbox{Reg}^i_1=&\left\{(\bar{Z}^i_{t},Z^{i,N}_{t})\text{ s.t. 
	}|2(X^{i,N}_{t}-\bar{X}^i_{t})-(C^{i,N}_{t}-\bar{C}^i_{t})|\geq\xi \text{ and } r^i_{t}\leq 
	R\right\},\\
	\mbox{Reg}^i_2=&\left\{(\bar{Z}^i_{t},Z^{i,N}_{t})\text{ s.t. 
	}|2(X^{i,N}_{t}-\bar{X}^i_{t})-(C^{i,N}_{t}-\bar{C}^i_{t})|<\xi \text{ and } r^i_{t}\leq 
	R_1\right\},\\
	\mbox{Reg}^i_3=&\left\{(\bar{Z}^i_{t},Z^{i,N}_{t})\text{ s.t. }r^i_{t}> 
	R\right\}.
\end{align*}

%
%
%Section
%
%
%
%\section*{Acknowledgements}
%
%Laetitia Colombani is a PhD student under the supervision of Patrick Cattiaux and Manon Costa, and Pierre {Le Bris} is a PhD student under the supervision of Arnaud Guillin and Pierre Monmarché. The authors would like to thank them, as well as Samir Salem, for their help throughout the redaction of the present article.
%
%This work has been (partially) supported by the Project EFI ANR-17-CE40-0030 of the French National Research Agency.

%
%
%Section
%
%

\section*{Index}

Throughout this article, we define many parameters and constants. For the sake of clarity, we list the main ones here so as to give the reader an index to refer to.
\begin{itemize}
	\item $\mathbf{X,C,Z}$: $\hspace{0.5cm}$  $X$ and $C$ are the processes we consider 
	(see~\eqref{eq:FN_MF} and~\eqref{eq:FN_limit}) and we often refer to $Z=(X,C)$,
	\item $\mathbf{\bar{\mu}_{t}=\text{Law}{(\bar{Z}_{t})}}$:$\hspace{0.5cm}$ the density of 
	the 
	non-linear 
	limit (see~\eqref{eq:FN_limit}),
	\item $\mathbf{\alpha, \beta, \gamma, \sigma_{X}, \sigma_{C}}$:$\hspace{0.5cm}$ 
	parameters of the problem (see~\eqref{eq:FN_MF}),
	\item $\mathbf{K_{X}, K_{C}, L_{X}, L_{C}, L_{X,\max}, L_{C,\max}}$:$\hspace{0.5cm}$ 
	$K_{X}$ 
	(resp. 
	$K_{C}$) is a Lipschitz continuous interaction kernel, with Lipschitz constant $L_{X}\in[0, 
	L_{X,\max}]$ (resp. $L_{C}\in[0, L_{C,\max}]$), as given in Assumption~\ref{hyp:K}. In the 
	case of uniform in time propagation of chaos, the inequalities $L_{X}$ and $L_{C}$ must 
	satisfy are listed in Subsection~\ref{subsec:hyp_{C}onstants},
	\item $\mathbf{\mathcal{W}_p}$:$\hspace{0.5cm}$ the usual Wasserstein distance 
	associated to the 
	$L^p$ distance (see~\eqref{eq:def_W}),
	\item $\mathbf{a, \tilde{a}, \mathcal{C}_{init,exp}}$:$\hspace{0.5cm}$ constants used to 
	give an exponential initial moment to the problem (see the assumptions of 
	Theorem~\ref{thm:unif} and Section~\ref{subsec:tildeH}),
	\item $\mathbf{\lambda, B, \tilde{B}, H, \tilde{H}, \alpha_{X}, \alpha_{C}, \beta_{X}, 
	\beta_{C}}$:$\hspace{0.5cm}$ $H$ (resp. $\tilde{H}$) is a Lyapunov functions given 
	in~\eqref{eq:def_H} (resp.~\eqref{eq:deftildeH}). Its main property involves parameters 
	$\lambda$ and $B$ (resp. $\lambda$ and $\tilde{B}$), as can for instance be seen 
	in~\eqref{eq:dyn_H} (resp.~\eqref{eq:dyn_{t}ilde_H_non_lin}). $\alpha_{X}, \alpha_{C}, 
	\beta_{X}$ 
	and $\beta_{C}$ are intermediate constants given in Lemma~\ref{lem:Lya_limit},
	\item $\mathbf{c}$:$\hspace{0.5cm}$ a contraction rate (see 
	Subsection~\ref{subsec:hyp_{C}onstants}),
	\item $\mathbf{r, f,  g, \phi, \Phi, G, \rho, \delta, R, \epsilon, \mathcal{C}_{f,1}, 
	\mathcal{C}_{f,2}}$:$\hspace{0.5cm}$ $f$ (see~\eqref{eq:def_f}) is a concave function, the 
	definition of which involves $g$, $\phi$, $\Phi$ (see 
	Subsection~\ref{subsec:hyp_{C}onstants}). Function $G$ (see~\eqref{eq:def_G}) is then 
	used to define $\rho$ (see~\eqref{eq:def_rho}), the semi-metric we consider in the end. All 
	these notations thus refer to the modified distance we consider. These functions will be 
	applied to a modification $r$ of the usual $L^1$ distance (see equation~\eqref{eq:def_r}). 
	Then, parameters $\delta$, $R$, $\epsilon$, $\mathcal{C}_{f,1}$, and $\mathcal{C}_{f,2}$ 
	are used to define such functions (see Subsection~\ref{subsec:hyp_{C}onstants} for some 
	explicit values),
	\item $\mathbf{R_0, \phi_{\min}}$:$\hspace{0.5cm}$ intermediate constants (see 
	Subsection~\ref{subsec:hyp_{C}onstants}),
	\item $\mathcal{C}_{init,2}$:$\hspace{0.5cm}$ uniform in time bound on the second moment of the processes (see Lemma~\ref{lem:borne_unif_moment_2}),
	\item $\mathbf{\mathcal{C}_1, \mathcal{C}_2, \mathcal{C}_{z}}$:$\hspace{0.5cm}$ 
	constants 
	used to quantify the control our modified distance has over the usual $L^1$ and $L^2$ 
	distance (see Lemma~\ref{lem:rho_1_2} for the control and 
	Subsection~\ref{subsec:hyp_{C}onstants} for explicit values),
	\item $\mathbf{\phi_{rc}, \phi_{sc}, \xi}$:$\hspace{0.5cm}$ $\phi_{rc}$ and $\phi_{sc}$ are 
	two 
	Lipschitz continuous functions used to define the coupling method, and their definitions 
	involve a parameter $\xi$ which converges to 0 in the end (see the beginning of 
	Section~\ref{sec:proof_{t}hm_unif}),
	\item $\mathbf{\mathcal{C}_{r,H}}$:$\hspace{0.5cm}$ used to explicit the control of the 
	Lyapunov function $H$ over the distance $r$ (see Lemma~\ref{prop:H}), 
\end{itemize}

%%%%%%%%%%%%%%%%%%%%%%%%%%%%%%%%%%%%%%%%%%%%%%%%%%%%%%%%%%%%%%%%%%%
%%                                                               %%
%% Supplementary Material, if any, should be provided in         %%
%% {supplement} environment  with title and short description.   %%
%%                                                               %%
%%%%%%%%%%%%%%%%%%%%%%%%%%%%%%%%%%%%%%%%%%%%%%%%%%%%%%%%%%%%%%%%%%%
%
%\begin{supplement}
%\stitle{Title of Supplement A.}
%\sdescription{Short description of Supplement A.}
%\end{supplement}
%\begin{supplement}
%\stitle{Title of Supplement B.}
%\sdescription{Short description of Supplement B.}
%\end{supplement}

%%%%%%%%%%%%%%%%%%%%%%%%%%%%%%%%%%%%%%%%%%%%%%%%%%%%%%%%%%%%%%%%%%%
%%                                                               %%
%% Use the two commands below for producing your bibliography    %%
%% with bibtex, then comment again the commands and include the  %%
%% content of the .bbl file in this file below the commands.     %%
%%                                                               %%
%%%%%%%%%%%%%%%%%%%%%%%%%%%%%%%%%%%%%%%%%%%%%%%%%%%%%%%%%%%%%%%%%%%

%\bibliographystyle{amsplainhyper}  
%\bibliography{biblio}

% add below the content of your .bbl file produced by bibtex.
%
\providecommand{\bysame}{\leavevmode\hbox to3em{\hrulefill}\thinspace}
\providecommand{\MR}{\relax\ifhmode\unskip\space\fi MR }
% \MRhref is called by the amsart/book/proc definition of \MR.
\providecommand{\MRhref}[2]{%
  \href{http://www.ams.org/mathscinet-getitem?mr=#1}{#2}
}
\providecommand{\bysame}{\leavevmode\hbox to3em{\hrulefill}\thinspace}

%%%%%%%%%%%%%%%%%%%%%%%%%%%%%%%%%%%%%%%%%%%%%%%%%%%%%%%%%%%%%%%%%%%
%%                                                               %%
%% You may add acknowledgments (optional).                       %%
%%                                                               %%
%%%%%%%%%%%%%%%%%%%%%%%%%%%%%%%%%%%%%%%%%%%%%%%%%%%%%%%%%%%%%%%%%%%

\ACKNO{During this study, Laetitia Colombani was a PhD student under the supervision of 
Patrick Cattiaux and Manon Costa, and Pierre {Le Bris} a PhD student under the supervision of 
Arnaud Guillin and Pierre Monmarché. The authors would like to thank them, as well as Samir 
Salem, for their help throughout the redaction of the present article.
	This work has been (partially) supported by the Project EFI ANR-17-CE40-0030 of the French National Research Agency.}

%%%%%%%%%%%%%%%%%%%%%%%%%%%%%%%%%%%%%%%%%%%%%%%%%%%%%%%%%%%%%%%%%%%
%%                                                               %%
%% You have reached the end of your document.                    %%
%%                                                               %%
%%%%%%%%%%%%%%%%%%%%%%%%%%%%%%%%%%%%%%%%%%%%%%%%%%%%%%%%%%%%%%%%%%%
\nocite{*}
\end{document}